\documentclass[a4paper,12pt]{article}
\usepackage{amsmath}
\usepackage{amssymb}
\usepackage{tabularx}
\usepackage{enumerate}
\usepackage{graphicx}
\usepackage{subfigure}
\usepackage{color}
\usepackage[pagewise]{lineno}
\usepackage{longtable}
\usepackage{epstopdf}
\usepackage{float}

\newtheorem{thm}{Theorem}[section]
\newtheorem{lem}[thm]{Lemma}

\newtheorem{exm}{Example}[section]

\newtheorem{rmk}{Remark}[section]

\newtheorem{pppp}{Proof.}

\newenvironment{pf}{\begin{pppp} \em}{\mbox{}\hfill\qed\end{pppp}}
\newcommand{\qed}{\hspace{1em}\mbox{\raisebox{0.65ex}{\fbox{}}}}

\numberwithin{equation}{section}

\newcommand{\be}{\begin{equation}}
\newcommand{\ee}{\end{equation}}
\newcommand\bes{\begin{eqnarray}} \newcommand\ees{\end{eqnarray}}
\newcommand{\bess}{\begin{eqnarray*}}
\newcommand{\eess}{\end{eqnarray*}}

\begin{document}

\thispagestyle{empty}

\title{A competition model with impulsive interventions and environmental perturbations  in moving environments\thanks{The work is partially supported by the NNSF of China (Grant No. 12271470, 61877052), Postgraduate Research \& Practice Innovation Program
of Jiangsu Province (KYCX22\_3446) and the support of CNPq/Brazil Proc.  $N^o$ $311562/2020-5$ and FAPDF grant 00193.00001133/2021-80.}}

\date{\empty}

\author{Yue Meng$^1$,  Zhigui Lin$^1\thanks{Corresponding author. Email: zglin@yzu.edu.cn (Z. Lin).}$ and Carlos Alberto Santos$^2$\\
{\footnotesize  1 School of Mathematical Science, Yangzhou University, Yangzhou 225002, China}\\
{\footnotesize 2 Department of Mathematics, University of Brasilia, BR-70910900 Brasilia, DF, Brazil}}

 \maketitle
\begin{quote}
\noindent
{\bf Abstract.} { 
In order to understand how impulsive interventions and environmental perturbations affect dynamics of competitors,  we focus on a diffusive competition model with  free boundaries and periodic pulses in a temporally heterogeneous environment with upward or downward advection. The dependence of the principal eigenvalue of corresponding periodic impulsive eigenvalue problem on advection rates, habitat sizes and pulses is investigated, which gives precise conditions that classify the dynamics into four types of competition outcomes including coexistence, co-extinction, two different competition exclusions for small or negative advection rates.  Some sufficient conditions on pulses or initial habitats for species spreading or vanishing,  and spreading speeds are then established. Our results  not only extend  the existing ones to the case with pulses, but also reveal the effects of human and natural factors, that is, impulsive interventions factors including positive or negative impulsive effect, pulse intensity and timing can significantly affect and alter the competition outcomes. The different performances of the superior and inferior affected by environmental perturbations are also reflected in the simulations.
}


\end{quote}

\section{Introduction}
 Some species are influenced by environmental factors or their own perception to migrate to favorable habitats, and their movement shows a directional trend, which can be described as an advection term in  mathematical models. Also, the habitats of species change naturally and are affected by the movement of species itself.  The unknown expansion frontiers of species can be modelled by free boundary problems.  Taking  advection and free boundary into consideration for one single species problem, Gu et al.  studied the problem
 \begin{eqnarray}
\left\{
\begin{array}{lll}
u_t=u_{xx}-q u_x+u(1-u),\; &\ t>0,\ x\in(g(t),h(t)), \\[2mm]
u(t,g(t))=0, g'(t)=-\mu u_x(t,g(t)),\; &\ t>0,\\[2mm]
u(t,h(t))=0, h'(t)=-\mu u_x(t,h(t)),\; &\ t>0,\\[2mm]
g(0)=-h_0,h(0)=h_0, \\[2mm]
u(0,x)=u_0(x),\; &\ x\in[-h_0,h_0]
\end{array} \right.
\end{eqnarray}
in \cite{gulin2014} and \cite{gulin2015}, where $qu_x$  is an advection term, and other meanings of parameters are omitted here.  The detailed  derivation and explanations of Stefan condition modeling the unknown fronts can be referred to \cite{lzg2007}. Gu et al.  obtained a spreading-vanishing dichotomy and estimates of spreading speed for small advection rates. The results were later extended for  advection rate $q\in(0,\infty)$ in \cite{ghjfa2015}. They found two critical values $c_0$ and $q^*$ of advection rates to give more complete dynamics of solution, that is, $(i)$ when advection rate $q$ is small (i.e. $q\in(0,c_0)$ ), a dichotomy result occurs, $(ii)$ when advection rate $q$ is moderate (i.e. $q\in[c_0,q^*)$ ), a trichotomy result exists,  and $(iii)$ when advection rate $q$ is large (i.e. $q\in[q^*,\infty)$ ), vanishing occurs, where $c_0$ is the minimal speed of traveling waves and more details about $c_0$, $q^*$ and results can be referred to Theorems 2.1-2.4 in \cite{ghjfa2015}. The negative advection case was later focused on in \cite{zhaowang2018}. The interested readers can also refer to  \cite{kaneko,liangliujin,monobe} and related references for more results about the effects of advection in free boundary problems modeling one single species.

Competition is a common relationship between two species. The Lotka-Volterra  competition model with advection
\begin{eqnarray}
\left\{
\begin{array}{lll}
u_t=d_1u_{xx}-\alpha_1 u_x+u(a_1-b_1u-c_1v), \\[2mm]
v_t=d_2v_{xx}-\beta_1 v_x+v(a_2-b_2v-c_2u)
\end{array} \right.
\label{i01}
\end{eqnarray}
modeling the dynamics of two competing species is well-known, where all coefficients in \eqref{i01} are positive, and $d_i$ and $a_i$ $(i=1,2)$ stand for the diffusion rates and intrinsic growth rates of species $u$ and $v$, $b_i$ and $c_i$ $(i=1,2)$ represent the intra-specific and inter-specific competition rates, respectively, and  $\alpha_1$ and $\beta_1$ are advection rates.  There has been increasingly keen interest in diffusive-advective competition problem \eqref{i01} (and its variations) with free boundaries, see \cite{duanjmaa, duandcds, xuren} for the case that species have the same  expanding fronts, \cite{ren} for the case that both species have the same advection rate and different free boundaries, and references therein.
We also remark some other interesting works about competition free boundary problems without advection here.  A free boundary competition problem that native species distributes in the entire space while invasive species invades with an unknown moving front was investigated by \cite{dulin2014, duwangzhou, wangj2015, wangzhang}.  Nonlocal diffusion has been recently introduced in \cite{caoade,liwang2020,PLL, zhangliwang,zhangzhou2022} and corresponding dynamics have been investigated.

The analysis and results on competition free boundary problems in homogeneous environments have been fully established and developed.
In reality, seasonal variations and climate changes make these environments heterogeneous and change spatially and temporally, which indicates that the parameters in \eqref{i01} depends on time $t$ and location $x$.  Several advection competition problems with free boundaries in heterogeneous environments have been investigated in \cite{duan2021, liuzh, zhouling2017} with intrinsic growth rates depends only on space, and in \cite{chenql} with all parameters are dependent in time and space.
Very recently, Khan et al. \cite{khan} investigated competition free boundary problems without advection in space-periodic or time-periodic environments via numerical simulations. Their observations revealed the effects of environmental heterogeneity and also showed that the heterogeneity of environment makes competition outcomes less predictable.
This leads us to wonder whether the desired competition outcomes or even the coexistence of two species can be achieved by virtue of human interventions.

Human interventions including harvesting, pesticide applications and natural enemies release are widely used in agriculture, fishing, pest management etc. (\cite{biological control}), and are usually implemented at a specific time. This  discrete-continuous dynamics of species can be described by impulsive reaction diffusion equations. For a  single species diffusive model, we refer to \cite{lewisli} by Lewis and Li with a seasonal pulse, \cite{fazly20} with a nonlocal discrete-time map, \cite{liangtangsy} with multiple periodic pulses, \cite{wuzhao19,wuzhao22} with nonlocal  dispersal, \cite{menglin21} with impulsive harvesting in a periodically evolving domain, \cite{menglin22} with a periodic pulse and  free boundary, etc. Taking account of advective environments, a diffusion-advection competition model with two different impulsive interventions have been recently studied in \cite{menglin23} and have obtained that there are a tradeoff between movements and impulsive interventions
to pursue the desired competition outcomes.

The paper is aim to investigate the combined effects of human and natural factors on the spreading or vanishing of competitors. For simplicity, let
\[
\hat{u}(t,x)=\frac{b_1}{a_1}u(\frac{t}{a_2},\sqrt{\frac{d_2}{a_2}}x),\ \hat{v}(t,x)=\frac{b_2}{a_2}v(\frac{t}{a_2},\sqrt{\frac{d_2}{a_2}}x),
\]
 denote
\[
D=\frac{d_1}{d_2},\gamma=\frac{a_1}{a_2},\alpha=\frac{\alpha_1}{\sqrt{a_2d_2}},
\beta=\frac{\beta_1}{\sqrt{a_2d_2}}, k=\frac{c_1a_2}{b_2a_1}, h=\frac{c_2a_1}{b_1a_2},
\]
and drop the hat sign,
then problem \eqref{i01} is transformed to the following form that  will be studied in this paper:
\begin{eqnarray}
\left\{
\begin{array}{lll}
u_t=Du_{xx}-\alpha u_x+\gamma u(1-u-kv), \\[2mm]
v_t=v_{xx}-\beta v_x+v(1-v-hu).
\end{array} \right.
\label{i02}
\end{eqnarray}
A natural phenomenon that there are two competing species initially occupying  different regions  might happen. We assume that the two competitors obey the competition relationships shown in \eqref{i02} and move along different  unknown fronts modeled by Stefan type free boundaries. And we only consider one-sided free boundaries in a one dimensional space as \cite{dyhlzg}, then, the left boundaries are fixed and satisfy Neumann boundary conditions at $x=0$, that is, species can only invade from the right end into the environments. The seasonal variations usually lead to the environments changed periodically as time increases, which will present the form that intrinsic growth rates of species are periodic in time mathematically. The following diffusive advective competition model  with different free boundaries involving time periodic environmental heterogeneity and impulsive interventions is proposed:
\begin{small}
 \begin{eqnarray}
\left\{
\begin{array}{lll}
u_t=Du_{xx}-\alpha u_x+\gamma u(1+\varepsilon_1(t)-u-kv),\; &\ t\in((n\tau)^+,(n+1)\tau],\ x\in(0,r(t)), \\[2mm]
v_t=v_{xx}-\beta v_x+v(1+\varepsilon_2(t)-v-hu),\; &\ t\in((n\tau)^+,(n+1)\tau],\ x\in(0,s(t)), \\[2mm]
u((n\tau)^+,x)=G(u(n\tau,x)),\; &\ x\in(0,r(n\tau)),\\[2mm]
v((n\tau)^+,x)=H(v(n\tau,x)),\; &\ x\in(0,s(n\tau)),\\[2mm]
u_x(t,0)=u(t,r(t))=v_x(t,0)=v(t,s(t))=0,\; &\ t\in(0,\infty),\\[2mm]
r'(t)=-\mu_1u_x(t,r(t)), s'(t)=-\mu_2v_x(t,s(t)),\; &\ t\in(n\tau,(n+1)\tau],\\[2mm]
r(0)=r_0,u(0,x)=u_0(x),\; &\ x\in[0,r_0],\\[2mm]
s(0)=s_0,v(0,x)=v_0(x),\; &\ x\in[0,s_0],
\end{array} \right.
\label{a01}
\end{eqnarray}
\end{small}
where $r(t)$ and $s(t)$ are the shifting boundaries with initial regions $r_0$ and $s_0$ to be determined, $\mu_1$ and $\mu_2$ are the expanding capability of species $u$ and $v$, respectively and  parameters $D$, $\gamma$, $k$, $h$, $\mu_1$, $\mu_2$, $r_0$ and $s_0$ are positive constants. Here, we mainly consider the advection  upward or downward along the gradient of resources with rates $\alpha,\beta\in \mathbb{R}$.  Initial functions $u_0(x)$ and $v_0(x)$ satisfy
\begin{eqnarray}
\left\{
\begin{array}{lll}
u_0(x)\in C^2([0,r_0]),\ u_0(x)>0 \ \textrm{for}\ x\in[0,r_0), \ u_0'(0)=u_0(r_0)=0,\\
v_0(x)\in C^2([0,s_0]),\ v_0(x)>0 \ \textrm{for}\ x\in[0,s_0),\ v_0'(0)=v_0(s_0)=0.
\end{array} \right.
\label{a02}
\end{eqnarray}
The terms $\varepsilon_i(t) (i=1,2)$ represent the environmental perturbations felt by species $u$ and $v$, respectively, which are caused by the seasonal variations of environments, and satisfy the condition

\medskip
$(\mathcal{H}1)$ \
 $\varepsilon_i(t+\tau)=\varepsilon_i(t)$, $\int_0^\tau \varepsilon_i(t)dt=0 $ and $|\varepsilon_i(t)|<1$,  $i=1,2$,
 \medskip\\
where the positive constant $\tau(>0)$ is time period. Different impulsive control strategies are carried out on two competing species at every time $t=n\tau$ with $n=0,1,2,...$ unless otherwise specified. The pulse functions $G(u)$ and $H(v)$ have the following properties:

\medskip
$(\mathcal{H}2)$  \ $G(u), H(v)\in C^2([0,+\infty))$, $G(0)=H(0)=0$, $G(u), H(v)>0$, $G'(u), H'(v)\geq 0$ for $u, v>0$, $G'(0),H'(0)>0$, $G(u)/u$ and $H(v)/v$ are nonincreasing with positive $u$ and $v$, respectively;

$(\mathcal{H}3)$\ there exist positive constants $D_i$, $\nu_i>1$ and small $\sigma_i$ $(i=1,2)$ such that $G(u)\geq G'(0)u-D_1u^{\nu_1}$ for $0\leq u\leq \sigma_1$ and $H(v)\geq H'(0)v-D_2v^{\nu_2}$ for $0\leq v\leq \sigma_2$.
\medskip\\
It should be mentioned that  $(\mathcal{H}1)$ and $(\mathcal{H}2)$ are the natural assumptions about the environmental perturbations terms and  pulse functions, respectively, and the assumption $(\mathcal{H}3)$ is the required condition in the construction of the lower solution, which is firstly proposed by Lewis and Li in \cite{lewisli} and can be also seen in \cite{menglin21}.

There are some more detailed explanations of model \eqref{a01}. Two competing species follow the  competition relationships modeled by the first two equations in  \eqref{a01} to grow and invade for $t\in ((n\tau)^+,(n+1)\tau]$. It means that this stage takes the initial values at initial time $t=(n\tau)^+$, where $t=(n\tau)^+$ is the right limit mathematically, and means the moment immediately after different human interventions take place.

If $\alpha=\beta=0$, $\varepsilon_i(t)=0 (i=1,2)$, $G(u)=u$ and $H(v)=v$, then model \eqref{a01} is reduced to the diffusive competition model with two different free boundaries, which was treated in \cite{guowu2015,liuwang2019,wangzhang, wu2015jde} and some related references. Guo and Wu \cite{guowu2015} considered the strong-weak competition case and Wu \cite{wu2015jde} took weak competition case into consideration, and they gave the existence and uniqueness of solutions and dynamics of species. Their works were improved by Wang and Zhang \cite{wangzhang}. Moreover, spreading speeds of species when two competitors coexist eventually were given in \cite{liuwang2019}.

In addition to some literatures cited above, there have been much more rich and extensive works about competition free boundaries problems. The introduction of seasonal environmental perturbations and periodic pulses may make analysis and results become complex and different,  and will raise some problems naturally.  What are the new criteria for classifying the dynamics of two competing species?  Whether and how can periodic pulses alter the competition outcomes, even lead to the coexistence of two competing species? How do advection, pulses and environmental perturbations affect spreading speeds when two competitors invade successfully? Whether do two competitors feel the same influence of environmental perturbations?

The propose of paper is to understand the effects of seasonal environmental perturbations and periodic pulses. The paper is organized as follows.

Section 2 is devoted to showing some preliminaries including existence and uniqueness of solution to problem \eqref{a01}, properties of principal eigenvalue and long time behaviors of one single problem \eqref{b01}. Since  periodic impulsive eigenvalue problem are related to advection and pulses, the explicit expressions of  principal eigenvalue can not be given. To obtain the relationships of  the principal eigenvalue of periodic impulsive eigenvalue problem with respect to factors including diffusion, length of habitat and advection,  we mainly investigate a fundamental eigenvalue problem instead by separating the variables of time and space. Two types of upper solutions are constructed according to different advection rates, and  dynamics of one single problem \eqref{b01} are then obtained for the small or negative advection case in Theorems 2.7 and 2.8.

In Section 3, we investigate the spreading and vanishing of competitors for the small or negative advection case, and obtain that competition outcomes of problem \eqref{a01} can be classified into four types: co-extinction, co-existence, spreading of species $u$ only and spreading of  species $v$ only, see Theorems 3.1 and 3.7.

Section 4 deals with some sufficient conditions on pulses  for spreading or vanishing, and  gives a minimal pulse intensity for the spreading of species under the linear impulsive function case. For small advection case, we also give some sufficient conditions on initial habitats and a sharp expanding capability under some specific initial region conditions. Meanwhile, some rough estimates of spreading speeds of competitors are obtained as advection is nonnegative.  Our results covers the results in aforementioned researches about competition problems without pulses and  are extended to the negative or small advection cases, which gives new perspectives about the effects of pulses on spreading or vanishing of competitors. And some techniques are modified and improved to overcome the difficulties induced by the combination of advection and pulses.

 Numerical simulations about the effects of seasonal environmental perturbations and periodic pulses are presented in Section 5. We mainly consider the small advection case and carry out impulsive interventions only on species $u$ for convenience. Some interesting and biologically valuable results are illustrated in the simulations and provide the suggestion about the design of human intervention strategies. Also, some observations about the effects of seasonal environmental perturbations different from those in \cite{khan} are shown in the simulations and show that the superior $u$ is more sensitive to environmental perturbations as a result of impulsive harvesting implemented only in the superior. Section 6 includes some discussions and views on the future work.

\section{Preliminary results}
Some preliminary results are shown in this section, which will be useful to deal with dynamics of the problem \eqref{a01}. Existence and uniqueness of solution to the problem \eqref{a01} and comparison principles are given in the first subsection. To investigate the long time behaviors of the solution to the problem \eqref{a01},  we then consider a general single species free boundary problem \eqref{b01},  study corresponding eigenvalue problem, and establish the results about  spreading and vanishing of species.

\subsection{Existence, uniqueness and comparison principles}
This subsection is devoted to giving the existence, uniqueness and some estimates of solution to problem \eqref{a01}. A comparison principle and some properties are also presented, which will play an important role in the discussion of long time behaviors of  the solution to the problem \eqref{a01}.
\begin{thm}
Problem \eqref{a01} admits a unique positive solution $(u,v;r,s)$ for all $t>0$. And, $u(t,x)\in PC^{1,2}((0,\infty)\times[0,r(t)]):=\{u(t,x)\in C^{1,2}
((n\tau,(n+1)\tau]\times[0,r(t)])\}$, $v(t,x)\in PC^{1,2}((0,\infty)\times[0,s(t)])$, and $r(t),s(t)\in C((0,\infty))\bigcap C^1((n\tau,(n+1)\tau])$. Furthermore, there exist positive constants $M_1$, $M_2$, $C_1$ and $C_2$ such that
\[
0<u(t,x)\leq M_1 \ \textrm{for}\ t>0,\ x\in[0,r(t)),
\]
\[
0<v(t,x)\leq M_2 \ \textrm{for}\ t>0,\ x\in[0,s(t)),
\]
\begin{equation}
 0<r'(t)\leq C_1\mu_1, \ 0<s'(t)\leq C_2\mu_2 \ \textrm{for}\ t>0 \ \textrm{and}\ t\neq n\tau.
\label{b051}
\end{equation}
\end{thm}
\begin{pf}
The existence and uniqueness of positive solution to the problem \eqref{a01} can be similarly obtained as Theorem 2.1 in \cite{menglin22} and Theorem 1 in \cite{guowu2015}, and are omitted here.

We only give some estimates of $u$ and $r$ in the sequel.
It is known as Theorem 2.2 in \cite{menglin22} that $u(t,x)\leq \max\{\|u_0\|_\infty, 1+\varepsilon_1^M\}$ with $\varepsilon_1^M:=\max_{t\in[0,\tau]}\varepsilon_1(t)$ provided that $0<G(u)/u\leq1$ for $u>0$. When $G(u)/u>1$ for $u>0$, it follows from comparison principle that $u(t,x)\leq Mw(t)$ for $t>0$ and $x\in[0,r(t)]$, where $M$ is a large constant such that $Mw(0)\geq u_0(x)$ for $x\in [0,r(t)]$, and $w(t)$ is the unique positive periodic solution to problem
\[
\left\{
\begin{array}{lll}
w_t=\gamma w(1+\varepsilon_1^M-w),\; &\  t\in((n\tau)^+,(n+1)\tau], \\[2mm]
w((n\tau)^+)=G'(0)w(n\tau),\\[2mm]
w(t+\tau)=w(t), \; &\  t>0.
\end{array} \right.
\]
Direct calculations give that
\[
w(t)=\frac{Ae^{\gamma A(t-n\tau)}G'(0)w(n\tau)}{A+(e^{\gamma A(t-n\tau)}-1)G'(0)w(n\tau)},\ t\in((n\tau^+),(n+1)\tau]
\]
with $A=1+\varepsilon_1^M$ and $w(n\tau)=\frac{A(e^{\gamma A\tau}-1/G'(0))}{e^{\gamma A\tau}-1}$, which is independent of $n$ because of  periodicity.
Since $G(u)$ satisfies $(\mathcal{H}2)$, we have $G'(0)>1$ and $e^{\gamma A\tau}G'(0)-1>0$, which yields that $w(n\tau)>A$ and
\[
\begin{array}{llllll}
w(t)=\frac{Ae^{\gamma A(t-n\tau)}(e^{\gamma A \tau}G'(0)-1)}{(e^{\gamma A \tau}-1)+(e^{\gamma A(t-n\tau)}-1)(e^{\gamma A \tau}G'(0)-1)}>\frac{Ae^{\gamma A(t-n\tau)}(e^{\gamma A \tau}G'(0)-1)}{(e^{\gamma A \tau}G'(0)-1)+(e^{\gamma A(t-n\tau)}-1)(e^{\gamma A \tau}G'(0)-1)}=A
\end{array}
\]
for $t\in((n\tau)^+,(n+1)\tau]$. Then, we have
 $w(t)>1+\varepsilon_1^M$ for $t>0$. And by observing the explicit expression of $w(t)$, we can obtain that $w(t)$ is bounded and there exists a positive constant $m_1$ such that $w(t)\leq m_1$ for $t>0$. Therefore, $u(t,x)\leq Mw(t)\leq  Mm_1=M_1$ for $t>0$ and $x\in[0,r(t))$. Similarly,  $v(t,x)\leq M_2$ for $t>0$ and $x\in[0,s(t))$.

By virtue of strong maximum principle, we obtain $u(t,x)>0$ and $u_x(t,r(t))<0$ for $t>0$ and $x\in[0,r(t))$. Then, $r'(t)=-\mu_1u_x(t,r(t))>0$ for $t>0$ and $t\neq n\tau$. $s'(t)>0$ can be similarly deduced. In what follows, we aim to construct an upper solution to prove $r'(t)\leq C_1\mu_1$ for $t>0$.  When $G(u)/u>1$,  define
\[
\Omega_K:=\{(t,x):t>0, r(t)-K^{-1}<x<r(t)\}
\]
and
\begin{equation}
\widetilde{u}(t,x)=Mw(t)[2K(r(t)-x)-K^2(r(t)-x)^2], \ (t,x)\in \Omega_K,
\label{b052}
\end{equation}
where $M$ and $w(t)$ are defined above.
Recalling  $w(t)>1+\varepsilon_1^M$, we have $w_t<0$, which yields
\[
\begin{array}{llllll}
\widetilde{u}_t&=Mw_t[2K(r(t)-x)-K^2(r(t)-x)^2]+Mw[2Kr'(t)-2K^2(r(t)-x)r'(t)]\\
&\geq M\gamma w(1+\varepsilon_1^M-w).
\end{array}
\]
By direct calculations, one has
\[
\begin{array}{llllll}
&-2KMw(t)\leq\widetilde{u}_x=2KMw(t)[K(r(t)-x)-1]\leq0,\\
&\widetilde{u}_{xx}=-2K^2Mw(t).
\end{array}
\]
Recalling $\widetilde{u}\leq Mw(t)$ and $1+\varepsilon_1^M<w(t)\leq m_1$ lead to
\[
\begin{array}{llllll}
\widetilde{u}_t-D\widetilde{u}_{xx}+\alpha\widetilde{u}_x
-\gamma\widetilde{u}(1+\varepsilon_1(t)-\widetilde{u}-kv)\\
\geq M\gamma w(1+\varepsilon_1^M-w)+2K^2DMw(t)-2K\alpha Mw(t)-\gamma(1+\varepsilon_1^M)\widetilde{u}\\
\geq M w(t)(-\gamma m_1+2DK^2-2\alpha K)\geq 0
\end{array}
\]
if $\alpha\geq0$ and $K\geq \frac{\alpha+\sqrt{\alpha^2+2D\gamma m_1}}{2D}$, while $\alpha<0$ leads to
\[
\begin{array}{llllll}
\widetilde{u}_t-D\widetilde{u}_{xx}+\alpha\widetilde{u}_x
-\gamma\widetilde{u}(1+\varepsilon_1(t)-\widetilde{u}-kv)\\
\geq M\gamma w(1+\varepsilon_1^M-w)+2K^2DMw(t)-\gamma(1+\varepsilon_1^M)\widetilde{u}\\
\geq M w(t)(-\gamma m_1+2DK^2)\geq 0
\end{array}
\]
provided that $K\geq \sqrt{\frac{\gamma m_1}{2D}}$.
 Owing to the monotonicity of $G(u)/u$,
\[
\widetilde{u}((n\tau)^+,x)=G'(0)\widetilde{u}(n\tau,x)\geq G(\widetilde{u}(n\tau,x)).
\]
And, we have
\[
\widetilde{u}(t,r(t))=0=u(t,r(t)),\ \widetilde{u}(t,r(t)-K^{-1})=Mw(t)\geq u(t,r(t)-K^{-1}).
\]
If we can select $K$ such that $\widetilde{u}(0,x)\geq u_0(x)$ for $x\in[r_0-K^{-1},r_0]$, then by virtue of comparison principle, we have $\widetilde{u}(t,x)\geq u(t,x)$ for $(t,x)\in\Omega_K$, which can deduce that
\[
u_x(t,r(t))\geq\widetilde{u}_x(t,r(t))=-2KMw(t)\]
and \[r'(t)=-\mu_1u_x(t,r(t))\leq 2K\mu_1Mw(t)\leq 2KM_1\mu_1.
\]

Choose
\[
K:=\left\{
\begin{array}{lll}
\max \{\frac{\alpha+\sqrt{\alpha^2+2D\gamma m_1}}{2D},-\frac{4\min_{x\in[0,r_0]} u_0'(x)}{3M(1+\varepsilon_1^M)}\} \ for \ \alpha\geq0,\\
\max\{\sqrt{\frac{\gamma m_1}{2D}},-\frac{4\min_{x\in[0,r_0]} u_0'(x)}{3M(1+\varepsilon_1^M)}\} \ for \ \alpha<0.
\end{array}\right.
\]
Then, for $x\in[r_0-K^{-1},r_0-(2K)^{-1}]$, $\widetilde{u}(0,x)\geq \frac{3}{4}Mw(0)\geq\frac{3}{4}M (1+\varepsilon_1^M)$ and $u_0(x)\leq-\min_{x\in[0,r_0]} u_0'(x)K^{-1}\leq\frac{3}{4}M (1+\varepsilon_1^M)$, which implies that
\[
\widetilde{u}(0,x)\geq u_0(x),\ x\in[r_0-K^{-1},r_0-(2K)^{-1}].
\]
And, for $x\in[r_0-(2K)^{-1},r_0]$, we have
$\widetilde{u}_x(0,x)=-2KMw(0)[1-K(r_0-x)]\leq -KMw(0)\leq-KM(1+\varepsilon_1^M)$.
One can deduce that $\widetilde{u}_x(0,x)\leq u_0'(x)$ for $x\in[r_0-(2K)^{-1},r_0]$, which together with $\widetilde{u}(0,r_0)=u_0(r_0)$ yields
\[
\widetilde{u}(0,x)\geq u_0(x),\ x\in[r_0-(2K)^{-1},r_0].
\]
Thus, \eqref{b051} can be obtained for the case $G(u)/u>1$. When $G(u)/u\leq1$,  \eqref{b051} can be verified by replacing $Mw(t)$ in the expression \eqref{b052} of $\widetilde{u}(t,x) $ by $M_1$, where $M_1$ is the bound of $u$. We here omit the proof for the case $G(u)/u\leq1$ since it is similar to Lemma 2.2 in \cite{dyhlzg}.
\end{pf}

We now present a comparison principle for problem \eqref{a01}. For $t\in [0,\tau]$, taking $t=0^+$ as a new initial time, due to
$\overline{u}(0^+,x)=G(\overline{u}(0,x))\geq G(u(0,x))=u(0^+,x)$ and  $\underline{v}(0^+,x)=H(\underline{v}(0,x))\leq H(v(0,x))=v(0^+,x)$, we can obtain  the follow comparison principle holds for $t\in (0^+,\tau]$ similarly as Lemma 2.3 in \cite{wu2015jde}. Then, the monotonicity of $G(u)$ and $H(v)$ with $u$ and $v$ give $\overline{u}(\tau^+,x)\geq u(\tau^+,x)$ and  $\underline{v}(\tau^+,x)\leq v(\tau^+,x)$.  For $t\in [\tau,2\tau]$ , we take $t=\tau^+$ as a new initial time and use comparison principle similarly as Lemma 2.3 in \cite{wu2015jde}, which yields the follow comparison principle holds for $t\in (\tau^+,2\tau]$. Repeating the above process, we can eventually obtain the following comparison principle for all $t>0$.
\begin{lem}
Let $\overline{r}, \underline{s}\in C((0,\infty))\bigcap C^1((n\tau,(n+1)\tau])$, $\overline{u}\in PC^{1,2}((0,\infty)\times[0,\overline{r}(t)])$ and $\underline{v}\in PC^{1,2}((0,\infty)\times[0,\underline{s}(t)])$. If
\begin{small}
\begin{eqnarray}
\left\{
\begin{array}{lll}
\overline{u}_t\geq D\overline{u}_{xx}-\alpha \overline{u}_x+\gamma \overline{u}(1+\varepsilon_1(t)-\overline{u}-k\underline{v}),\; &\ t\in((n\tau)^+,(n+1)\tau],\ x\in(0,\overline{r}(t)), \\[2mm]
\underline{v}_t\leq\underline{v}_{xx}-\beta \underline{v}_x+\underline{v}(1+\varepsilon_2(t)-\underline{v}-h\overline{u}),\; &\ t\in((n\tau)^+,(n+1)\tau],\ x\in(0,\underline{s}(t)), \\[2mm]
\overline{u}((n\tau)^+,x)\geq G(\overline{u}(n\tau,x)),\; &\ x\in(0,\overline{r}(n\tau)),\\[2mm]
\underline{v}((n\tau)^+,x)\leq H(\underline{v}(n\tau,x)),\; &\ x\in(0,\underline{s}(n\tau)),\\[2mm]
\overline{u}_x(t,0)=\overline{u}(t,\overline{r}(t))=\underline{v}_x(t,0)
=\underline{v}(t,\underline{s}(t))=0,\; &\ t\in(0,\infty),\\[2mm]
\overline{r}'(t)\geq-\mu_1\overline{u}_x(t,\overline{r}(t)), \underline{s}'(t)\leq-\mu_2\underline{v}_x(t,\underline{s}(t)),\; &\ t\in(n\tau,(n+1)\tau],\\[2mm]
\overline{r}(0)\geq r_0,\overline{u}(0,x)\geq u_0(x),\; &\ x\in[0,r_0],\\[2mm]
\underline{s}(0)\leq s_0,\underline{v}(0,x)\leq v_0(x),\; &\ x\in[0,s_0],
\end{array} \right.
\label{d2}
\end{eqnarray}
\end{small}
then
$r(t)\leq \overline{r}(t),\ s(t)\geq \underline{s}(t)$ for $t\in [0,\infty)$,
$u(t,x)\leq \overline{u}(t,x)$ for $t\in [0,\infty)$ and $x\in[0,r(t)]$, and $
v(t,x)\geq \underline{v}(t,x)$ for $t\in [0,\infty)$ and $x\in[0,\underline{s}(t)]$.
\end{lem}
The comparison principle still holds true for  \eqref{d2} with reversed inequalities.

To stress on the dependence of $(u,v,r,s)$ on $\mu_1$ and $\mu_2$, denote $u^{\mu_1,\mu_2}$, $v^{\mu_1,\mu_2}$, $r^{\mu_1,\mu_2}$ and $u^{\mu_1,\mu_2}$ in the following lemma.
\begin{lem}
Suppose all the parameters in the problem \eqref{a01} except $\mu_1$ and $\mu_2$ are fixed, and $\mu_1\leq \tilde{\mu}_1$, $\mu_2\geq \hat{\mu}_2$. Then
\[
\begin{array}{llllll}
r^{\tilde{\mu}_1,\hat{\mu}_2}(t) \geq r^{\mu_1,\mu_2}(t),\ s^{\tilde{\mu}_1,\hat{\mu}_2}(t) \leq s^{\mu_1,\mu_2}(t) \ & \textrm{for} \ t>0,\\ [2mm]
u^{\tilde{\mu}_1,\hat{\mu}_2}(t,x) \geq u^{\mu_1,\mu_2}(t,x) \ & \textrm{for}\ t>0,\ 0<x<r^{\mu_1,\mu_2}(t),\\ [2mm]
v^{\tilde{\mu}_1,\hat{\mu}_2}(t,x) \leq v^{\mu_1,\mu_2}(t,x) \ & \textrm{for}\ t>0,\ 0<x<s^{\tilde{\mu}_1,\hat{\mu}_2}(t).
\end{array}
\]
\end{lem}
\subsection{The principal eigenvalues}
To study the dynamics of competition free boundary problem \eqref{a01},
the following general single species problem with free boundary is considered
\begin{eqnarray}
\left\{
\begin{array}{lll}
u_t=Du_{xx}-\alpha u_x+\gamma u(a(t)-u),\; &\ t\in((n\tau)^+,(n+1)\tau],\ x\in(0,r(t)), \\[2mm]
u((n\tau)^+,x)=G(u(n\tau,x)),\; &\ x\in(0,r(n\tau)),\\[2mm]
u_x(t,0)=u(t,r(t))=0,\; &\ t\in(0,\infty),\\[2mm]
r'(t)=-\mu_1u_x(t,r(t)),\; &\ t\in(n\tau,(n+1)\tau],\\[2mm]
r(0)=r_0,u(0,x)=u_0(x),\; &\ x\in[0,r_0],
\end{array} \right.
\label{b01}
\end{eqnarray}
where $a(t)$ is a $\tau-$periodic function in time $t$ and $a(t)\in C^\alpha((n\tau,(n+1)\tau])$ with $\alpha \in (0,1)$.

If the domain is fixed, that is, $r(t)=l$, then problem \eqref{b01} becomes
\begin{eqnarray}
\left\{
\begin{array}{lll}
u_t=Du_{xx}-\alpha u_x+\gamma u(a(t)-u),\; &\ t\in((n\tau)^+,(n+1)\tau],\ x\in(0,l), \\[2mm]
u((n\tau)^+,x)=G(u(n\tau,x)),\; &\ x\in(0,l),\\[2mm]
u_x(t,0)=u(t,l)=0,\; &\ t\in(0,\infty),\\[2mm]
u(0,x)=u_0(x),\; &\ x\in[0,l].
\end{array} \right.
\label{b02}
\end{eqnarray}
For problem \eqref{b02}, we consider the corresponding periodic eigenvalue problem
\begin{eqnarray}
\left\{
\begin{array}{lll}
\phi_t=D\phi_{xx}-\alpha \phi_x+\gamma a(t) \phi+\lambda\phi,\; &\ t\in(0^+,\tau],\ x\in(0,l), \\[2mm]
\phi(0^+,x)=G'(0)\phi(0,x),\; &\ x\in(0,l),\\[2mm]
\phi_x(t,0)=\phi(t,l)=0,\; &\ t\in[0,\tau],\\[2mm]
\phi(0,x)=\phi(\tau,x),\; &\ x\in[0,l].
\end{array} \right.
\label{b03}
\end{eqnarray}
 To investigate periodic eigenvalue problem \eqref{b03}, the following fundamental eigenvalue problem is firstly taken into consideration
\begin{eqnarray}
\left\{
\begin{array}{lll}
-D\varphi_{xx}+\alpha \varphi_x=\lambda\varphi,\; &\  x\in(0,l), \\[2mm]
\varphi_x(0)=\varphi(l)=0.
\end{array} \right.
\label{b04}
\end{eqnarray}
Let $\lambda^*(D,\alpha,(0,l))$ be the principal eigenvalue of eigenvalue problem \eqref{b04}.
By variational method,
\[
\lambda^*=\inf_{\varphi\not\equiv0,\varphi\in \mathbb{H}(l)}\frac{\int_0^lDe^{-\frac{\alpha}{D}x}(\varphi_x)^2dx}
{\int_0^le^{-\frac{\alpha}{D}x}\varphi^2dx},
\]
where $\mathbb{H}(l):=\{\varphi>0:\varphi\in H^1((0,l)),\varphi'(0)=0,\varphi(l)=0 \}$.
Then properties of $\lambda^*$ with respect to $l$, $D$ and $\alpha$ are shown as follows.
\begin{lem}
Let $\lambda^*(D,\alpha,(0,l))$ be the principal eigenvalue of eigenvalue problem \eqref{b04}, then the following assertions hold.

$(i)$ $\lambda^*(D,\alpha,(0,l))$ is strictly decreasing with respect to $l$. Furthermore, $\lim\limits_{l\rightarrow0^+}\lambda^*(D,\alpha,(0,l))$\\$=+\infty$ and $\lim\limits_{l\rightarrow+\infty}\lambda^*(D,\alpha,(0,l))=\frac{\alpha^2}{4D}$ for any $\alpha\in \mathbb{R}$.

$(ii)$ $\lim\limits_{D\rightarrow0^+}\lambda^*(D,\alpha,(0,l))=+\infty$ for $\alpha\neq0$,
$\lim\limits_{D\rightarrow0^+}\lambda^*(D,\alpha,(0,l))=0$ for $\alpha=0$,
 and $\lim\limits_{D\rightarrow+\infty}\lambda^*(D,\alpha,(0,l))=+\infty$ for any $\alpha\in \mathbb{R}$.

$(iii)$ $\lim\limits_{\alpha\rightarrow0}\lambda^*(D,\alpha,(0,l))=\frac{D\pi^2}{4l^2}$.
\end{lem}
\begin{pf}
 We first prove the monotonicity of $\lambda^*(D,\alpha,(0,l))$ with  $l$. Given any $l_1<l_2$, we aim to show that $\lambda^*(D,\alpha,(0,l_1))>\lambda^*(D,\alpha,(0,l_2))$. Denote $(\lambda^*(D,\alpha,(0,l_1)),\varphi_1)$ and $(\lambda^*(D,\alpha,(0, l_2)),\varphi_2)$ be the principal eigen-pair of eigenvalue problem \eqref{b04} with $l=l_1$ and $l=l_2$, respectively. Multiplying the equations of $(\lambda^*(D,\alpha,(0,l_1)),\varphi_1)$ and  $(\lambda^*(D,\alpha,(0,l_2)),\varphi_2)$ by $e^{-\frac{\alpha}{D}x}\varphi_2$ and $e^{-\frac{\alpha}{D}x}\varphi_1$, respectively, subtracting  these two equations and integrating the results from $0$ to $l_1$, we have
\[
\begin{array}{llllll}
&(\lambda^*(D,\alpha,(0,l_1))-\lambda^*(D,\alpha,(0,l_2)))\int_0^{l_1}
e^{-\frac{\alpha}{D}x}\varphi_1\varphi_2dx\\
 &=\int_0^{l_1}e^{-\frac{\alpha}{D}x}[-D(\varphi_{1,xx}\varphi_2-\varphi_{2,xx}\varphi_1)
 +\alpha(\varphi_{1,x}\varphi_2-\varphi_{2,x}\varphi_1)]dx\\
 &=-D\int_0^{l_1}(e^{-\frac{\alpha}{D}x}\varphi_{1,x}\varphi_2-
 e^{-\frac{\alpha}{D}x}\varphi_{2,x}\varphi_1)_xdx\\
 &=-De^{-\frac{\alpha}{D}l_1}\varphi_{1,x}(l_1)\varphi_2(l_1)
\end{array}
\]
by using the boundary conditions in problem \eqref{b04}. Noting that $\varphi_{1,x}(l_1)<0$ by Hopf's boundary lemma, we have $(\lambda^*(D,\alpha,(0,l_1))-\lambda^*(D,\alpha,(0,l_2)))\int_0^{l_1}
e^{-\frac{\alpha}{D}x}\varphi_1
\varphi_2dx>0$, which deduces that $\lambda^*(D,\alpha,(0,l_1))>\lambda^*(D,\alpha,(0,l_2))$. Thus, $\lambda^*(D,\alpha,(0,l))$ is strictly decreasing with respect to $l$.

To investigate the asymptotic behaviors of $\lambda^*(D,\alpha,(0,l))$ with parameters $l$, $D$ and $\alpha$, we carried out careful calculations of eigenvalue problem \eqref{b04}.  It is known that $\lambda^*>0$. According to the roots of the equation
 \begin{equation}
 -D r^2+\alpha r-\lambda^*=0,
\label{b071}
\end{equation}
 we have three cases to be considered.

Case 1:  $\lambda^*<\frac{\alpha^2}{4D}$. In this case, equation \eqref{b071} has two different real roots $r_{1,2}=\frac{\alpha\pm\sqrt{\alpha^2-4D\lambda^*}}{2D}$, and the eigen-function satisfies
\begin{equation}
\varphi(x)=e^{\frac{\alpha}{2D}x}(c_1e^{\frac{\sqrt{\alpha^2-4D\lambda^*}}{2D}x}
+c_2e^{\frac{-\sqrt{\alpha^2-4D\lambda^*}}{2D}x}).
\label{bp03}
\end{equation}
Using the conditions $\varphi_x(0)=\varphi(l)=0$ in problem \eqref{b04} yields
that $c_2=-\frac{\alpha+z_1}
{\alpha-z_1}c_1$, and principal eigenvalue $\lambda^*$ of problem \eqref{b04} exists if and only if $\lambda^*$ satisfying
\begin{equation}
e^{\frac{l}{D}z_1}=\frac{\alpha+z_1}
{\alpha-z_1}
\label{bp01}
\end{equation}
is minimal and positive, and corresponding  eigenfunction $\varphi(x)$ is positive, where $z_1=\sqrt{\alpha^2-4D\lambda^*}$. Solving \eqref{bp01} by considering $y=e^{\frac{l}{D}z_1}$ and $y=\frac{\alpha+z_1}
{\alpha-z_1}$, we have the unique $\lambda^*$ corresponding to positive $\varphi(x)$ exists if and only if $\alpha>\frac{2D}{l}$, see Fig. \ref{sy} $(a)$.

Case 2: $\lambda^*>\frac{\alpha^2}{4D}$. Equation \eqref{b071} admits a pair of conjugate complex roots $r_{1,2}=\frac{\alpha\pm z_2i}{2D}$, and the eigenfunction $\varphi(x)$ of problem \eqref{b04} satisfies
\[
\varphi(x)=e^{\frac{\alpha}{2D}x}(c_3\cos(\frac{z_2}{2D}x)
+c_4\sin(\frac{z_2}{2D}x))
\]
with $z_2=\sqrt{4D\lambda^*-\alpha^2}$ and $c_4=-\alpha c_3/z_2$. If $\alpha=0$, using $\varphi(l)=0$ gives that $\cos(\frac{l\sqrt{4D\lambda^*}}{2D})=0$, and then $\lambda^*(D,0,(0,l))=\frac{D\pi^2}{4l^2}$.
If $\alpha\neq0$, $\varphi(l)=0$ implies that
\begin{equation}
\tan(\frac{l}{2D}z_2)=\frac{z_2}{\alpha}.
\label{bp02}
\end{equation}
It is known that the minimal value $\lambda^*$ satisfying \eqref{bp02} is principal eigenvalue of problem \eqref{b04} with positive eigenfunction $\varphi(x)$.
Solving \eqref{bp02}, we have the principal eigenvalue $\lambda^*$ if $0<\alpha<\frac{2D}{l}$ or $\alpha<0$, see Fig. \ref{sy} $(b)$ and $(c)$.

Case 3: $\lambda^*=\frac{\alpha^2}{4D}$. We have the roots $r_1$ and $r_2$ of equation \eqref{b071} is the same, that is $r_1=r_2=\frac{\alpha}{2D}$, and
\[
\varphi(x)=e^{\frac{\alpha}{2D}x}(c_5x+c_6)
\]
with $c_6=-\frac{2D}{\alpha}c_5$. In this case, it follows from $\varphi(l)=0$ that $\alpha=\frac{2D}{l}$.

\begin{figure}[ht]
\quad
\begin{minipage}{0.3\linewidth}
\centerline{\includegraphics[width=5cm]{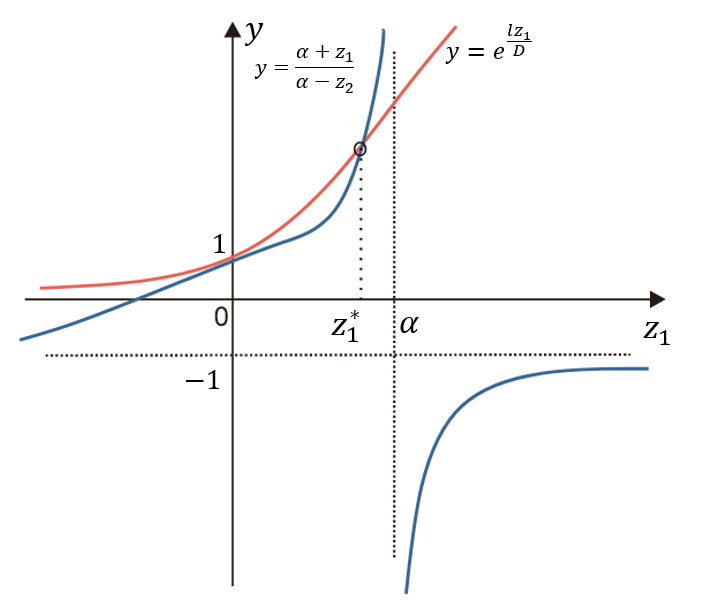}}
\centerline{\small{(a) Case 1: $\alpha>\frac{2D}{l}$}}
\end{minipage}
\quad
\begin{minipage}{0.3\linewidth}
\centerline{\includegraphics[width=5cm]{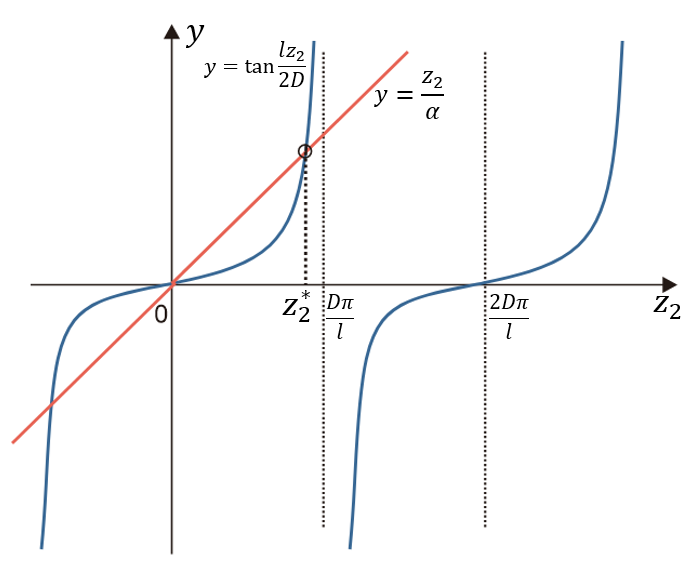}}
\centerline{\small{(b) Case 2: $0<\alpha<\frac{2D}{l}$}}
\end{minipage}
\quad
\begin{minipage}{0.3\linewidth}
\centerline{\includegraphics[width=5cm]{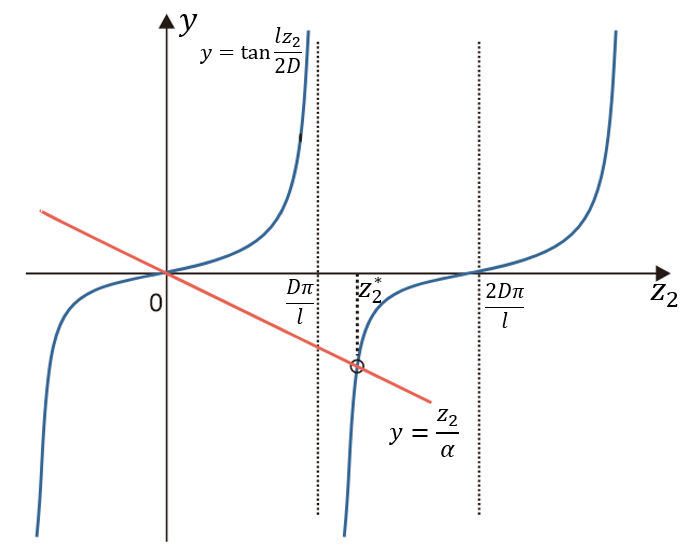}}
\centerline{\small{(c) Case 2: $\alpha<0$}}
\end{minipage}
\caption{\scriptsize Graphs $(a)-(c)$ are the schematic diagrams of solving \eqref{bp01} for the case $1$ and solving \eqref{bp02} for the case $2$, respectively. In $(a)$, the red line stands for the curve of $y=e^{\frac{l}{D}z_1}$, while the blue line is the hyperbolic curve of $y=\frac{\alpha+z_1}
{\alpha-z_1}$. The horizontal coordinate $z_1^*$ of the first positive intersection point in $(a)$ satisfies  \eqref{bp01}. The red and blue lines in $(b)$ and $(c)$ represent the line of $y=\frac{z_2}{\alpha}$ and tangent line of $y=\tan(\frac{l}{2D}z_2)$, respectively. And the horizontal coordinate of the first positive intersection point is $z_2^*$, which satisfies \eqref{bp02}.    }
\label{sy}
\end{figure}

Recalling that $\lambda^*(D,0,(0,l))=\frac{D\pi^2}{4l^2}$, the properties of $\lambda^*$ with $D$ and $l$ can directly be obtained by the explicit expression of $\lambda^*$ as $\alpha=0$. In what follows, we investigate long time behaviors of  $\lambda^*$ under the case $\alpha\neq0$.

If $l\rightarrow0^+$, then $\lambda^*$ satisfies \eqref{bp02}. It can be seen from Fig. $1$ $(b)$ and $(c)$ that the  periods  $\frac{D\pi}{l}$ and $\frac{2D\pi}{l}$ of $\tan(\frac{l}{2D}z_2)$ tends to $\infty$ as $l\rightarrow0^+$, which implies that the horizontal coordinate $z_2^*$ of the first positive intersection point satisfies $z_2^*=\sqrt{4D\lambda^*-\alpha^2}\rightarrow \infty$ as $l\rightarrow0^+$. We have
$\lim\limits_{l\rightarrow0^+}\lambda^*(D,\alpha,(0,l))=+\infty$. If $l\rightarrow\infty$, we consider two situations, $\alpha>\frac{2D}{l}$ and $\alpha<0$. When $\alpha>\frac{2D}{l}$, it follows from Case $1$ that $\lambda^*$ satisfies \eqref{bp01}. From Fig. $1$ $(a)$, we have the horizontal coordinate $z_2^*$ of the first positive intersection point satisfies $z_1^*=\sqrt{\alpha^2-4D\lambda^*}\rightarrow 0$ as $l\rightarrow\infty$, which means that $\lim\limits_{l\rightarrow+\infty}\lambda^*(D,\alpha,(0,l))=\frac{\alpha^2}{4D}$. When $\alpha<0$, we consider the Case 2 and $\lambda^*$ satisfies $\eqref{bp02}$. It can be directly seen from Fig. $1$ $(c)$ that $\lim\limits_{l\rightarrow +\infty}z_2^*=0$, which yields that $\lim\limits_{l\rightarrow+\infty}\lambda^*(D,\alpha,(0,l))=\frac{\alpha^2}{4D}$.

If $D\rightarrow0^+$, we consider Case $1$ and $\alpha<0$ in Case 2, where $\lambda^*$ satisfies \eqref{bp01} and \eqref{bp02}, respectively. It follows from Fig. $1$ $(a)$ and $(c)$ that  as $D\rightarrow0^+$, the horizontal coordinate  of the first positive intersection point $z_1^*\rightarrow 0^+$ and $z_2^*\rightarrow 0^+$. We have $\lambda^*\rightarrow\frac{\alpha^2}{4D}$, which together with $D\rightarrow0^+$ yields that $\lim\limits_{D\rightarrow0^+}\lambda^*(D,\alpha,(0,l))=+\infty$. We next claim that $\lim\limits_{D\rightarrow+\infty}\lambda^*=\infty$. If not, we assume that $\lambda^*$ tends to a finite value as $D\rightarrow\infty$. It is easy to see that $\lambda^*$ satisfies \eqref{bp02} for $D\rightarrow\infty$,   $\lim\limits_{D\rightarrow\infty}\tan{\frac{lz_2}{2D}}=0$ and $\lim\limits_{D\rightarrow\infty}\frac{z_2}{\alpha}=\pm\infty$, which leads to a contradiction. Thus, we have $\lim\limits_{D\rightarrow+\infty}\lambda^*(D,\alpha,(0,l))=+\infty$.

If $\alpha\rightarrow0^\pm$, then $\lambda^*$ satisfies \eqref{bp02}.  Notice that $\lim\limits_{\alpha\rightarrow0^\pm}\frac{z_2}{\alpha}=\pm\infty$, which gives that as $\alpha\rightarrow0$, $\frac{l}{2D}z_2\rightarrow\frac{\pi}{2}$, that is, $\frac{l\sqrt{4D\lambda^*}}{2D}\rightarrow\frac{\pi}{2}$. We obtain that $\lim\limits_{\alpha\rightarrow0}\lambda^*(D,\alpha,(0,l))=\frac{D\pi^2}{4l^2}$.
\end{pf}
\begin{rmk}
Lemma 2.4 implies that when $\alpha=0$, $\lambda^*$ is strictly monotone  increasing on $D>0$, but $\lambda^*$ fails  to be monotone on $D>0$ when advection $\alpha\neq0$ is introduced. The long time behaviors of $\lambda^*$ with large advection have not be given in this paper. Interested readers can refer to \cite{berestycki,chenlou,peng} and references therein for  principal eigenvalue of elliptic eigenvalue problems with large advection and different boundary conditions.
\end{rmk}

We then turn to define the principal eigenvalue of the problem \eqref{b03} as in \cite{wuzhao22}. The following auxiliary problem is considered
\begin{eqnarray}
\left\{
\begin{array}{lll}
w_t=Dw_{xx}-\alpha w_x+(\gamma a(t)-C^*)w,\; &\ t\in((n\tau)^+,(n+1)\tau],\ x\in(0,l), \\[2mm]
w((n\tau)^+,x)=G'(0)w(n\tau,x),\; &\ x\in[0,l],\\[2mm]
w_x(t,0)=w(t,l)=0,\; &\ t>0,
\end{array} \right.
\label{d1}
\end{eqnarray}
where $C^*:=\max_{t\in[0,\tau]}a(t)+\ln(\max\{ 1/G'(0),1\})+1$.
 Problem \eqref{d1} without  pulse  generates a strongly continuous semigroup $P(t)$ on $E$ that is strongly positive, that is, $P(t)E_+\subseteq E_+$, where $E =C([0,L],\mathbb{R})$ and $E_+=\{\phi\in S:\ \phi\geq0, \ \forall x\in [0,L]\}$.  If $\psi\geq0$ has a nonempty support in $t>0$, then $[P(t)\psi](x)>0$.

 The principal eigenvalue $\lambda_1$  of the  problem \eqref{b03} can be defined as in \cite{wuzhao22, menglin23}. Let $\mathcal{P}$ be the  Poincar\'{e}
 time map of the problem \eqref{d1}, then $\mathcal{P}=P(t)\circ L$, where $P(\tau)$ is the time-$\tau$ solution map of the system \eqref{d1} without  pulse  and $L$ is defined as a linear operator $L(\psi)(x)=G'(0)\psi(x)$.
Define $Q:=\mathcal{P}(\tau)$ and $r(Q)$ be the spectral radius of $Q$, then it is known that $Q$ is compact and strongly positive, and $r(Q)>0$ by the choice of $C^*$, which together with Lemma 2.2 in \cite{wuzhao22} and Krein-Rutman theorem gives  that $r(Q)$ is a principal eigenvalue of $Q$ with a strongly positive eigenfunction $\phi^*>0$.

From Theorem 2.7 in \cite{liangzhangzhao}, $\lambda_1:=-\frac{\ln r(Q)}{\tau}-C^*$ is the principal eigenvalue of problem \eqref{b03}, and $\phi(t,x)$ is corresponding eigenfunction with $\phi(t,x):=e^{(\lambda_1+C^*) t}\mathcal{P}(t)\phi^*(t,x)$, where $\phi^*(t,x)$ is the solution to  \eqref{d1} with initial function $\phi^*(0,x)$.

For the specific problem \eqref{b03}, the explicit expression of defined principal eigenvalue $\lambda_1$ can be given in the Theorem 2.5.
\begin{thm}
Let $\lambda_1:=\lambda_1(D,\alpha, \gamma a(t), G'(0),(0,l))$ be the principal eigenvalue of the eigenvalue problem \eqref{b03}. Then
\[
\lambda_1(D,\alpha, \gamma a(t), G'(0),(0,l))=\lambda^*(D,\alpha, (0,l))-\frac{1}{\tau}\ln G'(0)-\gamma \overline{a},
\]
where $\lambda^*(D,\alpha, (0,l))$ is the principal eigenvalue of the  problem \eqref{b04} and $\overline{a}:=\frac{1}{\tau}\int_{0^+}^\tau a(t)dt$. Moreover,

$(i)$ $\lambda_1(D,\alpha, \gamma a(t), G'(0),(0,l))$ is strictly decreasing with respect to $l$ and $G'(0)$, and $\lim\limits_{l\rightarrow0^+}\lambda_1=+\infty$ and $\lim\limits_{l\rightarrow+\infty}\lambda_1=\frac{\alpha^2}{4D}-\frac{1}{\tau}
\ln G'(0)-\gamma \overline{a}$ for any $\alpha\in \mathbb{R}$;

$(ii)$ $\lim\limits_{D\rightarrow0^+}\lambda_1=+\infty$ for $\alpha\neq0$, $\lim\limits_{D\rightarrow0^+}\lambda_1=-\frac{1}{\tau}
\ln G'(0)-\gamma \overline{a}$ for $\alpha=0$, and $\lim\limits_{D\rightarrow+\infty}\lambda_1=+\infty$ for any $\alpha\in \mathbb{R}$;

$(iii)$ $\lim\limits_{\alpha\rightarrow0}\lambda_1=\frac{D\pi^2}{4l^2}-\frac{1}{\tau}
\ln G'(0)-\gamma \overline{a}$.
\end{thm}
\begin{pf}
Let $\phi(t,x)=f(t)\varphi(x)$, where $\varphi(x)$ is the eigenfunction of the problem \eqref{b04} corresponding to $\lambda^*$. Then
\[
\frac{f_t}{f}=\frac{D\varphi_{xx}}{\varphi}-\frac{\alpha\varphi_x}{\varphi}+
\gamma a(t)+\lambda_1.
\]
By integrating this equation from $0^+$ to $\tau$  and  using  problem \eqref{b04}, one has
\[
\int_{0^+}^\tau \frac{f_t}{f}=(-\lambda^*+\lambda_1)\tau+\gamma\int_{0^+}^\tau a(t)dt.
\]
Recalling $f(0^+)=G'(0)f(0)$,  we deduce $\lambda_1=\lambda^*-\frac{1}{\tau}\ln G'(0)-\gamma \overline{a}$. Therefore, the properties of $\lambda_1$ can be obtained by using Lemma 2.4.
\end{pf}

\subsection{Dynamics for a single species free boundary problem}
Using similar methods as  in Theorems 3.2, 3.4 and 3.5 in \cite{menglin21}, the threshold type dynamics of fixed boundary problem \eqref{b02} are  first shown.
\begin{lem}
Let $\lambda_1(D,\alpha,\gamma a(t),G'(0),(0,l))$ be the principal eigenvalue of the problem \eqref{b03}. For any given $l$ and nonnegative nontrivial initial function $u_0$, the solution $u(t,x)$ of the problem \eqref{b02} satisfies the following assertions:

$(i)$ if $\lambda_1(D,\alpha,\gamma a(t),G'(0),(0,l))\geq 0$, then $\lim\limits_{t\rightarrow\infty}u(t,x)=0$ uniformly for $x\in[0,l]$;

$(ii)$ if $\lambda_1(D,\alpha,\gamma a(t),G'(0),(0,l))< 0$, then $\lim\limits_{m\rightarrow\infty}u(t+m\tau,x)=U_l(t,x)$ for $(t,x)\in[0,\tau]\times[0,l]$, where $U_l(t,x)$ is the unique positive solution to the periodic problem
\begin{eqnarray}
\left\{
\begin{array}{lll}
(U_l)_t=D(U_l)_{xx}-\alpha (U_l)_x+\gamma U_l(a(t)-U_l),\; &\ t\in(0^+,\tau],\ x\in(0,l), \\[2mm]
U_l(0^+,x)=G(U_l(0,x)),\; &\ x\in(0,l),\\[2mm]
(U_l)_x(t,0)=U_l(t,l)=0,\; &\ t\in[0,\tau],\\[2mm]
U_l(0,x)=U_l(\tau,x),\; &\ x\in[0,l].
\end{array} \right.
\label{b041}
\end{eqnarray}
Moreover,
if $\lambda_1(D,\alpha,\gamma a(t),G'(0),(0,\infty))< 0$, then $\lim\limits_{l\rightarrow\infty}U_l(t,x)=U^*(t)$ locally uniformly for $(t,x)\in[0,\tau]\times[0,\infty)$, where $U^*(t)$ is the unique positive periodic solution to the periodic problem
\begin{eqnarray}
\left\{
\begin{array}{lll}
U_t=\gamma U(a(t)-U),\; &\  t\in(0^+,\tau], \\[2mm]
U(0^+)=G(U(0)),\\[2mm]
U(0)=U(\tau).
\end{array} \right.
\label{b05}
\end{eqnarray}
Here $\lambda_1(D,\alpha,\gamma a(t),G'(0),(0,\infty)):=\lim\limits_
{l\rightarrow\infty}\lambda_1(D,\alpha,\gamma a(t),G'(0),(0,l))$.
\end{lem}
\begin{pf}
The threshold dynamics of the solution  for the bounded domain $(0,l)$ can be obtained as  in \cite{menglin21}. We show the proofs of $\lim\limits_{l\rightarrow\infty}U_l(t,x)=U^*(t)$ is locally uniformly for $(t,x)\in[0,\tau]\times[0,\infty)$ here. Due to  $\lambda_1(D,\alpha,\gamma a(t),G'(0),(0,\infty))< 0$, problem \eqref{b041} with $(0,l)$ replaced by $(0,\infty)$ admits a positive periodic solution $U_\infty(t,x)$. It follows from comparison principle that the positive solution $U_l(t,x)$ to   problem \eqref{b041} is strictly increasing with respect to $l$, which together with regularity theory for parabolic equations yields that $\lim\limits_{l\rightarrow\infty}U_l(t,x)=U_\infty(t,x)$ is locally uniformly for $(t,x)\in[0,\tau]\times[0,\infty)$. It remains to show  the uniqueness of positive periodic $U_\infty(t,x)$.

Let $U_1$ and $U_2$ be two positive periodic solutions to the problem \eqref{b041} with $l=\infty$. Without loss of generality, we suppose $U_1>U_2$ for $t\in[0,\tau]$ and $x\in[0,\infty)$. Define sector
\[
\Lambda=\{\sigma\in[0,1]: \sigma U_1\leq U_2 \ \textrm{for}\ t\in[0,\tau], x\in[0,\infty) \}.
\]
We claim that $1\in\Lambda$. By contradiction, suppose $\sigma_0=\sup\Lambda<1$. For  fixed $L$, direct calculations lead to
\[
\begin{array}{llllll}
&(U_2-\sigma_0U_1)_t-D(U_2-\sigma_0U_1)_{xx}+\alpha(U_2-\sigma_0U_1)_x\\
&=\gamma a(t)(U_2-\sigma_0U_1)-\gamma(U_2^2-\sigma_0U_1^2)\\
&\geq \gamma[a(t)-(U_2+\sigma_0U_1)](U_2-\sigma_0U_1)
\end{array}
\]
for $t\in(0^+,\tau]$ and $x\in(0,L)$, and
\[
(U_2-\sigma_0U_1)(0^+,x)=G(U_2(0,x))-\sigma_0G(U_1(0,x))\geq G(\sigma_0U_1(0,x))-\sigma_0G(U_1(0,x))\geq0
\]
for $x\in(0,L)$ after using the assumption $(\mathcal{H}2)$ of impulsive function $G(u)$. Boundary conditions satisfy $(U_2-\sigma_0U_1)_x(t,0)=0$ and $(U_2-\sigma_0U_1)(t,L)\geq0$.
We have $(U_2-\sigma_0U_1)(t,x)\geq0$ for $t\in[0,\tau]$ and $x\in[0,L]$ and $\min_{t\in[0,\tau]}(U_2-\sigma_0U_1)=0$ for some $(t_0,x_0)\in[0,\tau]\times[0,L]$. If $\min_{t\in[0,\tau]}(U_2-\sigma_0U_1)>0$, then a positive constant $\epsilon$ exists such that $U_2-\sigma_0U_1>\epsilon U_1$, which contradicts $\sigma_0=\sup\Lambda$.

If $x_0= 0$, then  it follows from Hopf's boundary lemma that $(U_2-\sigma_0U_1)_x(t,0)>0$, which is impossible. If $x_0=L$, we have $(U_2-\sigma_0U_1)(t,L)=0$, which together with Hopf's boundary lemma gives $(U_2-\sigma_0U_1)_x(t,L)<0$. In fact, $(U_2-\sigma_0U_1)_x(t,L)=0$. Thus, $x_0\in(0,L)$. Then, strong maximum principle implies $U_2-\sigma_0U_1\equiv0$ for $(t,x)\in[0,\tau]\times(0,L)$. Recalling the equation, we have $\gamma(U_2^2-\sigma_0U_1^2)=0$, which deduces that $\sigma_0=1$ and contradicts $\sigma_0<1$.  $1\in\Lambda$ and we have $U_1\leq U_2$ for $(t,x)\in[0,\tau]\times(0,L)$. Therefore, problem \eqref{b041} with $l=\infty$ admits a unique positive periodic solution. Noticing that if $\lambda_1(D,\alpha,\gamma a(t),G'(0),(0,\infty))< 0$, the problem \eqref{b04} admits a unique positive periodic solution $U^*(t)$, and $U^*(t)$ satisfies the problem \eqref{b041} with $l=\infty$, we can derive that $\lim\limits_{l\rightarrow\infty}U_l(t,x)=U^*(t)$ in $C_{loc}([0,\tau]\times[0,\infty))$.
\end{pf}

In the sequel, we study the long time behaviors of the solution to problem \eqref{b01} and  discuss the spreading or vanishing of species.

\begin{thm}
Assume $\alpha\leq0$. Let $u(t,x)$ be the solution of the problem \eqref{b01}. If the principal eigenvalue $\lambda_1$ of the problem \eqref{b03} satisfies $\lambda_1(D,\alpha,\gamma a(t),G'(0),(0,\infty))\geq 0$,  that is, $\frac{\alpha^2}{4D}-\frac{1}{\tau}\ln G'(0)-\gamma \overline{a}\geq0$, then $r_\infty:=\lim\limits_{t\rightarrow\infty}r(t)<\infty$ and $\lim\limits_{t\rightarrow\infty}\|u(t,\cdot)\|_{C([0,r(t)])}=0$.
\end{thm}
\begin{pf}
This theorem can be proved by Theorem 3.2 in \cite{menglin22} with some modifications. We give the sketch to show the effects of advection, environmental heterogeneity and pulses. An upper solution to the problem \eqref{b01} can be constructed as follows
\[
\widetilde{u}(t,x)=C_0e^{-\nu t}f(t)\varphi(\frac{\tilde{r}x}{\xi(t)}), \ t\in((n\tau)^+,(n+1)\tau],
\]
\[
\omega(t)=1+\eta-\frac{\eta}{2}e^{-\nu t},\ \xi(t)=\widetilde{r}\omega(t), \ t>0,
\]
where $\varphi(x)$ is the eigenfunction of the problem \eqref{b04}  corresponding to $\lambda^*(D,\alpha,(0,\widetilde{r}))$, $f(t)$ satisfies
\[
\left\{
\begin{array}{lll}
f_t=[-\lambda^*(D,\alpha,(0,\widetilde{r}))+\lambda_1(D,\alpha,\gamma a(t), G'(0),(0,\widetilde{r}))+\gamma a(t)]f,\; &\ t\in(0^+,\tau], \\[2mm]
f(0^+)=G'(0)f(0),\
f(0)=f(\tau),
\end{array} \right.
\]
$C_0$ is a sufficiently large constant such that $\widetilde{u}(0,x)=C_0f(0)\varphi(x/(1+\eta/2))\geq u_0(x)$ for $x\in[0,r_0]$, $\widetilde{r}$ can be selected such that $\xi'(t)\geq -\mu_1\widetilde{u}_x(t,\xi(t))$ for $t\in((n\tau)^+,(n+1)\tau]$, $\eta$ and $\nu$ are sufficiently small constants. We only verify the equation and the impulse condition here.

By virtue of Theorem 2.5 and $\lambda_1(D,\alpha,\gamma a(t),G'(0),(0,\infty))\geq 0$, we have $\lambda_1(D,\alpha,\gamma a(t),$\\$G'(0),(0,\widetilde{r}))>0$. It follows from Hopf's boundary lemma that $\phi(t,0)>0$ and $\phi_x(t,\widetilde{r})<0$ for $t\in [0,\tau]$, then there exists a positive constant $c_1$ such that $\phi_x(t,x)\leq c_1\phi(t,x)$ for $t\in[0,\tau]$ and $x\in[0,\widetilde{r}]$. Then, since $\alpha\leq0$, some calculations lead to
 \[
\begin{array}{llllll}
&\widetilde{u}_t-D\widetilde{u}_{xx}+\alpha\widetilde{u}_x
-\gamma\widetilde{u}(a(t)-\widetilde{u}) \\
&\geq C_0 e^{-\nu t}(-\nu f\varphi+f_t\varphi-xf\varphi_x\omega'\omega^{-2}-Df\varphi_{xx}\omega^{-2}
+\alpha f\varphi_x\omega^{-1}-\gamma a(t)f\varphi)\\
& =C_0e^{-\nu t}[-\nu f\varphi+(-\lambda^*+\lambda_1)f\varphi-xf\varphi_x\omega'\omega^{-2}
+f\omega^{-2}(\lambda^*\varphi-\alpha\varphi_x)
+\alpha f\varphi_x\omega^{-1}]\\
&\geq C_0e^{-\nu t}f\varphi\{-\nu+\lambda_1+\omega^{-2}[(\omega-1)\alpha c_1
-\frac{\eta\widetilde{r}c_1\nu}{2}]+(\omega^{-2}-1)\lambda^*\},
\end{array}
\]
which together with $\lim \limits_{\eta\rightarrow0}\omega=1$ yield
\[
\begin{array}{llllll}
\widetilde{u}_t-D\widetilde{u}_{xx}+\alpha\widetilde{u}_x
-\gamma\widetilde{u}(a(t))-\widetilde{u})
\geq C_0e^{-\nu t}f\varphi(-\nu+\lambda_1-\frac{\eta\widetilde{r}c_1\nu}{2})
\geq0,
\end{array}
\]
provided that $\nu$ is sufficiently small.
And, the assumption $(\mathcal{H}2)$ gives
\[
\begin{array}{llllll}
\widetilde{u}((n\tau)^+,x)&=C_0 e^{-\nu n\tau}f((n\tau)^+)\varphi(x)=G'(0)C_0 e^{-\nu n\tau}f(n\tau)\varphi(x)\\
&=G'(0)\widetilde{u}(n\tau,x)\geq G(\widetilde{u}(n\tau,x)).

\end{array}
\]
Therefore,  $\widetilde{u}(t,x)$ is an upper solution to the problem \eqref{b01} and
\[
r_\infty\leq\lim\limits_{t\rightarrow\infty}\xi(t)= \widetilde{r}(1+\eta)<\infty,\ \lim\limits_{t\rightarrow\infty}\|u(t,\cdot)\|_{C([0,r(t)])}\leq \lim\limits_{t\rightarrow\infty}\|\widetilde{u}(t,\cdot)\|_{C([0,r(t)])}=0.
\]
\end{pf}

\begin{thm}
Let $u(t,x)$ be the solution of the problem \eqref{b01}. If the principal eigenvalue $\lambda_1$ of the problem \eqref{b03} satisfies $\lambda_1(D,\alpha,\gamma a(t),G'(0),(0,\infty))<0$, that is, $\frac{\alpha^2}{4D}-\frac{1}{\tau}\ln G'(0)-\gamma \overline{a}<0$, then:

$(i)$ either $r_\infty \leq r_\star<\infty$ and $\lim\limits_{t\rightarrow\infty}\|u(t,\cdot)\|_{C([0,r(t)])}=0$ or $r_\infty=\infty$ and $\lim\limits_{m\rightarrow\infty}u(t+m\tau,x)=U^*(t)$ is locally uniformly for $t\in [0,\tau]$ and $x\in[0,\infty)$, where $r_\star$ satisfies $\lambda_1(D,\alpha,\gamma a(t),G'(0),$$(0,r_\star))$$=0$ and $U^*(t)$ is the unique positive periodic solution to problem \eqref{b05};

$(ii)$ there exists a $\mu_1^*>0$ such that spreading occurs when $\mu_1>\mu_1^*$ and vanishing occurs when $0<\mu_1\leq\mu_1^*$ provided that $r_0<r_\star$.
\end{thm}
\begin{pf}
 Assertion $(i)$ can be verified using comparison principle and the dynamic results for fixed boundary problem \eqref{b03} in Lemma 2.6, which is similar as the proofs in Theorems 3.3-3.4 in \cite{menglin22}.

$(ii)$  If $r_0<r_\star$, the results for the case $\mu_1>\mu_1^*$ can be proved by the proofs in  Theorems  4.5 in \cite{menglin22} with minor modifications, and are omitted here. We aim to construct a suitable upper solution to problem \eqref{b01} to prove that vanishing happens for $\mu_1\leq \mu_1^*$. Some difficulties are induced by the introduction of advection and pulses. The different construction of upper solution are chosen according to the situation $\alpha\leq0$ or $\alpha>0$.

$(a)$ $\alpha\leq0$. The upper solution can be selected as in Theorem 2.7, where $\widetilde{r}=r_0$ and other notations are unchanged.  We can find a $\mu_1^\vartriangle$ such that for $\mu_1\leq\mu_1^\vartriangle$, the following inequality
\[
-\mu_1\widetilde{u}_x(t,\xi(t))=-\mu_1\omega^{-1}C_0e^{-\nu t}f\varphi_x(r_0)\leq\frac{r_0\eta\nu}{2}e^{-\nu t}
\]
holds. Then, $\widetilde{u}$ is an upper solution and
vanishing happens for $\mu_1\leq\mu_1^\vartriangle$.

$(b)$ $\alpha>0$. The construction of the upper solution is inspired by the  proofs in Theorem 3.1 $(ii)$ of \cite{zhangqywang} and modified due to the involution of pulses.

By virtue of Theorem 2.5 and $r_0<r_\star$, we have $\lambda_1(D,\alpha,\gamma a(t),G'(0),(0,r_0))>0$, which together with Lemma 2.4 and the explicit expression of $\lambda_1$ yield $\lambda^*(D,\alpha, (0,r_0))>\alpha^2/(4D)$ and $\lambda^*(D,\alpha, (0,r_0))>\ln G'(0)/\tau+\gamma \overline{a}$.
Denote
\[
A=\frac{\alpha}{2D}, B=\frac{\sqrt{4D\lambda^*({D,\alpha, (0,r_\star))-\alpha^2}}}{2D}, \theta=\arctan\frac{B}{A}
\]
and let
\[
\psi(z)=B\cos z-A\sin z, \  z\in[0,\theta].
\]
Then, we have $0<\theta<\pi/2$, $r_\star=\theta/B$, and $\psi(z)>0$ for $z\in[0,\theta)$, $\psi(\theta)=0$ and $\psi'(z)<0$ for $z\in[0,\theta]$.
Let
\[
\widetilde{u}(t,x)=C_0g(t,x)f(t)\psi(\xi(t,x)), \ t\in((n\tau)^+,(n+1)\tau], x\in(0,\omega(t)),
\]
\[
\omega(t)=r_0(1+\eta-\frac{\eta}{2}e^{-\eta t}),\ g(t,x)=e^{-\eta t}e^{-A(r_0-\frac{r_0(1+\eta)x}{\omega(t)})},\  \xi(t,x)=\frac{\theta x}{\omega(t)}, \ t>0,
\]
where $0<\eta\ll1$, $f(t)$ is the positive solution to
  \[
\left\{
\begin{array}{lll}
f_t=(-\frac{1}{\tau}\ln G'(0)-\gamma \overline{a}+\gamma a(t))f,\; &\ t\in(0^+,\tau], \\[2mm]
f(0^+)=G'(0)f(0),\
f(0)=f(\tau),
\end{array} \right.
\]
 and $C_0$ is large enough such that
$\widetilde{u}(0,x)=C_0f(0)e^{-A(r_0-\frac{(1+\eta)x}{1+\eta/2})}\psi(\xi(0,x))
\geq u_0(x)$ for $t\in[0, r_0]$.

For $t\in((n\tau)^+,(n+1)\tau]$, direct calculations give
 \[
\begin{array}{llllll}
\widetilde{u}_t-D\widetilde{u}_{xx}+\alpha\widetilde{u}_x
-\gamma\widetilde{u}(a(t)-\widetilde{u})
\geq \widetilde{u}_t-D\widetilde{u}_{xx}+\alpha\widetilde{u}_x
-\gamma a(t)\widetilde{u}
:=q_1(t,x)+q_2(t,x),
\end{array}
\]
where
  \[
\begin{array}{llllll}
&q_1(t,x)=\widetilde{u}(-\eta-\frac{Ar_0(1+\eta)x\omega'(t)}{\omega^2(t)}
+\frac{A\alpha r_0(1+\eta)(2\omega(t)-r_0(1+\eta))}{2\omega^2(t)}-
\gamma\overline{a}-\frac{1}{\tau}\ln G'(0)+\frac{D\theta^2}{\omega^2(t)}),\\
&q_2(t,x)=-C_0\theta g(t,x)f(t)\psi'(\xi(t,x))(\frac{x\omega'(t)}{\omega^2(t)}+\frac{\alpha r_0(1+\eta)}{\omega^2(t)}-\frac{\alpha}{\omega(t)}).
\end{array}
\]
Recall that
\[
\frac{x\omega'(t)}{\omega^2(t)}\leq\frac{\eta^2}{2}, 0<\frac{x}{\omega(t)}<1, r_0(1+\eta/2)<\omega(t)<r_0(1+\eta),
\]
whence follows that $\frac{2\omega(t)-r_0(1+\eta)}{\omega^2(t)}\geq\frac{1}{r_0(1+\eta)^2}>0$. Then, we have
\[
\begin{array}{llllll}
q_1(t,x)&\geq \widetilde{u}(-\eta-\frac{Ar_0(1+\eta)\eta^2}{2}
+\frac{A\alpha r_0(1+\eta)(2\omega(t)-r_0(1+\eta))}{2\omega^2(t)}-
\gamma\overline{a}-\frac{1}{\tau}\ln G'(0)+\frac{D\theta^2}{\omega^2(t)})\\
& >\widetilde{u}(-\eta-\frac{Ar_0(1+\eta)\eta^2}{2}
+\frac{A\alpha}{2(1+\eta)}-
\gamma\overline{a}-\lambda^*(D,\alpha, (0,r_\star))+\frac{D\theta^2}{r_0^2(1+\eta)^2})\\
&:=\widetilde{u}Q(\eta)
\end{array}
\]
due to $\lambda^*(D,\alpha, (0,r_\star))>\ln G'(0)/\tau+\gamma \overline{a}$. $r_0<\theta/B$ gives
\begin{small}
\[
Q(0)=\frac{A\alpha}{2}+\frac{D\theta^2}{r_0^2}-\lambda^*(D,\alpha, (0,r_\star))=\frac{\alpha^2-4D\lambda^*(D,\alpha, (0,r_\star))}{4D}+\frac{D\theta^2}{r_0^2}=-DB^2+\frac{D\theta^2}{r_0^2}>0,
\]
\end{small}
which indicates that $q_1(t,x)>0$ if $0<\eta\ll1$. Meanwhile, $\psi'(\xi(t,x))<0$ and $\omega'(t)>0$ give $q_2(t,x)>0$. Then, we have $\widetilde{u}_t-D\widetilde{u}_{xx}+\alpha\widetilde{u}_x
-\gamma\widetilde{u}(a(t)-\widetilde{u})>0$ for $t\in ((n\tau)^+,(n+1)\tau]$ and $x\in(0,\omega(t))$.

For $t=(n\tau)^+$, it follows from the assumption $(\mathcal{H}2)$ that
\[
\begin{array}{llllll}
\widetilde{u}((n\tau)^+,x)&=C_0 g(n\tau,x)f((n\tau)^+)\psi(\xi(n\tau,x))=G'(0)C_0 g(n\tau,x)f(n\tau)\psi(\xi(n\tau,x)) \\
&=G'(0)\widetilde{u}(n\tau,x)\geq G(\widetilde{u}(n\tau,x)).

\end{array}
\]

Notice that $\omega(0)=r_0(1+\eta/2)>r_0$ and   $\widetilde{u}(t,\omega(t))=C_0e^{-\eta t+A\eta r_0}f(t)\psi(\theta)=0$. Thanks to $Br_0<\theta$, we have
 $\widetilde{u}_x(t,0)=C_0e^{-(Ar_0+\eta t)} f(t)\frac{A[Br_0(1+\eta)-\theta]}{\omega(t)}<0$ provided that $0<\eta\ll1$.
 And, in order to satisfy $\omega'(t)\geq-\mu_1\widetilde{u}_x(t,\omega(t))$, by calculations we have
 \[
\begin{array}{llllll}
-\mu_1\widetilde{u}_x(t,\omega(t))&=\frac{\mu_1 C_0\theta g(t,\omega(t))f(t)\psi'(\theta)}{\omega(t)}=\frac{\mu_1 C_0\theta }{\omega(t)}f(t)e^{-\eta t+Ar_0\eta}(B\sin \theta+A \cos \theta)\\
&\leq \frac{\mu_1 C_0\theta }{r_0}f(t)e^{-\eta t+Ar_0\eta}(B\sin \theta+A \cos \theta)\leq \frac{1}{2}r_0\eta^2 e^{-\eta t}=\omega'(t),
\end{array}
\]
whence we can find a $\mu^\vartriangle_1$ such that the above inequality holds for $\mu_1\leq \mu^\vartriangle_1$.
Therefore, $\widetilde{u}(t,x)$ is an upper solution to the problem \eqref{b01},
$
r_\infty\leq\lim\limits_{t\rightarrow\infty}\omega(t)= r_0(1+\eta)<\infty
$, and vanishing happens for $0<\mu_1\leq \mu^\vartriangle_1$.
\end{pf}

To investigate the sufficient conditions about pulses and expanding capability for spreading and vanishing of two competitors in the later section, the dynamics of  one single species free boundary problem \eqref{b01} are divided according to different pulses.

\begin{thm}
For initial domain $(0,r_0)$ and whole area $(0,\infty)$, denote $g_\star:=G'_1(0)$ and $g^\star:=G'_2(0)$ with impulsive functions $G_1$ and $G_2$, which satisfy \[
\lambda_1(D,\alpha,\gamma a(t),g_\star,(0,\infty))=0 \ \textrm{and} \  \lambda_1(D,\alpha,\gamma a(t),g^\star,(0,r_0))=0, \]respectively.

$(i)$ If $G'(0)\geq g^\star$, then spreading occurs for any $\mu_1>0$.

$(ii)$ If $g_\star<G'(0)<g^\star$, then there exists a $\mu_1^*>0$ such that spreading occurs when $\mu_1>\mu_1^*$ and vanishing occurs when $0<\mu_1\leq\mu_1^*$.

$(iii)$ If  $\alpha\leq0$ and $G'(0)\leq g_\star$, then vanishing occurs for any $\mu_1>0$.
\end{thm}
\begin{pf}
If $G'(0)\geq g^\star$, Theorem 2.5 gives that $\lambda_1(D,\alpha,\gamma a(t),G'(0),(0,r_0))\leq0$. Noting that $\lambda_1(D,\alpha,\gamma a(t),G'(0),(0,l))$ is decreasing with respect to $l$, we have $\lambda_1(D,\alpha,\gamma a(t),$\\$G'(0),(0,\infty))<0$. By contradiction, if $r_\infty<\infty$ and vanishing happens, it follows from Theorem 2.8 $(i)$ that $r_\infty\leq r_\star$, which deduces that $\lambda_1(D,\alpha,\gamma a(t), G'(0),(0,\infty))\geq 0$, that leads to a contradiction.

  If $g_\star<G'(0)<g^\star$, it follows from Theorem 2.5 that $\lambda_1(D,\alpha,\gamma a(t),G'(0),(0,\infty))<0$ and $\lambda_1(D,\alpha,\gamma a(t),G'(0),(0,r_0))>0$. Similarly as  the proofs in Theorem 2.8 $(ii)$, the assertion $(ii)$ can be obtained and are omitted.

  If $G'(0)\leq g_\star$, we have $\lambda_1(D,\alpha,\gamma a(t),G'(0),(0,\infty))\geq0$. $\alpha\leq0$ and Theorem 2.7 implies that $\lim\limits_{t\rightarrow\infty}\|u(t,\cdot)\|_{C([0,r(t)])}=0$ and vanishing occurs.
\end{pf}

\section{Classifications of spreading and vanishing}
Recalling that $r'(t)>0$ and $s'(t)>0$ in Theorem 2.1,  there exist $r_\infty\in(r_0,\infty]$ and $s_\infty\in(s_0,\infty]$ such that
\[
\lim\limits_{t\rightarrow\infty}r(t)=r_\infty \ \textrm{and} \ \lim\limits_{t\rightarrow\infty}s(t)=s_\infty.
\]

This section deals with the criteria governing the spreading or vanishing modeled by problem \eqref{a01}.

It follows from Theorem 2.7 that if $\alpha\leq0$ and $\alpha^2/(4D)-\ln G'(0)/\tau-\gamma\geq0$, then $\|u(t,\cdot)\|_{C([0,r(t)])}$ tends to 0 as $t\rightarrow\infty$ and problem \eqref{a01} can be regarded as one single species free boundary problem for $v$. And if $\beta\leq0$ and $\beta^2/4-\ln H'(0)/\tau-1\geq0$, then $\lim\limits_{t\rightarrow\infty}\|v(t,\cdot)\|_{C([0,s(t)])}=0$ and the problem \eqref{a01} can be regarded as one single species free boundary problem for $u$. Here, $\alpha^2/(4D)-\ln G'(0)/\tau-\gamma$ and $\beta^2/4-\ln H'(0)/\tau-1$ are the explicit expressions of the $\lambda_1(D,\alpha,\gamma (1+\varepsilon_1(t)), G'(0),(0,\infty))$ and $\lambda_1(1,\beta, 1+\varepsilon_2(t), H'(0),(0,\infty))$, respectively. Combined with Theorems 2.7  and 2.8, then the dynamics of two competing species modeled by the problem \eqref{a01} can be divided into the following situations.
\begin{thm}
The following assertions hold.

$(i)$ If $\alpha,\beta\leq0$, $\alpha^2/(4D)-\ln G'(0)/\tau-\gamma\geq0$ and $\beta^2/4-\ln H'(0)/\tau-1\geq0$, then both species $u$ and $v$ vanish eventually;

$(ii)$ if $\alpha^2/(4D)-\ln G'(0)/\tau-\gamma<0$, $\beta\leq0$ and $\beta^2/4-\ln H'(0)/\tau-1\geq0$, then co-extinction occurs for $r_\infty\leq r_*<\infty$,
and species $u$ spreads and $v$ vanishes for $r_\infty=\infty$, where $r_*$  satisfies $\lambda_1(D,\alpha,1+\varepsilon_1(t), G'(0),(0,r_*))=0$;

$(iii)$  if $\alpha\leq0$, $\alpha^2/(4D)-\ln G'(0)/\tau-\gamma\geq0$ and $\beta^2/4-\ln H'(0)/\tau-1<0$, then co-extinction occurs for $s_\infty\leq s_*<\infty$, and species $u$ vanishes and $v$ spreads for $s_\infty=\infty$,
where $s_*$ satisfies $\lambda_1(1,\beta,1+\varepsilon_2(t), H'(0),(0,s_*))=0$.
\end{thm}

Naturally, we next consider the case
\[
(\mathcal{A}1)\ \alpha^2/(4D)-\ln G'(0)/\tau-\gamma<0, \ \beta^2/4-\ln H'(0)/\tau-1<0.
\]
It follows from Lemma 2.6 $(ii)$ that the problem \eqref{b05}  with $(\gamma, a(t),G(u))$  replaced with $(\gamma, 1+\varepsilon_1(t), G(u))$ and $(1, 1+\varepsilon_2(t), H(v))$ admits a unique positive periodic solution $U^*(t)$ and $V^*(t)$, respectively, where $U^*(t)$ and $V^*(t)$ satisfy
\begin{eqnarray}
\left\{
\begin{array}{lll}
U^*_t=\gamma U^*(1+\varepsilon_1(t)-U^*),\; &\  t\in(0^+,\tau], \\[2mm]
U^*(0^+)=G(U^*(0)),\ U^*(0)=U^*(\tau),
\end{array} \right.
\label{c011}
\end{eqnarray}
and
\begin{eqnarray}
\left\{
\begin{array}{lll}
V^*_t=V^*(1+\varepsilon_2(t)-V^*),\; &\  t\in(0^+,\tau], \\[2mm]
V^*(0^+)=H(V^*(0)), \ V^*(0)=V^*(\tau),
\end{array} \right.
\label{c01}
\end{eqnarray}
 respectively.
 Then, in order to give the sufficient conditions about habitat sizes for the coexistence, see Theorem 3.7 $(ii)$, we  assume $\lambda_1(D,\alpha,\gamma (1+\varepsilon_1(t)-kV^*),G'(0),(0,\infty))<0$ and $\lambda_1(1,\beta,1+\varepsilon_2(t)-hU^*, H'(0),(0,\infty))<0$.

Since $V^*(t)$ satisfies problem \eqref{c01},
dividing the first equation by $V^*$ and integrating from $0^+$ to $\tau$, we obtain  $\overline{V^*}:=\frac{1}{\tau}\int_{0^+}^\tau V^*(t)dt=1+\frac{1}{\tau}\ln\frac{H(V^*(0))}{V^*(0)}$. Recalling Theorem 2.5 and $(\mathcal{H}2)$, we have
\begin{small}
\[
\begin{array}{llllll}
&\lambda_1(D,\alpha,\gamma (1+\varepsilon_1(t)-kV^*),G'(0), (0,\infty))=\frac{\alpha^2}{4D}-\frac{1}{\tau}\ln G'(0)
 -\gamma+k\gamma\overline{V^*}\\
 &=\frac{\alpha^2}{4D}-\frac{1}{\tau}\ln G'(0)
 +\gamma(k-1)+\frac{k\gamma}{\tau}\ln\frac{H(V^*(0))}{V^*(0)}
 \leq\frac{\alpha^2}{4D}-\frac{1}{\tau}\ln G'(0)
 +\gamma(k-1)+\frac{k\gamma}{\tau}\ln H'(0).
\end{array}
\]
\end{small}
And \[
\lambda_1(1,\beta,1+\varepsilon_2(t)-hU^*,H'(0), (0,\infty))\leq \beta^2/4-\ln H'(0)/\tau+(h-1)+h\ln G'(0)/(\gamma\tau)
\]
can be obtained similarly.
Therefore, we give the following assumption
\[
\begin{array}{llllll}
(\mathcal{A}2)\ &\alpha^2/4D-\ln G'(0)/\tau
 +\gamma(k-1)+k\gamma\ln H'(0)/\tau<0, \\
  &\beta^2/4-\ln H'(0)/\tau+(h-1)+h\ln G'(0)/(\gamma\tau)<0.
\end{array}
\]

In the sequel, we assume that assumptions $(\mathcal{A}1)$ and $(\mathcal{A}2)$ hold unless otherwise  specified.

\begin{lem} If $r_\infty\leq r_*$, then $\lim\limits_{t\rightarrow\infty}\|u(t,\cdot)\|_{C([0,r(t)])}=0$,
where $r_*$ satisfies $\lambda_1(D,\alpha,\gamma (1+\varepsilon_1(t)), G'(0),(0,r_*))=0$.
\end{lem}
\begin{pf}
For any given $l\in[r_\infty,r_*]$, let $\overline{u}(t,x)$ be the solution to the problem
\[
\left\{
\begin{array}{lll}
\overline{u}_t=D\overline{u}_{xx}-\alpha \overline{u}_x+\gamma \overline{u}(1+\varepsilon_1(t)-\overline{u}),\; &\ t\in((n\tau)^+,(n+1)\tau],\ x\in(0,l), \\[2mm]
\overline{u}((n\tau)^+,x)=G(\overline{u}(n\tau,x)),\; &\ x\in(0,l),\\[2mm]
\overline{u}_x(t,0)=\overline{u}(t,l)=0,\; &\ t\in(0,\infty),\\[2mm]
\overline{u}(0,x)=M_1,\; &\ x\in[0,l],
\end{array} \right.
\]
where $M_1$, defined in Theorem 2.1, is the upper bound of $u(t,x)$, which together with comparison principle  yield that $u(t,x)\leq \overline{u}(t,x)$ for $t\in[0,\infty)$ and $x\in[0,r(t)]$. Recalling that $l\in[r_\infty,r_*]$ and Theorem 2.5, $\lambda_1(D,\alpha,\gamma (1+\varepsilon_1(t)),G'(0),(0,l))\geq0$. Then, Lemma 2.6 $(i)$ gives that $\lim\limits_{t\rightarrow\infty}\|\overline{u}(t,\cdot)\|_{C([0,l])}=0$, and therefore $\lim\limits_{t\rightarrow\infty}\|u(t,\cdot)\|_{C([0,r(t)])}=0$.
\end{pf}

\begin{thm}
 If $r_\infty=\infty$ and $s_\infty\leq s_*$, then the solution $(u,v)$ to the problem \eqref{a01} satisfies $\lim\limits_{m\rightarrow\infty}(u,v)(t+m\tau,x)=(U^*(t),0)$ locally uniformly in $[0,\tau]\times[0,\infty)$, where $s_*$ satisfies $\lambda_1(1,\beta,1+\varepsilon_2(t), H'(0),(0,s_*))=0$.
\end{thm}
\begin{pf}
 It follows from assumption $(\mathcal{A}1)$ that $\lambda_1(D,\alpha,\gamma (1+\varepsilon_1(t)),G'(0),(0,\infty))<0$. There exists a sufficiently small $\delta>0$ such that $\lambda_1(D,\alpha,\gamma (1+\varepsilon_1(t)-k\delta),G'(0),(0,\infty))<0$. Since $r_\infty=\infty$, there exists a large integer $N$ such that $r(N\tau)>r_*^\delta$, where $r_*^\delta$ satisfies $\lambda_1(D,\alpha,\gamma (1+\varepsilon_1(t)-k\delta),G'(0),(0,r_*^\delta))=0$. Similar as Lemma 3.2, $s_\infty\leq s_*$ gives $\lim\limits_{t\rightarrow\infty}\|v(t,\cdot)\|_{C([0,s(t)])}=0$, which deduces that there exists a large integer $N_\delta\geq N$ such that $0<v(t,x)\leq \delta$ for $t\geq N_\delta\tau$ and $0\leq x \leq s(t)$. For $n\geq N_\delta$, consider problem
\begin{small}
\begin{eqnarray}
\left\{
\begin{array}{lll}
\underline{u}_t=D\underline{u}_{xx}-\alpha \underline{u}_x+\gamma \underline{u}(1+\varepsilon_1(t)-k\delta-\underline{u}),\; &\ t\in((n\tau)^+,(n+1)\tau],\ x\in(0,\underline{r}(t)), \\[2mm]
\underline{u}((n\tau)^+,x)=G(\underline{u}(n\tau,x)),\; &\ x\in(0,\underline{r}(n\tau)),\\[2mm]
\underline{u}_x(t,0)=\underline{u}(t,\underline{r}(t))=0,\; &\ t\in(0,\infty),\\[2mm]
\underline{r}'(t)=-\mu_1\underline{u}_x(t,\underline{r}(t)),\; &\ t\in(n\tau,(n+1)\tau],\\[2mm]
\underline{r}(N_\delta\tau)=r(N_\delta\tau),\underline{u}(N_\delta\tau,x)
=u(N_\delta\tau,x),\; &\ x\in[0,\underline{r}(N_\delta\tau)].
\end{array} \right.
\label{c02}
\end{eqnarray}
\end{small}
By virtue of comparison principle,
\begin{equation}
u(t,x)\geq \underline{u}(t,x),\ r(t)\geq \underline{r}(t)\ \textrm{for} \ t\geq N_\delta\tau\ \textrm{and} \ 0\leq x\leq \underline{r}(t),
\label{c03}
\end{equation}
where $\underline{u}(t,x)$ is the solution of the problem \eqref{c02}. Recalling $\lambda_1(D,\alpha,\gamma (1+\varepsilon_1(t)-k\delta),G'(0),(0,\infty))<0$ and $\underline{r}(N_\delta\tau)=r(N_\delta\tau)>r_*^\delta$, it follows from Theorem 2.8 $(i)$ that $\lim\limits_{m\rightarrow\infty}\underline{u}(t+m\tau,x)=U_\delta^*(t)$ is locally uniformly in $[0,\tau]\times[0,\infty)$, where $U_\delta^*(t)$ satisfies the problem \eqref{b05} with  $ a(t)$ replaced with $ 1+\varepsilon_1(t)-k\delta$, which combined with \eqref{c03} yield
\begin{equation}
\liminf\limits_{m\rightarrow\infty}u(t+m\tau,x)\geq U_\delta^*(t)
\label{c04}
\end{equation}
locally uniformly in $[0,\tau]\times[0,\infty)$.

Let $\overline{u}(t,x)$ be the solution of the problem
\begin{small}
\begin{eqnarray}
\left\{
\begin{array}{lll}
\overline{u}_t=D\overline{u}_{xx}-\alpha \overline{u}_x+\gamma \overline{u}(1+\varepsilon_1(t)-\overline{u}),\; &\ t\in((n\tau)^+,(n+1)\tau],\ x\in(0,\overline{r}(t)), \\[2mm]
\overline{u}((n\tau)^+,x)=G(\overline{u}(n\tau,x)),\; &\ x\in(0,\overline{r}(n\tau)),\\[2mm]
\overline{u}_x(t,0)=\overline{u}(t,\overline{r}(t))=0,\; &\ t\in(0,\infty),\\[2mm]
\overline{r}'(t)=-\mu_1\overline{u}_x(t,\overline{r}(t)),\; &\ t\in(n\tau,(n+1)\tau],\\[2mm]
\overline{r}(0)=r_0,\overline{u}(0,x)
=\max\limits_{x\in[0,r_0]}u_0(x),\; &\ x\in[0,r_0],
\end{array} \right.
\label{c041}
\end{eqnarray}
\end{small}
which by comparison principle leads to
\begin{equation}
\overline{u}(t,x)\geq u(t,x),\ \overline{r}(t)\geq r(t)\ \textrm{for} \ t\geq 0\ \textrm{and} \ 0\leq x\leq r(t).
\label{c05}
\end{equation}
As $r_\infty=\infty$, then $\overline{r}_\infty=\infty$. Theorem 2.8 $(i)$ and \eqref{c05} give that
\[
\limsup \limits_{m\rightarrow\infty}u(t+m\tau,x)\leq \limsup \limits_{m\rightarrow\infty}\overline{u}(t+m\tau,x)= U^*(t)
\]
locally uniformly in $[0,\tau]\times[0,\infty)$. Recalling \eqref{c04} and by the continuous dependence of $U_\delta^*(t)$ on $\delta$, we have $\lim\limits_{\delta\rightarrow0}U_\delta^*(t)=U^*(t)$ for $t\in[0,\tau]$. The proof is completed.
\end{pf}

The coexistence of two competing species can occur, which is shown in the following theorem.
\begin{thm}
If $r_\infty=s_\infty=\infty$, then
\[
U_*(t)\leq\liminf \limits_{m\rightarrow\infty}u(t+m\tau,x)\leq \limsup \limits_{m\rightarrow\infty}u(t+m\tau,x)\leq U^*(t),
\]
\[
V_*(t)\leq\liminf \limits_{m\rightarrow\infty}v(t+m\tau,x)\leq \limsup \limits_{m\rightarrow\infty}v(t+m\tau,x)\leq V^*(t)
\]
locally uniformly in $[0,\tau]\times[0,\infty)$, where $U^*(t)$ and $V^*(t)$ are defined in \eqref{c011} and \eqref{c01}, $U_*(t)$ and $V_*(t)$ are the unique positive periodic solutions to the problems
\begin{eqnarray}
\left\{
\begin{array}{lll}
(U_*)_t=\gamma U_*(1+\varepsilon_1(t)-U_*-kV^*(t)),\; &\  t\in(0^+,\tau], \\[2mm]
U_*(0^+)=G(U_*(0)), \ U_*(0)=U_*(\tau),
\end{array} \right.
\label{c012}
\end{eqnarray}
and
\begin{eqnarray}
\left\{
\begin{array}{lll}
(V_*)_t=V_*(1+\varepsilon_2(t)-V_*-hU^*(t)),\; &\  t\in(0^+,\tau], \\[2mm]
V_*(0^+)=H(V_*(0)),\ V_*(0)=V_*(\tau),
\end{array} \right.
\label{c013}
\end{eqnarray}
 respectively.
\end{thm}
\begin{pf}
Similarly as above, we can deduce that
\begin{equation}
\limsup \limits_{m\rightarrow\infty}u(t+m\tau,x)\leq U^*(t), \ \limsup \limits_{m\rightarrow\infty}v(t+m\tau,x)\leq V^*(t).
\label{c06}
\end{equation}
The assumption $(\mathcal{A}2)$ gives that $\lambda_1(D,\alpha,\gamma (1+\varepsilon_1(t)-kV^*),G'(0),(0,\infty))<0$, which implies that there exists $\delta_0$ such that $\lambda_1(D,\alpha,\gamma [1+\varepsilon_1(t)-k(V^*+\delta)],G'(0),(0,\infty))<0$ for $0<\delta\leq\delta_0$. For fixed $0<\delta\leq\delta_0$, there exists $l_\delta$ such that $\lambda_1(D,\alpha,\gamma [1+\varepsilon_1(t)-k(V^*+\delta)],G'(0),(0,l))<0$ for $l\geq l_\delta$. Due to \eqref{c06} and $s_\infty=\infty$, for fixed $0<\delta\leq\delta_0$ and $l\geq l_\delta$, there exists a large integer $N_\delta$ such that $s(t)>l$ and $v(t,x)\leq V^*(t)+\delta$ for $t\geq N_\delta\tau$ and $0<x<l$. For $n\geq N_\delta$ and $l\geq l_\delta$, let $\underline{u}(t,x)$ be the solution to
\begin{small}
\[
\left\{
\begin{array}{lll}
\underline{u}_t=D\underline{u}_{xx}-\alpha \underline{u}_x+\gamma \underline{u}[1+\varepsilon_1(t)-k(V^*(t)+\delta)-\underline{u}],\; &\ t\in((n\tau)^+,(n+1)\tau],\ x\in(0,l), \\[2mm]
\underline{u}((n\tau)^+,x)=G(\underline{u}(n\tau,x)),\; &\ x\in(0,l),\\[2mm]
\underline{u}_x(t,0)=\underline{u}(t,l)=0,\; &\ t\in[N_\delta\tau,\infty),\\[2mm]
\underline{u}(N_\delta\tau,x)
=u(N_\delta\tau,x),\; &\ x\in[0,l],
\end{array} \right.
\]
\end{small}
which by comparison principle gives that $u(t,x)\geq\underline{u}(t,x)$ for $t\geq N_\delta\tau$ and $0\leq x \leq l$. Since $\lambda_1(D,\alpha,\gamma [1+\varepsilon_1(t)-k(V^*+\delta)],G'(0),(0,l))<0$, it follows from Lemma 2.6 $(ii)$ that $\lim\limits_{m\rightarrow\infty}\underline{u}(t+m\tau,x)=\underline{U_l}(t,x)$, where $\underline{U_l}(t,x)$ is the unique positive solution to the problem
\begin{small}
\[
\left\{
\begin{array}{lll}
\underline{U_l}_t=D\underline{U_l}_{xx}-\alpha \underline{U_l}_x+\gamma \underline{U_l}[1+\varepsilon_1(t)-k(V^*(t)+\delta)-\underline{U_l}],\; &\ t\in(0^+,\tau],\ x\in(0,l), \\[2mm]
\underline{U_l}((n\tau)^+,x)=G(\underline{U_l}(n\tau,x)),\; &\ x\in(0,l),\\[2mm]
\underline{U_l}_x(t,0)=\underline{U_l}(t,l)=0,\; &\ t\in[0,\tau],\\[2mm]
\underline{U_l}(0,x)
=\underline{U_l}(\tau,x),\; &\ x\in[0,l].
\end{array} \right.
\]
\end{small}
It is known that $\underline{U_l}$ is increasing with respect to $l$ and $\lim\limits_{l\rightarrow\infty}\underline{U_l}(t,x)=U^\delta_*(t)$  is locally uniformly in $[0,\tau]\times[0,\infty)$, where $U^\delta_*(t)$ is the unique positive solution of
\begin{eqnarray}
\left\{
\begin{array}{lll}
(U^\delta_*)_t=\gamma U^\delta_*(1+\varepsilon_1(t)-kV^*(t)-k\delta-U^\delta_*),\; &\  t\in(0^+,\tau], \\[2mm]
U^\delta_*(0^+)=G(U^\delta_*(0)),\\[2mm]
U^\delta_*(0)=U^\delta_*(\tau).
\end{array} \right.
\label{c07}
\end{eqnarray}
As $\delta\rightarrow0$, we have $U^\delta_*(t)\rightarrow U_*(t)$. Since $\delta$ is arbitrary, $\liminf\limits_{m\rightarrow\infty}u(t+m\tau,x)\geq U_*(t)$ locally uniformly in $[0,\tau]\times[0,\infty)$, which together with \eqref{c06} yield that $U_*(t)\leq\liminf \limits_{m\rightarrow\infty}u(t+m\tau,x)\leq \limsup \limits_{m\rightarrow\infty}u(t+m\tau,x)\leq U^*(t)$ uniformly for $t\in[0,\tau]$ and locally uniformly for $x\in[0,\infty)$. The result of $v(t,x)$ can also be similarly obtained.
\end{pf}

We aim to understand what dynamics of species will be as the limit of moving front belongs to different  threshold ranges in the sequel.
\begin{lem}
If $r_\infty>r^*$, then $r_\infty=\infty$ and $\liminf\limits_{m\rightarrow\infty}u(t+m\tau,x) \geq U_*(t)$ locally uniformly in $[0,\tau]\times[0,\infty)$, where $r^*$ satisfies $\lambda_1(D,\alpha,\gamma (1+\varepsilon_1(t)-kV^*),G'(0),(0,r^*))=0$.
\end{lem}
\begin{pf}
Similarly as in the proof of Theorem 3.3, we have $\limsup\limits_{m\rightarrow\infty}v(t+m\tau,x)\leq V^*(t)$ in $C_{loc}([0,\tau]\times[0,\infty))$. For fixed $\delta$,
there exists $r_\delta^*>0$ such that $\lambda_1(D,\alpha,\gamma [1+\varepsilon_1(t)-k(V^*+\delta)],G'(0),(0,r_\delta^*))=0$, and then $r_\delta^*<r^*$ by the properties of $\lambda_1$ in Theorem 2.5. As $r_\infty>r^*$, we have $r_\infty>r_\delta^*$. Then, there exists a large integer $N$ such that $r(N\tau)>r_\delta^*$ and $v(t,x)\leq V^*(t)+\delta$ locally uniformly in $[N\tau,\infty)\times[0,\infty)$. For $n\geq N$, let $(\hat{u}(t,x),\hat{r}(t))$ be the solution to
\begin{small}
\[
\left\{
\begin{array}{lll}
\hat{u}_t=D\hat{u}_{xx}-\alpha \hat{u}_x+\gamma \hat{u}[1+\varepsilon_1(t)-kV^*(t)-k\delta-\hat{u}],\; &\ t\in((n\tau)^+,(n+1)\tau],\ x\in(0,\hat{r}(t)), \\[2mm]
\hat{u}((n\tau)^+,x)=G(\hat{u}(n\tau,x)),\; &\ x\in(0,\hat{r}(n\tau)),\\[2mm]
\hat{u}_x(t,0)=\hat{u}(t,\hat{r}(t))=0,\; &\ t\in[N\tau,\infty),\\[2mm]
\hat{r}'(t)=-\mu_1\hat{u}_x(t,\hat{r}(t)), \hat{r}(N\tau)=r(N\tau)\; &\ t\in(n\tau,(n+1)\tau],\\[2mm]
\hat{u}(N\tau,x)
=u(N\tau,x),\; &\ x\in[0,\hat{r}(N\tau)],
\end{array} \right.
\]
\end{small}
which by comparison principle gives that
\begin{equation}
u(t,x)\geq \hat{u}(t,x) \ \textrm{and}\ r(t)\geq \hat{r}(t)\ \textrm{for}\ t\geq N\tau, 0\leq x \leq \hat{r}(t).
\label{c08}
\end{equation}
Recalling that $\hat{r}(N\tau)\geq r_\delta^*$ and Theorem 2.8 $(i)$, then $\hat{r}_\infty=\infty$ and $\lim\limits_{m\rightarrow\infty}\hat{u}(t+m\tau,x)=U_\delta(t)$ locally uniformly in $[0,\tau]\times[0,\infty)$, where $U_\delta(t)$ is the unique positive solution to problem \eqref{c07}.  \eqref{c08} gives that $r_\infty\geq \hat{r}_\infty=\infty$ and $\lim\limits_{m\rightarrow\infty}u(t+m\tau,x)\geq U_\delta(t)$ locally uniformly in $[0,\tau]\times[0,\infty)$. Since $\delta$ is arbitrary and let $\delta\rightarrow0$, the proof is completed.

\end{pf}

\begin{lem}
If $r_\infty\in (r_*,r^*]$, then $s_\infty=\infty$ and  $\lim\limits_{m\rightarrow\infty}(u,v)(t+m\tau,x)=(0,V^*(t))$ locally uniformly in $[0,\tau]\times[0,\infty)$.
\end{lem}
\begin{pf}
We first claim that $s_\infty>s_*$. By contradiction, suppose $s_\infty\leq s_*$. It can be proved similarly as in Lemma 3.2 that $\lim\limits_{t\rightarrow\infty}\|v(t,\cdot)\|_{C([0,s(t)])}=0$. Then, for sufficiently small $\delta$, there exists $\underline{r}_\delta>0$ such that $\lambda_1(D,\alpha,\gamma (1+\varepsilon_1(t)-k\delta),G'(0),(0,\underline{r}_\delta))=0$ and $\underline{r}_\delta<r_*$, which together with $r_\infty>r_*$ give that $r_\infty>\underline{r}_\delta$ and there exists a large integer $N_\delta$ such that $r(N_\delta\tau)>\underline{r}_\delta$ and $v(t,x)\leq \delta$ for $t\geq N_\delta \tau$ and $0\leq x \leq \underline{r}(t)$. Considering the problem \eqref{c02} for $n\geq N_\delta$, the comparison principle implies that $r(t)\geq \underline{r}(t)$ for $t\geq N_\delta\tau$. By Theorem 2.8 $(i)$, we have $\underline{r}_\infty=\infty$. Then, $r_\infty\geq\underline{r}_\infty=\infty$, which contradicts with $r_\infty\leq r^*$.

Next,  if $\lim\limits_{t\rightarrow\infty}\|u(t,\cdot)\|_{C([0,r(t)])}=0$ is true, then similarly as the proofs in Theorem 3.3, we can prove that $s_\infty=\infty$ and  $\lim\limits_{m\rightarrow\infty}(u,v)(t+m\tau,x)=(0,V^*(t))$ locally uniformly in $[0,\tau]\times[0,\infty)$. In what follows, we aim to prove $\lim\limits_{t\rightarrow\infty}\|u(t,\cdot)\|_{C([0,r(t)])}=0$, which is  inspired by the proofs of Lemma 3.5 in \cite{guowu2015}. If it is not true, suppose  $\delta_0=\limsup\limits_{t\rightarrow\infty}\|u(t,\cdot)\|_{C([0,r(t)])}$\\$>0$.
Then, there exists a sequence $\{(t_i,x_i)\}_{i=1}^\infty\in((n\tau)^+,(n+1)\tau]\times[0,r(t_i)]$ satisfying $t_i\rightarrow\infty$ as $i\rightarrow\infty$ such that $u(t_i,x_i)\geq \delta_0/2$ for $i\in \mathbb{N}$. Due to $0\leq x_i<r_\infty<\infty$, up to  a subsequence, we have $\lim\limits_{i\rightarrow\infty}x_i=x_0\in[0,r_\infty]$. We claim that $x_0\in[0,r_\infty)$. In fact, if $x_0=r_\infty$, then $\lim\limits_{i\rightarrow\infty}x_i-r(t_i)=0$. A direct calculation yields that
\[
\begin{array}{llllll}
|\frac{u(t_i,x_i)}{x_i-r(t_i)}|&=|\frac{u(t_i,x_i)-u(t_i,r(t_i))}{x_i-r(t_i)}|
\\[2mm]
&=|u_x(t_i, \overline{x})|\leq c, \ \overline{x}\in(x_i,r(t_i)),
\end{array}
\]
which together with $|\frac{u(t_i,x_i)}{x_i-r(t_i)}|\geq |\frac{\delta_0}{2(x_i-r(t_i))}|$ gives $|\frac{\delta_0}{2(x_i-r(t_i))}|\leq c$. This is in contradiction with $\lim\limits_{i\rightarrow\infty}x_i-r(t_i)=0$.
Due to $r_\infty<\infty$, define
\[
(w_1,w_2)(y,t):=(u,v)(x,t), \ y:=\frac{x}{r(t)},\ q(t):=\frac{s(t)}{r(t)},
\]
then the problem \eqref{a01} turns to
\begin{small}
\[
\left\{
\begin{array}{lll}
w_{1,t}=\frac{D}{r^2(t)}w_{1,yy}-(\frac{\alpha}{r(t)}-\frac{r'(t)y}{r(t)}) w_{1,y}\\
\ \ \ \ \ \ \ \ \ \ +\gamma w_1(1+\varepsilon_1(t)-w_1-kw_2),\; &\ t\in((n\tau)^+,(n+1)\tau],\ y\in(0,1), \\[2mm]
w_{2,t}=\frac{1}{r^2(t)}w_{2,yy}-(\frac{\beta}{r(t)}-\frac{r'(t)y}{r(t)}) w_{2,y}\\
\ \ \ \ \ \ \ \ \ \ +w_2(1+\varepsilon_2(t)-w_2-hw_1),\; &\ t\in((n\tau)^+,(n+1)\tau],\ y\in(0,q(t)), \\[2mm]
w_1((n\tau)^+,y)=G(w_1(n\tau,y)),\; &\ y\in(0,1),\\[2mm]
w_2((n\tau)^+,y)=H(w_2(n\tau,y)),\; &\ y\in(0,q(t)),\\[2mm]
w_{1,y}(t,0)=w_1(t,1)=w_{2,y}(t,0)=w_2(t,q(t))=0,\; &\ t\in(0,\infty),\\[2mm]
q(0)=\frac{s_0}{r_0},w_1(0,y)=u_0(r_0y),w_2(0,y)=u_0(r_0y),\; &\ y\in[0,\infty).
\end{array} \right.
\]
\end{small}
We consider the following two situations.

$(i)$ If $r_\infty\leq s_\infty$, then there exists a large integer $N$ such that $r(t)\leq s(t)$ for all $t\geq N\tau$, which gives $q(t)\geq 1$ for $t\geq N\tau$. We can derive that $\|w_1\|_{C^{1,2}((n\tau,(n+1)\tau]\times[0,1])}
+\|w_2\|_{C^{1,2}((n\tau,(n+1)\tau]\times[0,1])}\leq C$ for all $n\geq N$, where $C$ is a positive constant. Let
\[
U_i(t,y)=w_1(t+t_i,y),\ V_i(t,y)=w_2(t+t_i,y)
\]
for $t\in[0,\tau]$ and $y\in[0,1]$, where $t_i=t_i^*+i \tau$ with $t_i^*\in[0,\tau)$, $i\in\mathbb{N}$ and up to a subsequence, $t_i^*\rightarrow t_0$ as $i\rightarrow\infty$. Then, $r_\infty<\infty$ implies that $\lim\limits_{i\rightarrow\infty}r'(t_i)=0$, and we have $(U_i,V_i)\rightarrow (\overline{u},\overline{v})$ for $t\in[0,\tau]$ and $y\in[0,1]$, where $\overline{u}$ satisfies
\begin{small}
\[
\left\{
\begin{array}{lll}
\overline{u}_{t}=\frac{D}{r^2_\infty}\overline{u}_{yy}-\frac{\alpha}{r_\infty} \overline{u}_{y}+\gamma \overline{u}(1+\varepsilon_1(t+t_0)-\overline{u}-k\overline{v}),\; &\ t\in(0^+,\tau),\ y\in(0,1), \\[2mm]
\overline{u}(0^+,y)=G(\overline{u}(0,y)),\; &\ y\in(0,1),\\[2mm]
\overline{u}_{y}(t,0)=\overline{u}(t,1)=0,\; &\ t\in(0,\tau),\\[2mm]
\overline{u}(0,\frac{x_0}{r_\infty})>0.
\end{array} \right.
\]
\end{small}
It follows from strong maximum principle and Hopf's boundary lemma that $\overline{u}(t,y)>0$ and $\overline{u}_y(t,1)<0$ for $t\in(0,\tau)$ and $y\in(0,1)$. Then, there exists a $\sigma>0$ such that $\overline{u}_y(t,1)\leq -\sigma$ for $t\in(\tau/6,\tau)$, whence follows that
\[
\begin{array}{llllll}
\lim\limits_{i\rightarrow\infty}r'(t_i+\tau/3)&=
\lim\limits_{i\rightarrow\infty}-\mu_1 u_x(t_i+\tau/3,r(t_i+\tau/3))\\
&=-\mu_1\lim\limits_{i\rightarrow\infty}\frac{(U_i)_y(\tau/3,1)}{r(t_i+\tau/3)}
\geq\frac{\mu_1\sigma}{r_\infty}.
\end{array}
\]
This contradicts $\lim\limits_{i\rightarrow\infty}r'(t_i)=0$.

$(ii)$ If $r_\infty>s_\infty$, then $s_\infty<\infty$ and there exists a large $n\in\mathbb{N}$ such that $s(t)<r(t)$ for all $t\geq N\tau$. Similarly as above, we can prove this situation is also impossible. Therefore, we have $\lim\limits_{t\rightarrow\infty}\|u(t,\cdot)\|_{C([0,r(t)])}=0$.
\end{pf}

Based on the above analysis, competition outcomes can be divided in Theorem 3.7.
\begin{thm}
Assume assumptions $(\mathcal{A}1)$ and $(\mathcal{A}2)$ hold. Then there are four cases of spreading and vanishing to the problem \eqref{a01}.

$(i)$ If $r_\infty\leq r_*$ and $s_\infty\leq s_*$, then species $u$ and $v$ vanish.

$(ii)$ If $r_\infty> r^*$ and $s_\infty> s^*$, then species $u$ and $v$ coexist.

$(iii)$ If (a) $r_\infty> r^*$ and $s_\infty\leq s_*$ or (b) $s_*<s_\infty\leq s^*$, then species $u$ and $v$ satisfies $\lim\limits_{m\rightarrow\infty}(u,v)(t+m\tau,x)=(U^*(t),0)$ locally uniformly in $[0,\tau]\times[0,\infty)$.

$(iv)$ If (c) $r_\infty\leq r_*$ and $s_\infty> s^*$ or (d) $r_*<r_\infty\leq r^*$, then species $u$ and $v$ satisfies $\lim\limits_{m\rightarrow\infty}(u,v)(t+m\tau,x)=(0,V^*(t))$ locally uniformly in $[0,\tau]\times[0,\infty)$.
\end{thm}
\begin{pf}
Three different intervals for $r_\infty$ will be taken into consideration to show this theorem.

$(1)$ $r_\infty\leq r_*$. It follows from Lemma 3.2  that $\lim\limits_{t\rightarrow\infty}\|u(t,\cdot)\|_{C([0,r(t)])}=0$. In this situation, $s_*<s_\infty\leq s^*$ is impossible. In fact, if $s_*<s_\infty\leq s^*$, we can prove that $r_\infty>r_*$ similarly as in Lemma 3.6, which leads to a contradiction. Then, if $s_\infty\leq s_*$,  species $v$ vanishes similarly by Lemma 3.2 and case $(i)$ holds. If $s_\infty>s^*$, then we have $s_\infty=\infty$ similarly as Lemma 3.5, which together with Theorem 3.3 parallelly deduces that species $v$ spreads eventually and case $(c)$ occurs.

$(2)$ $r_\infty>r^*$. By using Lemma 3.5 gives $r_\infty=\infty$. If $s_\infty\leq s_*$, then by virtue of Theorem 3.3, we have species $u$ spreads and $v$ vanishes, which indicates that case $(a)$ holds. If $s_\infty>s^*$, then $s_\infty=\infty$ can be proved similarly as Lemma 3.5, which implies coexistence happens by Theorem 3.4, and case $(ii)$ occurs. If $s_*<s_\infty\leq s^*$,  then species $u$ spreads and $v$ vanishes, and case $(b)$ takes place similarly as Lemma 3.6.

$(3)$ $r_*<r_\infty\leq r^*$. From Lemma 3.6, it is known that $s_\infty>s_*$ and  $\lim\limits_{m\rightarrow\infty}(u,v)(t+m\tau,x)=(0,V^*(t))$ locally uniformly in $[0,\tau]\times[0,\infty)$, which shows that case $(d)$ holds.
\end{pf}

\section{Sufficient conditions and spreading speeds}
We first propose some sufficient conditions for spreading and vanishing of species, and then give some estimates of spreading speeds.

To investigate the effects of impulsive interventions on dynamics of species, we first give some sufficient conditions about pulses and expanding capabilities for species spreading or vanishing.

Denote $g^*$, $g^{**}$, $g_*$ and $g_{**}$ satisfy
\[
\begin{array}{llllll}
&\lambda_1(D,\alpha,\gamma (1+\varepsilon_1(t)),g^*, (0,r_0))=0, \ \lambda_1(D,\alpha,\gamma (1+\varepsilon_1(t)-kV^*(t)),g^{**}, (0,r_0))=0,\\
&\lambda_1(D,\alpha,\gamma (1+\varepsilon_1(t)),g_*,(0,\infty))=0, \ \lambda_1(D,\alpha,\gamma (1+\varepsilon_1(t)-kV^*(t)),g_{**}, (0,\infty))=0.
\end{array}
\]

\begin{thm}
$(i)$ If $\alpha\leq0$ and $0<G'(0)\leq g_*$, then the species $u$ vanishes for all $\mu_1>0$. If $g_*<G'(0)<g^*$, then
there exists a $\mu_1^*>0$ such that the species $u$ vanishes for $0<\mu_1<\mu_1^*$, where $\mu_1^*$ depends on $u_0$ and $r_0$.

$(ii)$ If $G'(0)\geq g^{**}$, then the  species $u$ spreads for all $\mu_1>0$.  If $g_{**}<G'(0)<g^{**}$, then there exists a $\mu_1^{**}>0$ such that the species $u$ spreads for $\mu_1>\mu_1^{**}$, where $\mu_1^{**}$ depends on $u_0$, $v_0$ and $r_0$.
\end{thm}
\begin{pf}
$(i)$ It follows from the comparison principle that the solution $(u(t,x),r(t))$ to the problem \eqref{a01} satisfies $u(t,x)\leq \overline{u}(t,x)$ and $r(t)\leq \overline{r}(t)$ for $t\in[0,\infty)$ and $x\in[0,r(t)]$, where $(\overline{u}(t,x),\overline{r}(t))$ is the solution of the problem \eqref{c041}. It is known in Theorem 2.9 $(iii)$ that if $\alpha\leq0$ and $G'(0)<g_*$, then $\lim\limits_{t\rightarrow\infty}\|\overline{u}(t,\cdot)\|_{C([0,\overline{r}(t)])}=0$ for $t\geq0$, and the species $u$ vanishes. Meanwhile, if $g_*<G'(0)<g^*$, then Theorem 2.9 $(ii)$ gives that there exists a $\mu_1^*>0$ such that $\overline{r}_\infty<\infty$ and $\overline{u}$ vanishes for $0<\mu_1<\mu_1^*$, where $\mu_1^*$ depends on $u_0$ and $r_0$. Therefore, we have that vanishing of the species $u$ occurs for $0<\mu_1<\mu_1^*$.

$(ii)$ We aim to construct a lower solution $(\hat{u}(t,x), \hat{r}(t))$ to prove assertion $(ii)$. Similarly as in the proof of Theorem 3.3, we have $\limsup\limits_{m\rightarrow\infty}v(t+m\tau,x)\leq V^*(t)$   locally uniformly for $t\in[0,\tau]$ and $x\in[0,\infty)$. Then, for any sufficiently small $\delta>0$, there exists a large $N\in \mathbb{N}$ such that $v(t+m\tau,x)\leq V^*(t)+\delta$ for $m\geq N$, $t\in [0,\tau]$ and $x\in[0,\infty)$. Let $(w_1,w_2)$ satisfy the problem
\begin{small}
\[
\left\{
\begin{array}{lll}
w_{1,t}=Dw_{1,xx}-\alpha w_{1,x}+\gamma w_1(1+\varepsilon_1(t)-w_1-kw_2),\; &\ t\in((n\tau)^+,(n+1)\tau],\ x\in(0,r_0), \\[2mm]
w_{2,t}=w_{2,xx}-\beta w_{2,x}+w_2(1+\varepsilon_2(t)-w_2-hw_1),\; &\ t\in((n\tau)^+,(n+1)\tau],\ x\in(0,s_0), \\[2mm]
w_1((n\tau)^+,x)=G(w_1(n\tau,x)),\; &\ x\in(0,r_0),\\[2mm]
w_2((n\tau)^+,x)=H(w_2(n\tau,x)),\; &\ x\in(0,s_0),\\[2mm]
w_{1,x}(t,0)=w_{2,x}(t,0)=0,\; &\ t\in(0,\infty),\\[2mm]
w_1(t,r_0)=0,w_1(0,x)=u_0(x),\; &\ x\in[0,r_0],\\[2mm]
w_2(t,s_0)=M_2, w_2(0,x)=M_2,\; &\ x\in[0,s_0],
\end{array} \right.
\]
\end{small}
where $M_2$ is the bound of $v$ defined in Theorem 2.1. Noticing that $u(N\tau,x)\geq w_1(N\tau,x)$ and $w_1(N\tau,x)$ is independent of $\mu_1$ for $x\in [0,r_0]$, one can verify that for $n> N$, $(\hat{u}(t,x), \hat{r}(t))$ satisfying
\begin{small}
\[
\left\{
\begin{array}{lll}
\hat{u}_t=D\hat{u}_{xx}-\alpha \hat{u}_x+\gamma \hat{u}[1+\varepsilon_1(t)-kV^*(t)-k\delta-\hat{u}],\; &\ t\in((n\tau)^+,(n+1)\tau],\ x\in(0,\hat{r}(t)), \\[2mm]
\hat{u}((n\tau)^+,x)=G(\hat{u}(n\tau,x)),\; &\ x\in(0,\hat{r}(n\tau)),\\[2mm]
\hat{u}_x(t,0)=\hat{u}(t,\hat{r}(t))=0,\; &\ t\in[N\tau,\infty),\\[2mm]
\hat{r}'(t)=-\mu_1\hat{u}_x(t,\hat{r}(t)), \hat{r}(N\tau)=r_0\; &\ t\in(n\tau,(n+1)\tau],\\[2mm]
\hat{u}(N\tau,x)
=w_1(N\tau,x),\; &\ x\in[0,r_0]
\end{array} \right.
\]
\end{small}
is a lower solution. If $G'(0)\geq g^{**}$, it follows from Theorem 2.9 $(i)$ that $\hat{r}_\infty=\infty$ for all $\mu_1>0$, which together with the comparison principle yields that $r_\infty\geq \hat{r}_\infty=\infty$ and the species $u$ spreads for all $\mu_1>0$. If $g_{**}<G'(0)<g^{**}$, then  Theorem 2.9 $(ii)$ gives there exists  a $\mu_{1}^{**}>0$ such that $\hat{r}_\infty=\infty$ for $\mu_1>\mu_1^{**}$, where $\mu_1^{**}$ depends on $w_1(N\tau,x)$ and $r_0$.  By virtue of the comparison principle, we have  $r_\infty\geq \hat{r}_\infty=\infty$ and the spreading of the species $u$ occurs for $\mu_1>\mu_1^{**}$, where $\mu_1^{**}$ depends on $u_0$, $v_0$ and $r_0$.
\end{pf}

The results of $v$ is parallel to the results of $u$ shown in Theorem 4.1. And let $p^*$ and $p_*$ satisfying
\[
\begin{array}{llllll}
\lambda_1(1,\beta, 1+\varepsilon_2(t), p^*, (0,s_0))=0, \
\lambda_1(1,\beta,1+\varepsilon_2(t), p_*,(0,\infty))=0.
\end{array}
\]
\begin{thm}
If either $\beta\leq0$ and $0<H'(0)\leq p_*$, or $p_*<H'(0)<p^*$ and $0<\mu_2<\mu_2^*$ hold, then vanishing of the species $v$ occurs, where $\mu_2^*$ is defined similarly as in Theorem 4.1 $(i)$. In this situation, the following assertions hold.

$(i)$ If $G'(0)\geq g^*$, then spreading of the species $u$ occurs for any $\mu_1>0$.

$(ii)$ If $g_*<G'(0)<g^*$, then there exists a $\mu_1^\diamond>0$ such that spreading of the species $u$ occurs when $\mu_1>\mu_1^\diamond$ and vanishing of the species $u$ occurs when $0<\mu_1\leq\mu_1^\diamond$.

$(iii)$ If $\alpha\leq0$ and $0<G'(0)\leq g_*$, then vanishing of the species $u$ occurs for any $\mu_1>0$.
\end{thm}

\begin{pf}
If  $\beta\leq0$ and $0<H'(0)\leq p_*$, or $p_*<H'(0)<p^*$ and $0<\mu_2<\mu_2^*$ holds, then it follows in a similar way from Theorem 4.1 $(i)$ that the species $v$ vanishes. Then, the problem \eqref{a01} can be viewed as one single species free boundary problem for $u$. By virtue of Theorem 2.9, this theorem can be obtained.
\end{pf}

The minimal impulsive intensity for the spreading of the species $u$ can be given under a special case of impulsive function.

\begin{thm}
Assume $G(u)=gu$ with $g>0$. For any given $r_0$,
there exists
$$g_{min}:=\min\{g^\bigtriangleup>0: u \ spreads \ successfully \ for \ g=g^\bigtriangleup\}$$
such that $u$ spreads successfully for the problem \eqref{a01} regardless of $u_0$, $v_0$, $s_0$, $\mu_1$ and $\mu_2$  if and only if $g\geq g_{min}$.
\end{thm}
\begin{pf}
Define $A:=\{g^\bigtriangleup>0: u \ spreads \ succesfully \ for \ g=g^\bigtriangleup\}$. Theorem 4.1 $(ii)$ implies that $g^{**}\in A$, which gives that $A\neq\emptyset$. Then, define $g_{min}:=\inf A$. If some $g^\bigtriangleup\in A$, then for all $g>g^\bigtriangleup$, we have $u((n\tau)^+,x)=gu(n\tau,x)>g^\bigtriangleup u(n\tau,x)$, which together with  Lemma 2.2 yield that for all $g>g^\bigtriangleup$, $u$ spreads successfully, which means $g\in A$ for all $g>g^\bigtriangleup$. Therefore, if $g>g_{min}$, then $u$ spreads successfully regardless of $u_0$, $v_0$, $s_0$, $\mu_1$ and $\mu_2$.

We show $g_{min}\in A $. If $g_{min}\not\in A$, then $r_\infty<\infty$. It follows from Theorem 2.1 that we can find a  time $t^*$ such that $r(t^*)>r_0$. Taking time $t^*$ as a new initial time, for fixed initial habitat $r(t^*)$, there exists a new $g^*_{min}$ such that for any $g>g^*_{min}$, the species $u$ spreads successfully and $r_\infty=\infty$. Noticing that $g_{min}>g^*_{min}$ from Theorem 2.5, then we have $r_\infty=\infty$, which leads to a contradiction.
\end{pf}

For a special case, the effects of pulse timing on competition outcomes when pulse intensities are fixed can be obtained.
\begin{rmk}
Suppose $G(u)=gu$ and $H(v)=pv$ with $g,p>0$. Assume $\alpha\leq0$. Under the case without pulses, that is, $g=p=1$, competition exclusion occurs in sense that the species $u$ spreads and $v$ vanishes for some specified parameters in \eqref{a01}. When different pulses are introduced, we assume such parameters be fixed and choose fixed $g<1$ and $p\geq1$. Then we investigate what kind of competition outcome will be as pulse timing $\tau$ is sufficiently small or large.

$(i)$ If $\tau\rightarrow0^+$, we have $\lambda_1(D,\alpha,\gamma (1+\varepsilon_1(t)),  g,\tau,(0,\infty))=\frac{\alpha^2}{4D}-\frac{1}{\tau}\ln g-\gamma\rightarrow+\infty$ and $\lambda_1(1,\beta, 1+\varepsilon_2(t), p, \tau, (0,\infty))=\frac{\beta^2}{4}-\frac{1}{\tau}\ln p-1\rightarrow-\infty$, which together with Theorem 3.1 $(iii)$ yield the opposite competition exclusion occurs, that is,  the species $v$ spreads and $u$ vanishes.

$(ii)$ If $\tau\rightarrow\infty$, $\lim\limits_{\tau\rightarrow\infty}\lambda_1(D,\alpha,\gamma (1+\varepsilon_1(t)),  g,\tau,(0,\infty))=\frac{\alpha^2}{4D}-\gamma$ and $\lim\limits_{\tau\rightarrow\infty}\lambda_1(1,\beta, 1+\varepsilon_2(t), p, \tau, (0,\infty))=\frac{\beta^2}{4}-1$, which implies that  long time behaviors of $u$ and $v$ will be consistent with those for $g=p=1$, then competition outcomes can not be changed, and the species $u$ spreads and $v$ vanishes.

Moreover, if $\tau$ is not small or large, the competition exclusion can be altered to a coexistence state, which can be similarly shown as in Fig. $3$ $(2a-2c)$.
\end{rmk}

 The following is devoted to giving some sufficient conditions of expanding  capability according to the initial habitat for species spreading and vanishing. Under the assumptions $(\mathcal{A}1)$ and $(\mathcal{A}2)$,  the minimal initial region for the spreading of species can be  introduced in the following theorem. The definition is similar as Theorem 2 in \cite{wu2015jde} and the proof is omitted here. From now on we assume  $(\mathcal{A}1)$ and $(\mathcal{A}2)$ hold unless otherwise specified.

\begin{thm}
There exists
$$r_{min}:=\min\{r^\bigtriangleup>0: u \ spreads \ successfully \ for \ r_0=r^\bigtriangleup\}$$
such that $u$ spreads successfully for the problem \eqref{a01} regardless of $u_0$, $v_0$, $s_0$, $\mu_1$ and $\mu_2$  if and only if $r_0\geq r_{min}$. Moreover, $r_*<\hat{r}_*\leq r_{min}\leq r^*$, where $r_*$ satisfies $\lambda_1(D,\alpha,\gamma(1+\varepsilon_1(t)),G'(0), (0, r_*))=0$, $r^*$ satisfies $\lambda_1(D,\alpha,\gamma(1+\varepsilon_1(t)-kV^*(t)),G'(0),(0,r^*))=0$, and $\hat{r}_*$ satisfies $\lambda_1(D,\alpha, \gamma(1+\varepsilon_1(t)-kV_*(t)),G'(0),(0,\hat{r}_*))=0$
with $V^*(t)$ and $V_*(t)$  defined in \eqref{c01} and \eqref{c013}.
\end{thm}

\begin{lem}
$(i)$ If $0<r_0<r_*$, then there exists a $\mu_1^\circ>0$ such that $r_\infty<\infty$ for $0<\mu_1<\mu_1^\circ$, where $\mu_1^\circ$ depends on $u_0$ and $r_0$.

$(ii)$ If $r_0\geq r_{min}$, then $r_\infty=\infty$ for all $\mu_1>0$.

$(iii)$ If $0<r_0<r_{min}$, then there exists a $\mu_{1*}>0$ such that $r_\infty=\infty$ for $\mu_1>\mu_{1*}$, where $\mu_{1*}$ depends on $u_0$, $v_0$ and $r_0$.

$(iv)$ If $r_0<r_{min}$, then there exists a $\mu_1^\star>0$ such that $r_\infty=\infty$ for $\mu_1>\mu_1^\star$ and $r_\infty<\infty$ for $0<\mu_1\leq\mu_1^\star$, where $\mu_1^\star>0$ depends on $u_0$, $v_0$, $r_0$, $s_0$ and $\mu_2$.
\end{lem}
\begin{pf}
Assertions $(i)$ and $(iii)$ can be proved by the proofs of Theorem 4.1 with minor modifications, and assertion $(ii)$ can be directly obtained by the definition of $r_{min}$.
By the similar methods in Theorem 3 in \cite{wu2015jde}, the existence of sharp value $\mu_1^\star$ for the situation $r_0<r_{min}$  in $(iv)$ can be proved.
\end{pf}

\begin{thm}
$(i)$ If $r_*\leq r_0 < \hat{r}_*$ and $s_\infty=\infty$, then there exists a $\bar{\mu}_1>0$ such that $r_\infty<\infty$ for $0<\mu_1<\bar{\mu}_1$, where $\bar{\mu}_1>0$ depends on $u_0$, $v_0$, $r_0$, $s_0$ and $\mu_2$.

$(ii)$ If $r_0\geq r_*$, $0<s_0<s_*$ and $0<\mu_2<\mu_2^\circ$, then $u$ spreads successfully and $r_\infty=\infty$ for any $\mu_1>0$, where $\mu_2^\circ$ is similarly defined in Lemma 4.5 $(i)$.
\end{thm}
\begin{pf}
 $(i)$ Due to $s_\infty=\infty$, we can prove that $\liminf\limits_{m\rightarrow\infty}v(t+m\tau,x)\geq V_*(t)$  uniformly for $t\in[0,\tau]$ and locally uniformly for $x\in[0,\infty)$ similarly as in Lemma 3.5. Then, for $l>\hat{r}_*^\delta$ and sufficiently small $\delta>0$, there exists a $N\in \mathbb{N}$ such that for $m\geq N$, $v(t+m\tau,x)\geq V_*(t)-\delta$ for $t\in[0,\tau]$ and $x\in[0,l]$, where $\hat{r}_*^\delta$ satisfies $\lambda_1(D,\alpha,\gamma[1+\varepsilon_1(t)-k(V_*(t)-\delta)],G'(0),
(0,\hat{r}_*^\delta))=0$ and $\hat{r}_*^\delta>\hat{r}_*$. To stress on the dependence of $u$ and $r$ on $\mu_1$ and $\mu_2$,  we denote $r^{\mu_1,\mu_2}(t)$ and $u^{\mu_1,\mu_2}(t,x)$. Recalling $r_0<\hat{r}_*^\delta$, it follows from $r'(t)\leq C_1\mu_1$ for $t\in(n\tau,(n+1)\tau]$ and Lemma 2.3 that there exists a sufficiently small $\mu_1^\delta$ such that $r^{\mu_1,\mu_2}(N\tau)<\hat{r}_*^\delta$ for $0<\mu_1\leq \mu_1^\delta$. For $n\geq N$ and $0<\mu_1\leq \mu_1^\delta$, let $(\tilde{u},\tilde{r})$ satisfies
\begin{small}
\[
\left\{
\begin{array}{lll}
\tilde{u}_t=D\tilde{u}_{xx}-\alpha \tilde{u}_x+\gamma \tilde{u}[1+\varepsilon_1(t)-k(V_*(t)-\delta)-\tilde{u}],\; &\ t\in((n\tau)^+,(n+1)\tau],\ x\in(0,\tilde{r}(t)), \\[2mm]
\tilde{u}((n\tau)^+,x)=G(\tilde{u}(n\tau,x)),\; &\ x\in(0,\tilde{r}(n\tau)),\\[2mm]
\tilde{u}_x(t,0)=\tilde{u}(t,\tilde{r}(t))=0,\; &\ t\in[N\tau,\infty),\\[2mm]
\tilde{r}'(t)=-\mu_1\tilde{u}_x(t,\tilde{r}(t)),\; &\ t\in(n\tau,(n+1)\tau],\\[2mm]
\tilde{r}(N\tau)\in(r^{\mu_1^\delta,\mu_2}(N\tau),\hat{r}_*^\delta),\
\tilde{u}(N\tau,x)>u^{\mu_1^\delta,\mu_2}(N\tau,x),\; &\ x\in[0,\tilde{r}(N\tau)].
\end{array} \right.
\]
\end{small}
By virtue of Theorem 2.8 $(ii)$, we have there exists a $\mu_1^*>0$ such that for $0<\mu_1\leq\mu_1^*$, $\tilde{r}_\infty\leq \hat{r}_*^\delta$ and $\tilde{u}$ occurs vanishing. It follows from Lemma 2.3  that the initial values satisfies
\[u^{\mu_1,\mu_2}(N\tau,x)\leq u^{\mu_1^\delta,\mu_2}(N\tau,x) <\tilde{u}(N\tau,x), \ r^{\mu_1,\mu_2}(N\tau)\leq r^{\mu_1^\delta,\mu_2}(N\tau) <\tilde{r}(N\tau)\]
for $0<\mu_1<\bar{\mu}_1:=\min\{\mu_1^\delta,\mu_1^*\}$ and $x\in[0,r^{\mu_1,\mu_2}(N\tau)]$, which together with the comparison principle implies that for $0<\mu_1<\bar{\mu}_1$, we have $u^{\mu_1,\mu_2}(t,x)\leq \tilde{u}(t,x)$ and $r^{\mu_1,\mu_2}(t)\leq\tilde{r}(t)$ for $t\geq N\tau$ and $x\in[0,r^{\mu_1,\mu_2}(t)]$. Therefore, $r_\infty^{\mu_1,\mu_2}\leq\tilde{r}_\infty<\infty$ for $0<\mu_1<\bar{\mu}_1$.

$(ii)$ Due to $0<s_0<s_*$ and $0<\mu_2<\mu_2^\circ$, we obtain in a similar way as done in  Lemma 4.5 $(i)$ that $s_\infty<\infty$ and $\lim\limits_{t\rightarrow\infty}\|v(t,\cdot)\|_{C([0,s(t)])}=0$. We aim to prove that $r_\infty=\infty$. By contradiction, if $r_\infty<\infty$, then by $r_0\geq r_*$, there are two situations  for $r_\infty$: either $r_*<r_\infty\leq r^*$ or $r_\infty>r^*$. If $r_*<r_\infty\leq r^*$, then Lemma 3.6 implies that $s_\infty=\infty$, which contradicts $s_\infty<\infty$. The other situation is also impossible due to Lemma 3.5. Therefore, $r_\infty=\infty$ and $u$ always spreads successfully.
\end{pf}

In summary,  different conditions for the existence of the sharp value $\mu_1^\star$ under assumptions $(\mathcal{A}1)$ and $(\mathcal{A}2)$ are given as follows.
\begin{thm}
There exists a $\mu_1^\star>0$ such that $r_\infty<\infty$ for $0<\mu_1<\mu_1^\star$ and $r_\infty=\infty$ for $\mu_1\geq\mu_1^\star$ provided that either one situation holds:

$(i)$ $0<r_0<r_*$;

$(ii)$ $r_*\leq r_0<\hat{r}_*$ and $s_0\geq s_{min}$;

$(iii)$ $r_*\leq r_0< \hat{r}_*$, $0<s_0<s_{min}$ and $\mu_2>\mu_{2*}$.
\end{thm}
\begin{pf}
Recalling Lemma 4.5 $(i)$ and $(iv)$, if $(i)$ holds, then $\mu_1^\star\geq \mu_1^\circ>0$. If either $s_0\geq s_{min}$ or $0<s_0<s_{min}$ and $\mu_2>\mu_{2*}$ holds, we have $s_\infty=\infty$ similarly by Lemma 4.5 $(ii)$ and $(iii)$ respectively, which together with $r_*\leq r_0<\hat{r}_*$ and Theorem 4.6 $(i)$ give the existence of $\mu_1^\star$.
\end{pf}

The effects of initial values on dynamics of species under assumptions $(\mathcal{A}1)$ and $(\mathcal{A}2)$ are shown in Theorem 4.8.
\begin{thm}
If $r_0<r_*$, then $u$ vanishes eventually for sufficiently small $\|u_0(x)\|_{C([0,r_0])}$.
 Moreover, the following assertions hold.

(i) If $s_0<s_*$, then $v$ vanishes eventually for sufficiently small $\|v_0(x)\|_{C([0,s_0])}$.

(ii) If $s_0\geq s_*$, then $v$ spreads successfully for any $\|v_0(x)\|_{C([0,s_0])}$.
\end{thm}
\begin{pf}
 If $r_0<r_*$ and $\|u_0(x)\|_{C([0,r_0])}$ is sufficiently small, it can be similarly proved that vanishing occurs for $u$ by Theorem 2.8 $(ii)$ with some minor modifications. Then, the problem \eqref{a01} can be regarded as one single species free boundary problem for $v$, which together with Theorem 2.8 yields  $(i)$ and $(ii)$ hold, respectively.
\end{pf}

The rest of section aims to investigate spreading speeds of species when both two competing species spread eventually for the case that $\alpha,\beta\geq0$.

Consider the following semi wave problem
\begin{eqnarray}
\left\{
\begin{array}{lll}
U_t=DU_{xx}-K(t)U_x+\gamma U(a(t)-U),\; &\ t\in(0^+,\tau],\ x\in(0,\infty), \\[2mm]
U(0^+,x)=G(U(0,x)),\; &\ x\in(0,\infty),\\[2mm]
U(t,0)=0,\; &\ t\in(0,\tau),\\[2mm]
U(0,x)=U(\tau,x),\; &\ x\in(0,\infty),
\end{array} \right.
\label{d01}
\end{eqnarray}
where $K(t)\in E:=\{ K(t): K(t)\ is \ \tau-peridoic \ in \ time \ t \ and \ is \ a \ H\ddot{o}lder\  continuous\\ \ function  \}$.

It is known in [Lemma 5.2, 24] that the problem \eqref{d01} admits a unique positive solution $U^{K}(t,x)$ if and only if
\[
\frac{1}{4D}[\frac{1}{\tau}\int_0^\tau K(t)dt]^2-\gamma \bar{a}-\frac{1}{\tau}\ln G'(0)<0.
\]
Under the case that $\alpha\geq0$, the estimates of spreading speed of  free boundary problem \eqref{b01} for one single species is first given.
\begin{lem}
Assume $\alpha\geq0$. If $\alpha^2/(4D)-\ln G'(0)/\tau-\gamma \bar{a}<0$, then there exists a unique $K_0^a(t)\in E$ such that $K_0^a(t)=\mu_1 U_x^{K_0^a}(t,0)$ for any $\mu_1>0$. And, \begin{equation}
\frac{1}{\tau}\int_0^\tau K_0^a(t)dt<2\sqrt{D(\gamma \bar{a}+\frac{1}{\tau}\ln G'(0))}.
\label{d02}
\end{equation}
Moreover, if $r_\infty=\infty$, then
\begin{equation}
\lim\limits_{t\rightarrow\infty}\frac{r(t)}{t}=\frac{1}{\tau}\int_0^\tau K_0^a(t)dt+\alpha.
\label{d03}
\end{equation}
\end{lem}

\begin{pf}
Since the existence and uniqueness of $K_0^a(t)$, and \eqref{d02} can be similarly proved as in Lemma 5.3 in \cite{menglin22}, respectively, we only prove \eqref{d03} here to show the effects of advection and pulses.

Due to $r_\infty=\infty$, it follows from Theorem 2.8 $(i)$ that $\lim\limits_{m\rightarrow\infty}u(t+m\tau,x)=U^*(t)$ locally uniformly for $t\in [0,\tau]$ and $x\in[0,\infty)$. Then, for any small $\delta>0$, there  exists a large $N\in\mathbb{N}$ such that $u(t,x)\leq (1-\delta)^{-1}U^*(t)$ for $t\geq N\tau$ and $x\in[0,\infty)$. Similar as Lemma 5.1 in \cite{menglin22}, one can prove that the unique positive solution $U^{K_0^a}(t,x)$ to the problem \eqref{d01} satisfies $\lim\limits_{x\rightarrow\infty}U^{K_0^a}(t,x)=U^*(t)$ uniformly for $t\in[0,\tau]$, and $U_x^{K_0^a}(t,x)>0$ for $t\in[0,\tau]$ and $x\in[0,\infty)$. Then, for any small $\delta>0$, there exists a large $l_0>0$ such that $U^{K_0^a}(t,x)>(1-\delta)U^*(t)$ for $t\in[0,\tau]$ and $x\in[l_0,\infty)$.

Let
\[
\tilde{u}(t,x)=(1-\delta)^{-2}U^{K_0^a}(t,\zeta(t)-x),\ t\geq 0,\ 0\leq x\leq \zeta(t),
\]
\[
\zeta(t)=(1-\delta)^{-2}\int_0^t(K_0^a(s)+\alpha)ds+l_0+r(N\tau),\ t\geq 0.
\]
Direct calculations give that
\[
\zeta'(t)=(1-\delta)^{-2}(K_0^a(t)+\alpha),\ -\mu_1\tilde{u}_x(t,\zeta(t))=\mu_1(1-\delta)^{-2}U_x^{K_0^a}(t,0)
=(1-\delta)^{-2}K_0^a(t),
\]
which implies that $\zeta'(t)\geq \mu_1\tilde{u}_x(t,\zeta(t))$ for $t\geq0$.
Recalling $U_x^{K_0^a}(t,x)>0$ for $t\in[0,\tau]$ and $\alpha\geq0$, then
\[
\begin{array}{llllll}
&\tilde{u}_t-D\tilde{u}_{xx}+\alpha\tilde{u}_x
=(1-\delta)^{-2}(U_t^{K_0^a}+U_x^{K_0^a}\zeta'(t)-DU_{xx}^{K_0^a}-\alpha U_x^{K_0^a})\\ [2mm]
 &=(1-\delta)^{-2}[U_t^{K_0^a}+U_x^{K_0^a}(1-\delta)^{-2}(K_0^a(t)+\alpha)
 -DU_{xx}^{K_0^a}-\alpha U_x^{K_0^a}] \\ [2mm]
 &\geq (1-\delta)^{-2}(U_t^{K_0^a}+U_x^{K_0^a}K_0^a(t)-DU_{xx}^{K_0^a})
 =(1-\delta)^{-2}\gamma U^{K_0^a}(a(t)-U^{K_0^a}) \\ [2mm]
 &\geq \gamma \tilde{u}(a(t)-\tilde{u}), \ t\in((n\tau)^+,(n+1)\tau], \ x\in(0,\zeta(t)).
\end{array}
\]
Since $G(u)/u$ is nonincreasing with respect to $u$, one has
\[
\begin{array}{llllll}
\tilde{u}((n\tau)^+,x)&=(1-\delta)^{-2}U^{K_0^a}((n\tau)^+,\zeta(t)-x)
=(1-\delta)^{-2}G(U^{K_0^a}(n\tau,\zeta(t)-x))\\ [2mm]
& \geq G((1-\delta)^{-2}U^{K_0^a}(n\tau,\zeta(t)-x))=G(\tilde{u}(n\tau,x)),\
x\in(0,\zeta(t)).
\end{array}
\]
It is clear that $\tilde{u}_x(t,0)=-(1-\delta)^{-2}U^{K_0^a}_x(t,\zeta(t))\leq0$ and $\tilde{u}(t,\zeta(t))=0$.
And initial value $(\tilde{u}(0,x),\zeta(0))$ satisfies $\zeta(0)=l_0+r(N\tau)>r(N\tau)$ and
\[
\begin{array}{llllll}
\tilde{u}(0,x)&=(1-\delta)^{-2}U^{K_0^a}(0,\zeta(0)-x)
\geq(1-\delta)^{-2}U^{K_0^a}(0,l_0)\\ [2mm]
& \geq (1-\delta)^{-2}U^*(0)\geq u(N\tau,x),\
x\in[0,r(N\tau)].
\end{array}
\]
It follows from the comparison principle that $u(t+N\tau,x)\leq \tilde{u}(t,x)$ and $r(t+N\tau)\leq \zeta(t)$ for $t\geq0$ and $0\leq x\leq r(t+N\tau)$. Then, we have
\[
\begin{array}{llllll}
\limsup\limits_{t\rightarrow\infty}\frac{r(t)}{t}&\leq
\limsup\limits_{t\rightarrow\infty}\frac{\zeta(t-N\tau)}{t} \\ [2mm]
& =(1-\delta)^{-2}\lim\limits_{t\rightarrow\infty}
\frac{\int_0^t(K_0^a(s)+\alpha)ds+l_0+r(N\tau)}{t}\\ [2mm]
&=(1-\delta)^{-2}(\frac{1}{\tau}\int_0^\tau K_0^a(t)dt+\alpha),
\end{array}
\]
which yields $\limsup\limits_{t\rightarrow\infty}\frac{r(t)}{t}\leq\frac{1}{\tau}\int_0^\tau K_0^a(t)dt+\alpha$ owing to the arbitrariness of  small $\delta$.

Next, we can construct a lower solution $(\hat{u}(t,x),\eta(t))$ satisfying
\[
\hat{u}(t,x)=\sqrt{1-\delta}U^{K_0^a}(t,\eta(t)-x),\ t\geq 0,\ 0\leq x\leq \eta(t),
\]
\[
\eta(t)=(1-\delta)\int_0^t(K_0^a(s)+\alpha)ds+r_0,\ t\geq 0
\]
to prove $\liminf\limits_{t\rightarrow\infty}\frac{r(t)}{t}\geq\frac{1}{\tau}\int_0^\tau K_0^a(t)dt+\alpha$, which can be proved by  the proofs of Theorem 5.4 in \cite{menglin22} with some minor modifications and is omitted.
\end{pf}

Based on the above results about spreading speeds for single species, the estimates of spreading speeds of species in competition free boundary problem \eqref{a01} for the case $\alpha,\beta\geq0$ can be obtained.
\begin{thm}
Assume $\alpha,\beta\geq0$. If $r_\infty=\infty$, then
\[
\frac{1}{\tau}\int_0^\tau K_0^{a_2}(t)dt+\alpha\leq\liminf\limits_{t\rightarrow\infty}\frac{r(t)}{t}\leq \limsup\limits_{t\rightarrow\infty}\frac{r(t)}{t}\leq\frac{1}{\tau}\int_0^\tau K_0^{a_1}(t)dt+\alpha,
\]
where $a_1(t)=1+\varepsilon_1(t)$ and $a_2(t)=1+\varepsilon_1(t)-kV^*$.
\end{thm}

\begin{pf}
Let $(\tilde{u},\tilde{r})$ be the solution of the problem \eqref{c041}. It follows from comparison principle that $r(t)\leq \tilde{r}(t)$ for $t\geq0$,  which with Lemma 4.9 implies that
\[\limsup\limits_{t\rightarrow\infty}\frac{r(t)}{t}\leq
\limsup\limits_{t\rightarrow\infty}\frac{\tilde{r}(t)}{t}
=\frac{1}{\tau}\int_0^\tau K_0^{a_1}(t)dt+\alpha,\]
where $a_1(t)=1+\varepsilon_1(t)$.

Similarly as in the proof of Theorem 3.3, we have $\limsup\limits_{m\rightarrow\infty}v(t+m\tau,x)\leq V^*(t)$  locally uniformly in $[0,\tau]\times[0,\infty)$. Then, for any sufficiently small $\delta>0$, there exists a large $N\in \mathbb{N}$ such that $v(t+m\tau,x)\leq V^*(t)+\delta$ for $m\geq N$, $t\in [0,\tau]$ and $x\in[0,\infty)$. It is easy to verify that $u(t,x)\geq \hat{u}(t,x)$ and $r(t)\geq \hat{r}(t)$ for $t\geq N\tau$ and $x\in[0,\hat{r}(t)]$, where $(\hat{u}(t,x),\hat{r}(t))$ satisfies the problem \eqref{b01} with $a(t)$ replaced by $a_2^\delta(t)=1+\varepsilon_1(t)-k(V^*+\delta)$. One can deduce that
\[
\liminf\limits_{t\rightarrow\infty}\frac{r(t)}{t}\geq
\liminf\limits_{t\rightarrow\infty}\frac{\hat{r}(t)}{t}
=\frac{1}{\tau}\int_0^\tau K_0^{a_2^\delta}(t)dt+\alpha,
\]
which together with the arbitrariness of $\delta$ yields that $\liminf\limits_{t\rightarrow\infty}\frac{r(t)}{t}\geq \frac{1}{\tau}\int_0^\tau K_0^{a_2}(t)dt+\alpha$,
where $a_2(t)=1+\varepsilon_1(t)-kV^*$.
\end{pf}

\section{Numerical simulations}
In this section, we aim to show the effects of periodic pulses and environmental perturbations on the spreading of two competitors in advective environments numerically. Noting that environmental perturbation functions $\varepsilon_1(t)$ and $\varepsilon_2(t)$ satisfy $(\mathcal{H}1)$, we take
\[
\varepsilon_i(t)=\sigma_i\sin(\frac{2\pi}{\tau}t), \  i=1,2
\]
similarly as in \cite{khan}, where $\sigma_1$ and $\sigma_2$ represent the magnitudes of environmental perturbations on competitors $u$ and $v$, respectively, and $\sigma_i<1$. In order to investigate effects of pulses intuitively and simply, we only choose linear impulsive functions $G(u)=c_1u$ and $H(u)=c_2u$ with $c_1,c_2>0$, where impulsive effects $c_1$ and $c_2$ represent the intensities of impulsive interventions. In all simulations, some parameters
\[
D=1.2,\ \gamma=1.5,\ \alpha=0.8,\ \beta=0.6,\ \mu_1=0.6,\ \mu_2=0.1,\ r_0=s_0=1
\]
are fixed, and all parameters in \eqref{a01} will be chosen to satisfy assumptions $(\mathcal{A}1)$ and $(\mathcal{A}2)$.

\begin{figure}[ht]
\quad
\begin{minipage}{0.2\linewidth}
\centerline{\includegraphics[width=4cm]{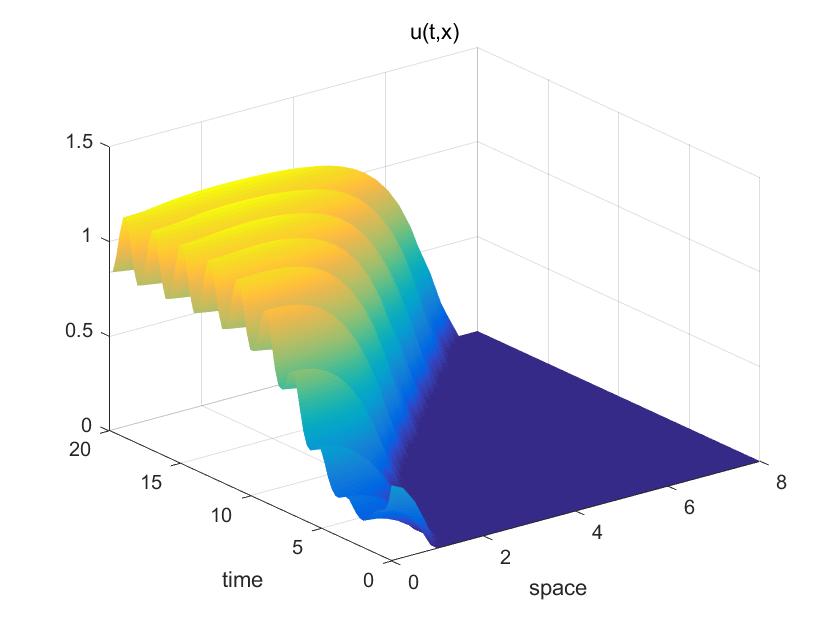}}
\centerline{\small{(a)}}
\end{minipage}
\quad
\begin{minipage}{0.2\linewidth}
\centerline{\includegraphics[width=4cm]{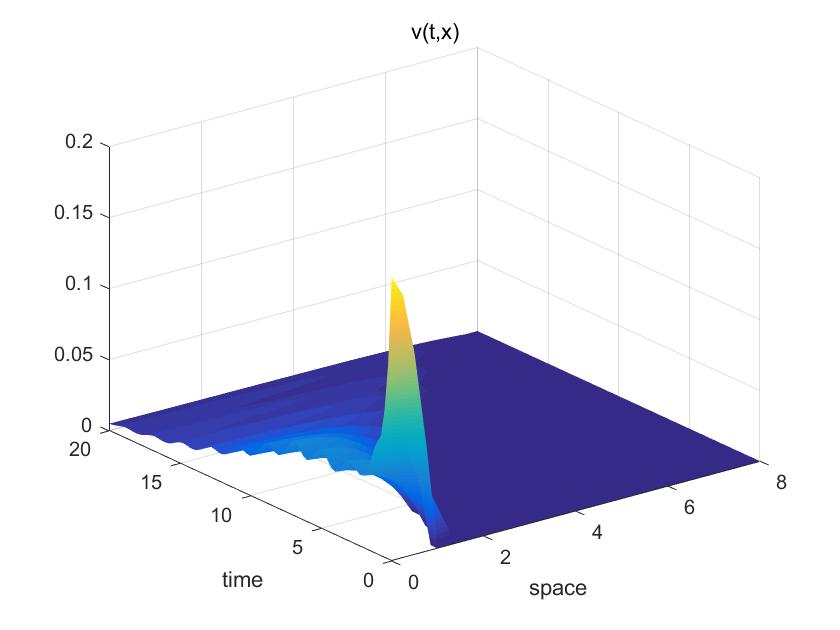}}
\centerline{\small{(b)}}
\end{minipage}
\quad
\begin{minipage}{0.2\linewidth}
\centerline{\includegraphics[width=4cm]{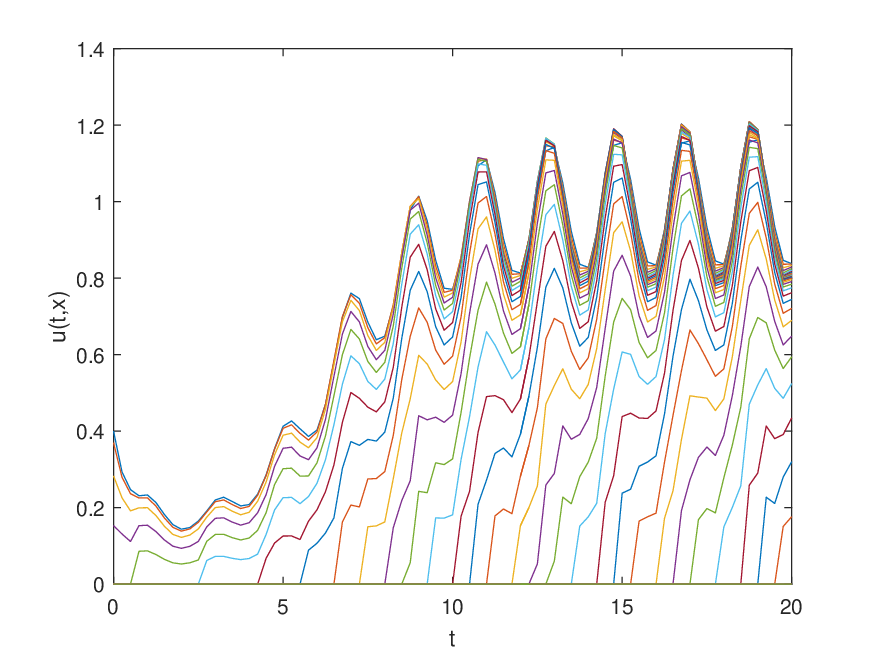}}
\centerline{\small{(c)}}
\end{minipage}
\quad
\begin{minipage}{0.2\linewidth}
\centerline{\includegraphics[width=4cm]{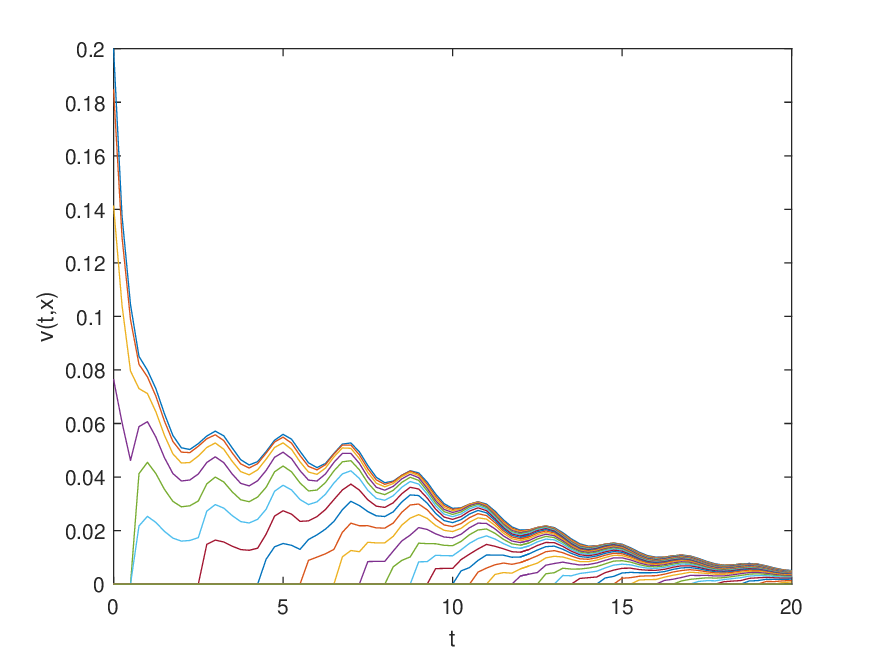}}
\centerline{\small{(d)}}
\end{minipage}
\caption{\scriptsize The dynamics of the competitors $u$ and $v$ without pulses. Graphs $(c)$ and $(d)$ are the right view of graphs $(a)$ and $(b)$, respectively, which indicate that competitor $u$ spreads and $v$ vanishes eventually.   }
\label{tu1}
\end{figure}

\begin{figure}[ht]
\begin{minipage}{0.3\linewidth}
\centerline{\includegraphics[width=5cm]{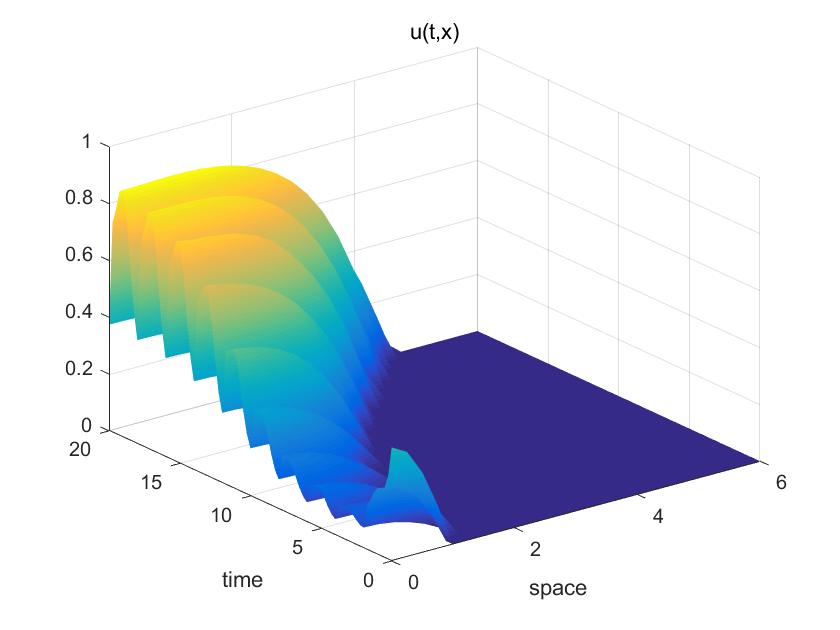}}
\centerline{\small{(1a)}}
\end{minipage}
\quad
\begin{minipage}{0.3\linewidth}
\centerline{\includegraphics[width=5cm]{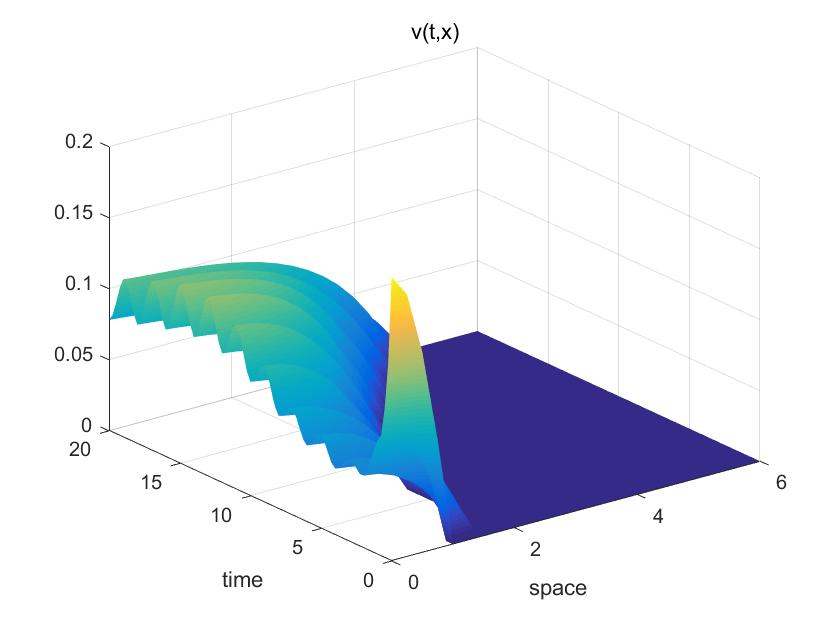}}
\centerline{\small{(1b)}}
\end{minipage}
\quad
\begin{minipage}{0.3\linewidth}
\centerline{\includegraphics[width=5cm]{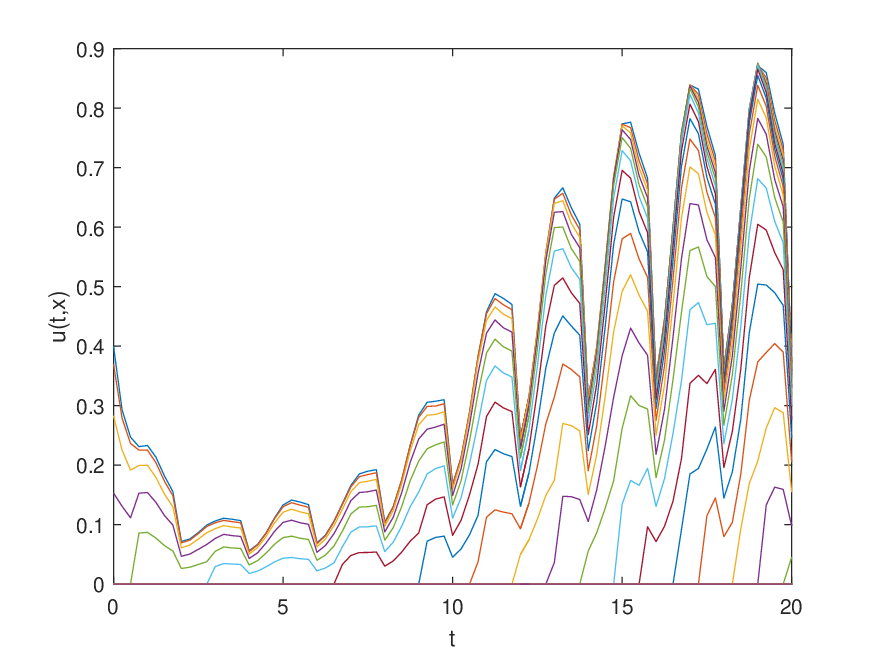}}
\centerline{\small{(1c)}}
\end{minipage}
\begin{minipage}{0.3\linewidth}
\centerline{\includegraphics[width=5cm]{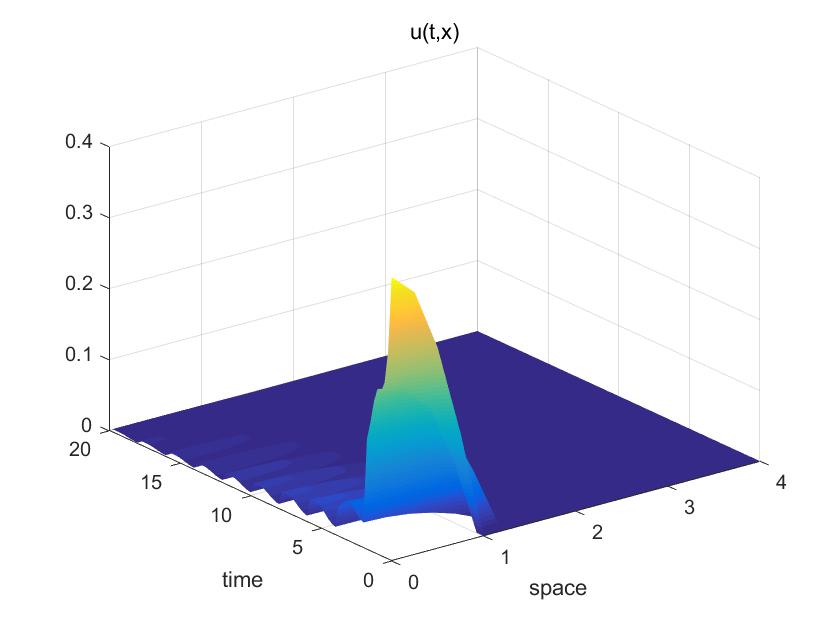}}
\centerline{\small{(2a)}}
\end{minipage}
\quad
\begin{minipage}{0.3\linewidth}
\centerline{\includegraphics[width=5cm]{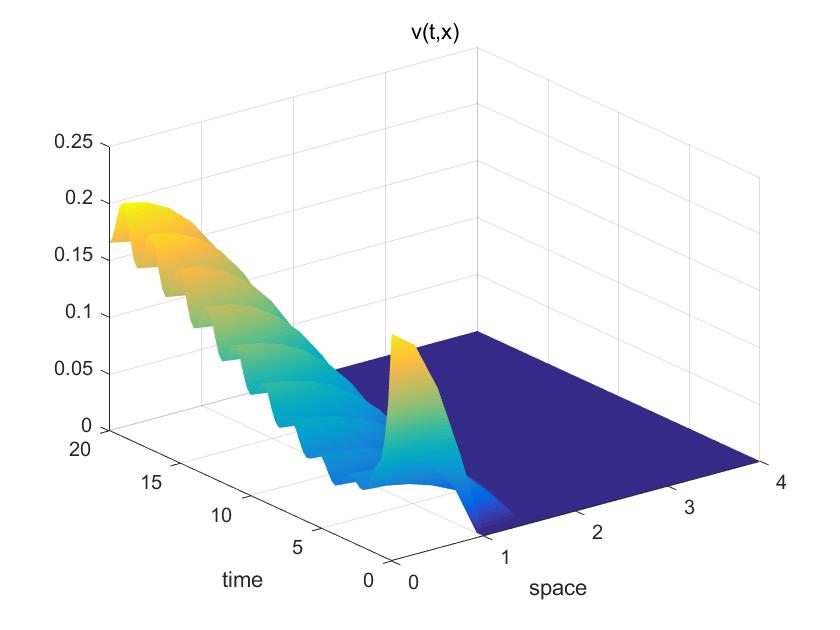}}
\centerline{\small{(2b)}}
\end{minipage}
\quad
\quad
\begin{minipage}{0.3\linewidth}
\centerline{\includegraphics[width=5cm]{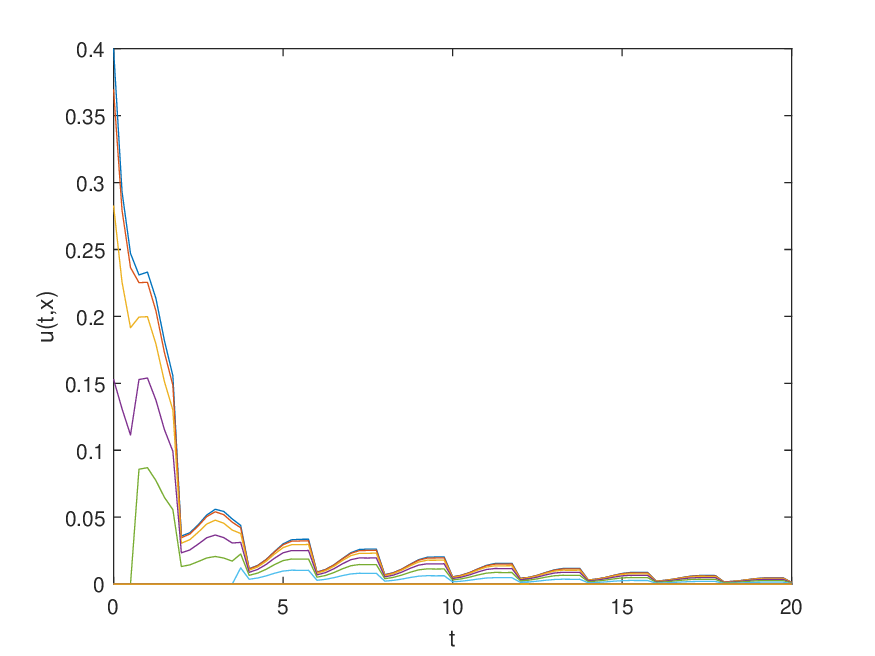}}
\centerline{\small{(2c)}}
\end{minipage}
\caption{\scriptsize The dynamics of competitors $u$ and $v$ with impulsive harvesting applied to superior $u$ at every $\tau=2$.  Graphs $(1c)$ and $(2c)$ are the right view of the spatiotemporal distribution of u in graphs $(1a)$ and $(2a)$, respectively. The different impulsive harvesting carried out at $c_1=0.5$ and $c_1=0.25$ are shown in graphs $1(a)-1(c)$ and $2(a)-2(c)$, respectively. Graphs $1(a)- 1(c)$ imply that two competitors $u$ and $v$ coexist eventually while graphs $2(a)-2(c)$ indicate that the inferior $v$ will out-compete.}
\label{tu2}
\end{figure}

Recalling Theorem 3.7, the dynamics of two competitors with periodic pulses in \eqref{a01} can be classified into four competition outcomes. Since $\int_0^\tau\varepsilon_i(t)dt=0(i=1,2)$,  periodic environmental perturbation terms $\varepsilon_1(t)$ and $\varepsilon_2(t)$ have no effects on threshold values $\lambda_1$ defined in Theorem 2.5, which  determines competition outcomes. Next, the focus of attention is on whether and how periodic impulsive interventions alter the competition outcomes.
\begin{exm}(Effects of pulses on competition outcomes)
Fix $\sigma_1=\sigma_2=0.5$ and $\tau=2$.
To focus on the effects of periodic impulsive interventions, we first fix $c_2=1$ and then choose different impulsive harvesting $(c_1\leq1)$ only on the superior $u$  for simplicity.

We first choose $c_1=1$, which means that no impulsive interventions occur. One can see from Fig. \ref{tu1} that the competitor $u$ stabilizes to a positive periodic steady state while the competitor $v$ decays to zero eventually without impulsive interventions.

 Choose $c_1=0.5$. It follows from Fig. \ref{tu2}(1a-1c) that the competitors $u$ and $v$ reach  a coexistence state eventually.  It is shown in the Fig. \ref{tu2}(1c) that impulsive harvesting is implemented in the superior $u$ at a frequency of $\tau=2$, which slows down the spread of superior $u$ and maintains the coexistence of two competitors from the comparison of Figs.  \ref{tu1} and \ref{tu2}(1a-1c) .

Choose $c_1=0.25$. The impulsive harvesting carried out on  the superior $u$ is strong enough such that the inferior $v$ turns into victory in the competition, which are exhibited in Fig. \ref{tu2}(2a-2c).
\end{exm}

Example 5.1 shows that negative impulsive effect $(c_1<1)$ only on the superior $u$ can suppress the spread of $u$ and alter the competition outcomes. We can also obtain that positive  impulsive effect $(c_2>1)$ only on the inferior $v$ have the similar effects on the competition outcomes.

\begin{figure}[ht]
\centering
\subfigure[$\tau=1$]{ {
\includegraphics[width=0.35\textwidth]{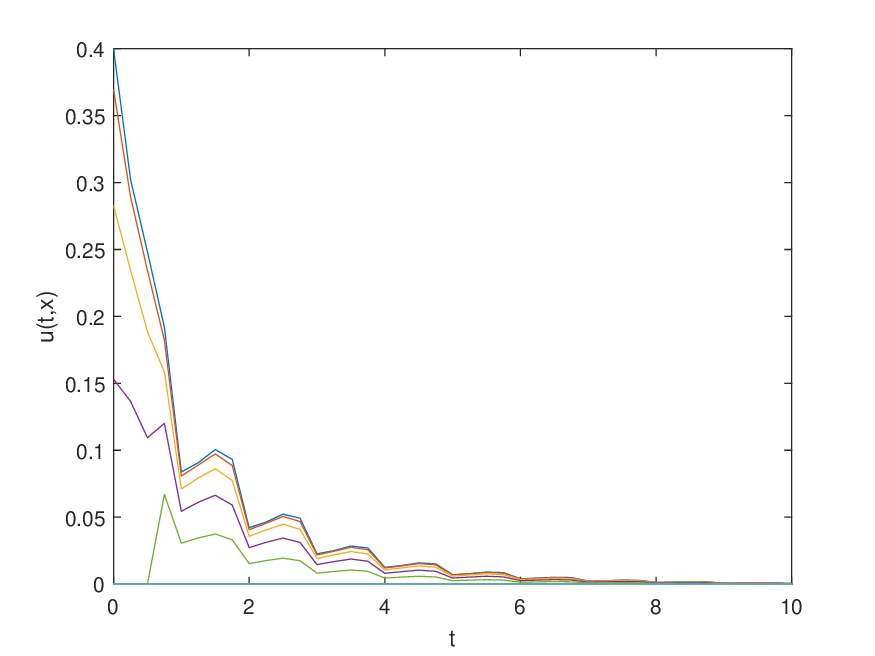}
} }
\subfigure[$\tau=1$]{ {
\includegraphics[width=0.35\textwidth]{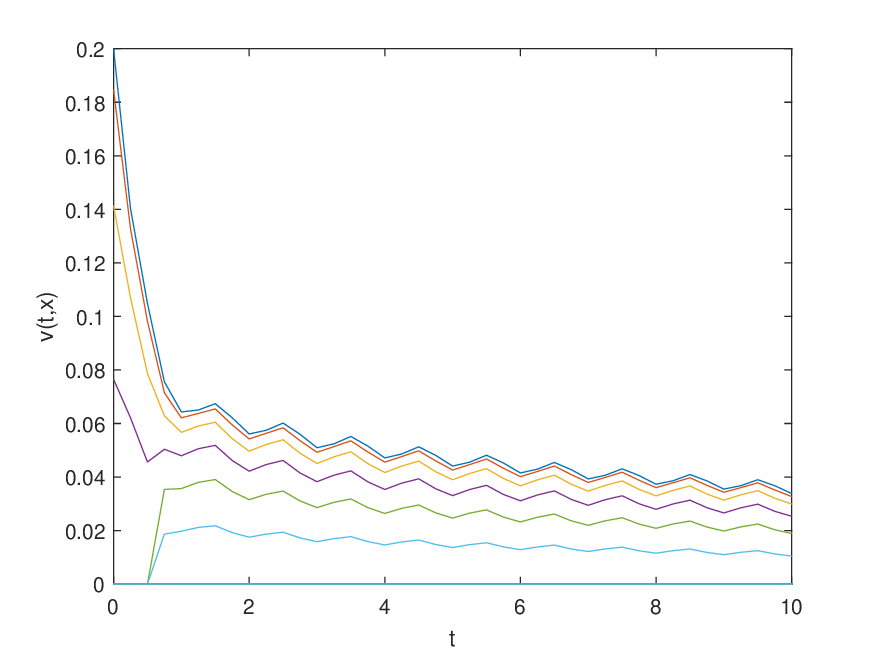}
} }
\subfigure[$\tau=4$]{ {
\includegraphics[width=0.35\textwidth]{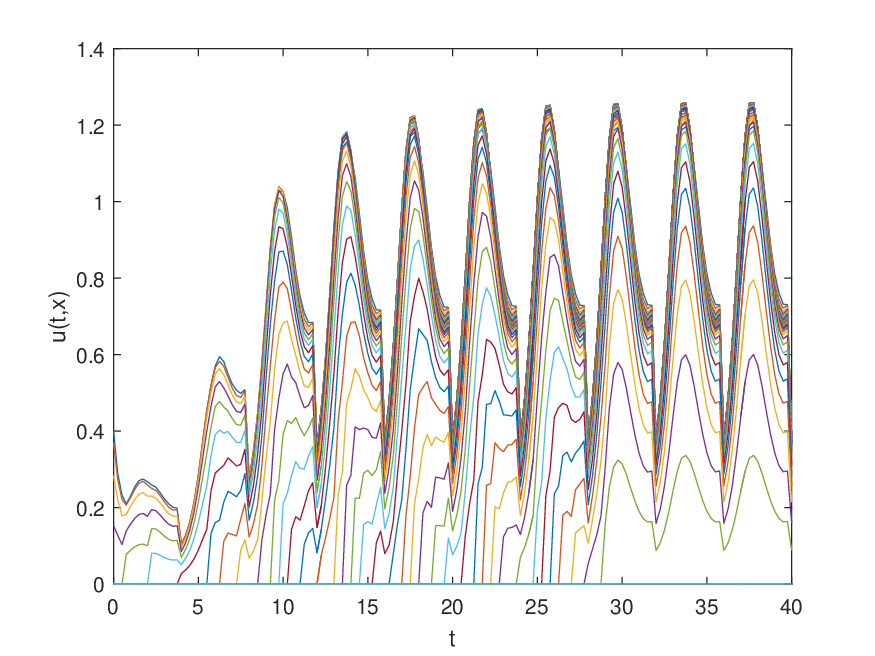}
} }
\subfigure[$\tau=4$]{ {
\includegraphics[width=0.35\textwidth]{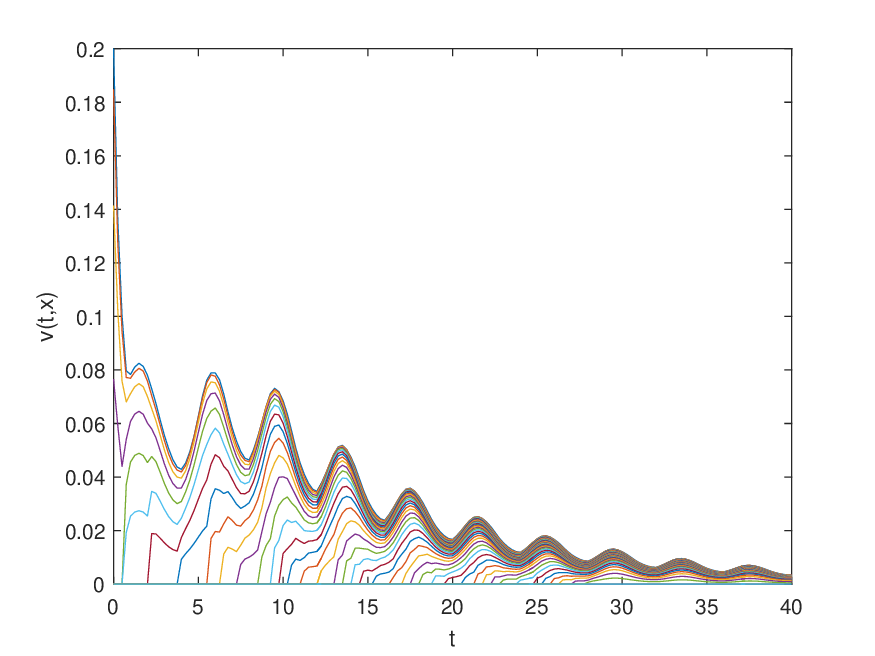}
} }
\caption{\scriptsize The dynamics of competitors $u$ and $v$ with impulsive effect $c_1=0.5$ takes place at every time $\tau=1$ in graphs $(a)$ and $(b)$ and $\tau=4$ in graphs $(c)$ and $(d)$, respectively. All graphs are the right view of the spatiotemporal distribution of competitors $u$ and $v$, which have the same interpretations as those in Figs. \ref{tu1} and \ref{tu2}.}
\label{tut1}
\end{figure}

 Competition outcomes have been shown for a sufficiently large or small $\tau$ in Remark 4.1 under a special case.  In order to understand how impulse timing $\tau$ affects the spreading of two competitors numerically, we consider the case that impulsive harvesting is implemented only in the superior $u$ for convenience,  choose a large or small time period $\tau$ and fix other parameters.
\begin{exm}(Effects of time period $\tau$ on competition outcomes)
Fix $\sigma_1=\sigma_2=0.5$, $c_1=0.5$ and $c_2=1$.

Choose $\tau=1$. Fig. \ref{tut1}(a) shows that superior $u$ is harvested at every time $\tau=1$ and will vanish eventually, while the inferior $v$ will out-compete  from Fig. \ref{tut1}(b).

Choose $\tau=4$. The superior $u$ will spread and the inferior $v$ will go to extinction in  Fig. \ref{tut1}(c-d). This competition outcome is similar to the case without pulses.
\end{exm}

Comparisons of  the results showed in Fig. \ref{tu1} and Example 5.2 reveal that  impulse timing $\tau$ is also a significant factor  that affects the long time behaviors of competitors, that is, if impulsive harvesting is applied to the superior $u$ more frequently, the competition exclusion in Fig. \ref{tu1} will be transformed to  the opposite competition exclusion (Fig. \ref{tut1}(a-b)), while competition outcomes in Fig. \ref{tu1} will not be changed for a larger $\tau$ (Fig. \ref{tut1}(c-d)), which is consistent with Remark 4.1. Also, one can see from Figs. \ref{tu1} and \ref{tu2} (1a-1c) that if the time period $\tau$ of impulsive harvesting is moderate, this kind of impulsive intervention will maintain the coexistence of two competitors.

Similar as Example 5.2, when positive impulsive intervention is carried out only on the inferior $v$, we can obtain  that the more frequently  positive intervention is applied, the more beneficial it is for the inferior species.

Examples 5.1 and 5.2 indicate that suitable impulsive strategy including positive or negative impulsive effect, pulse intensity and timing can be carried out to achieve the desired competition outcomes and help maintain ecological balance and gain sustainable development.

\begin{figure}[ht]
\centering
\subfigure[$c_1=0.5$]{ {
\includegraphics[width=0.35\textwidth]{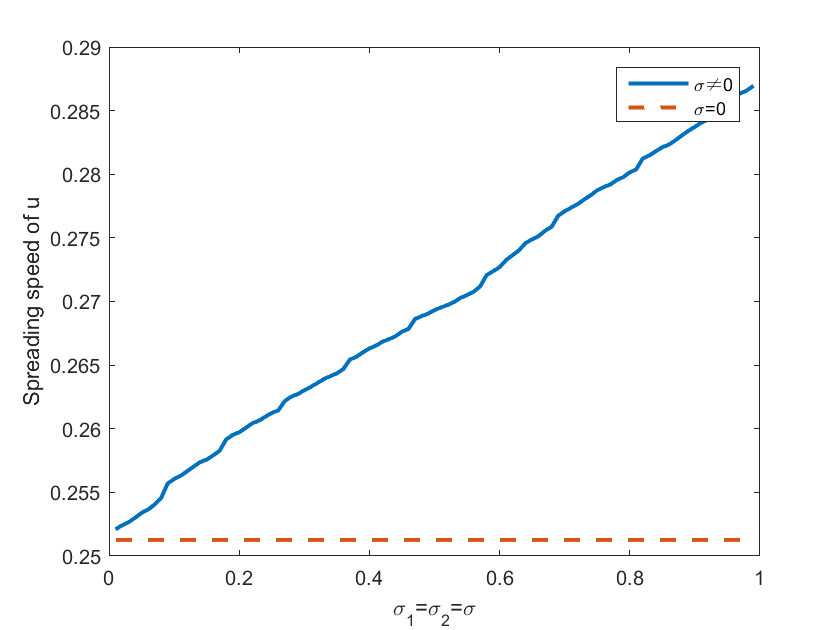}
} }
\subfigure[$c_1=0.5$]{ {
\includegraphics[width=0.35\textwidth]{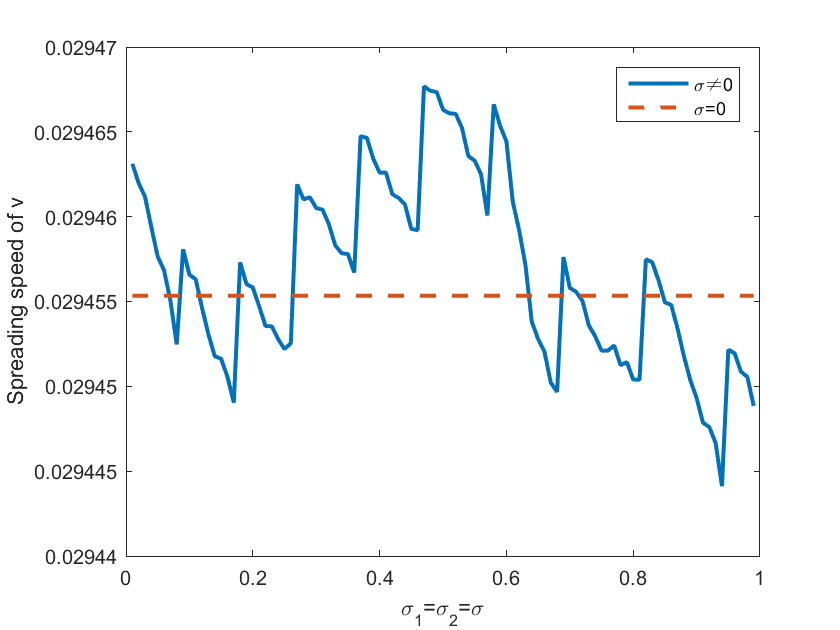}
} }
\subfigure[$c_1=0.4$]{ {
\includegraphics[width=0.35\textwidth]{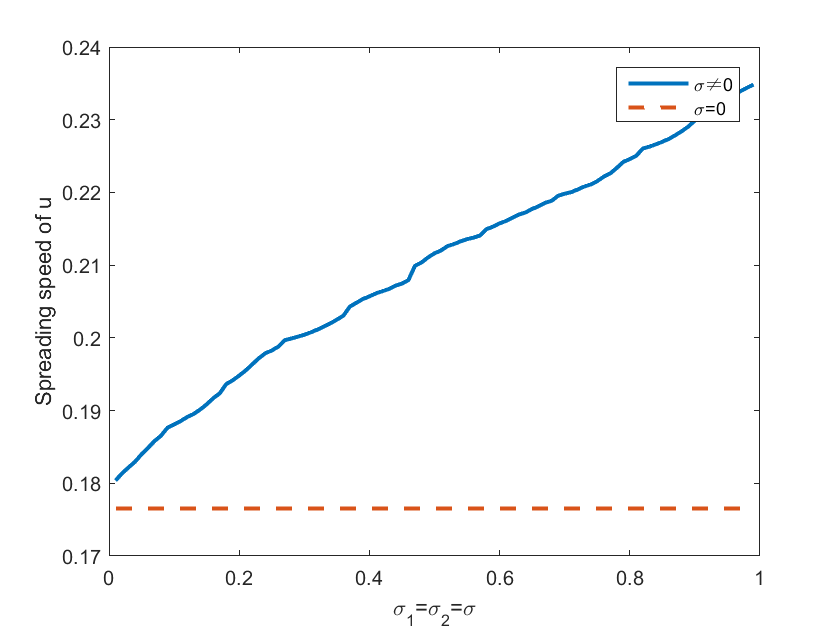}
} }
\subfigure[$c_1=0.4$]{ {
\includegraphics[width=0.35\textwidth]{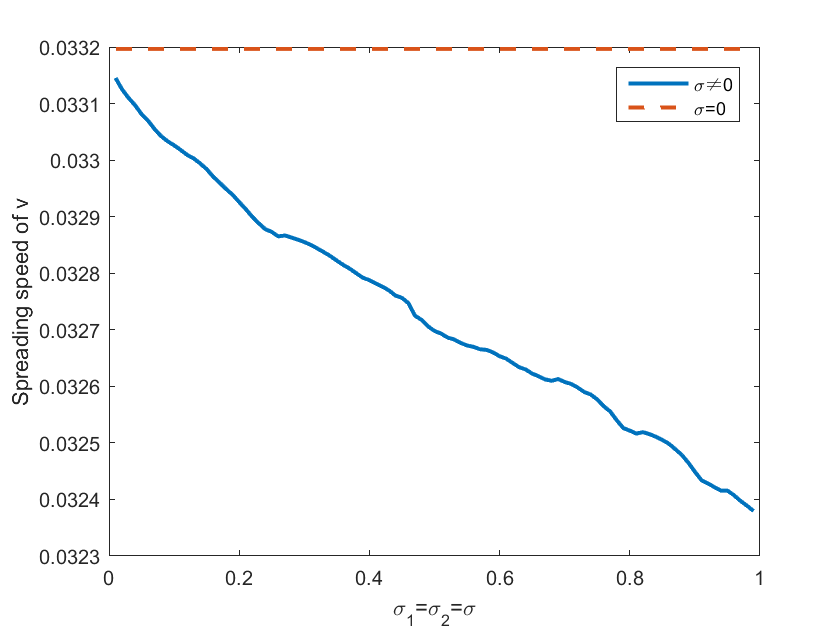}
} }
\caption{\scriptsize Spreading speeds of $u$ and $v$ with $\sigma_1=\sigma_2=\sigma\in[0.01,0.99]$ for impulsive effects $c_1=0.5$ in graphs $(a-b)$ and $c_1=0.4$ in graphs $(c-d)$. The red dashed lines stand for the spreading speeds of $u$  and $v$ in homogeneous environments. Graphs $(a)$ and $(c)$ show that the spreading speeds of $u$ with $\sigma\in[0.01,0.99]$ is bigger than those in homogeneous environments. }
\label{tus1}
\end{figure}

\begin{figure}[ht]
\centering
\subfigure[$c_1=0.5$]{ {
\includegraphics[width=0.35\textwidth]{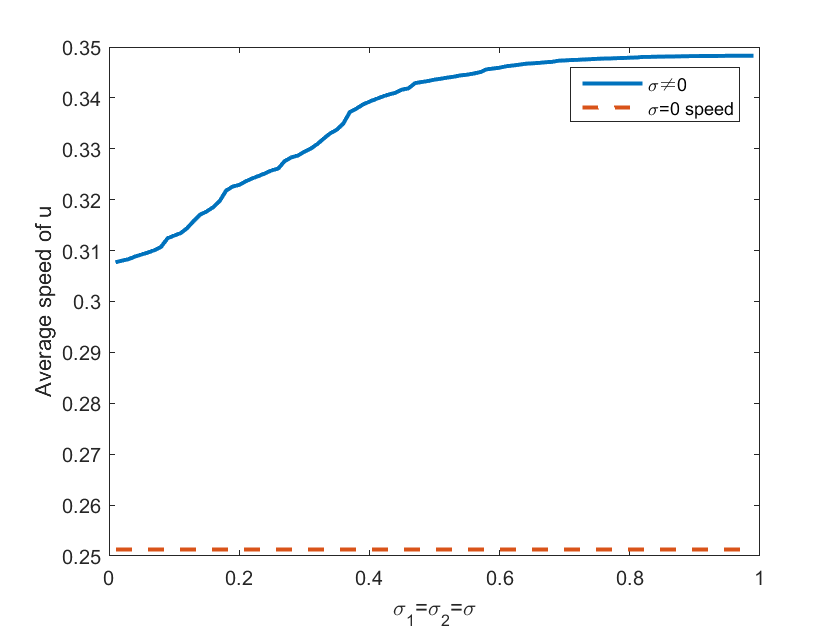}
} }
\subfigure[$c_1=0.5$]{ {
\includegraphics[width=0.35\textwidth]{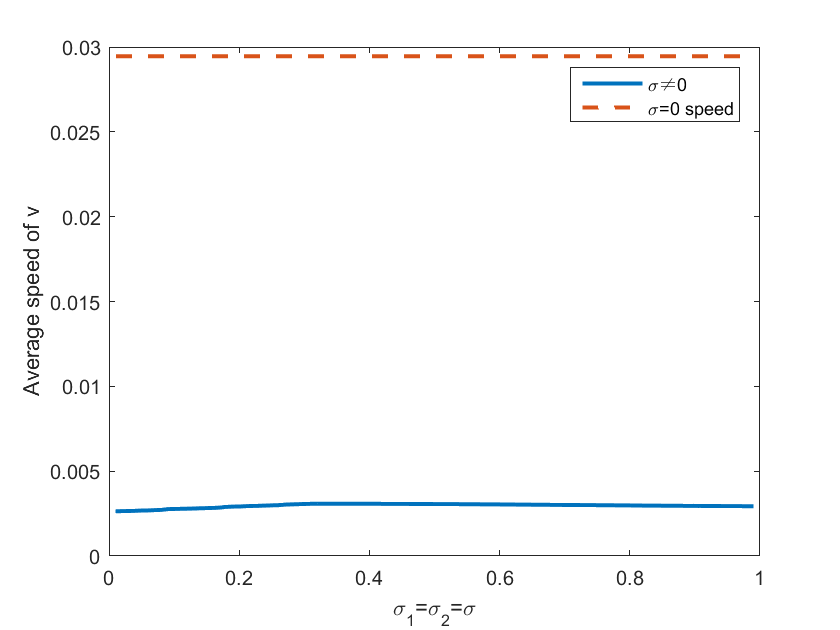}
} }
\subfigure[$c_1=0.3$]{ {
\includegraphics[width=0.35\textwidth]{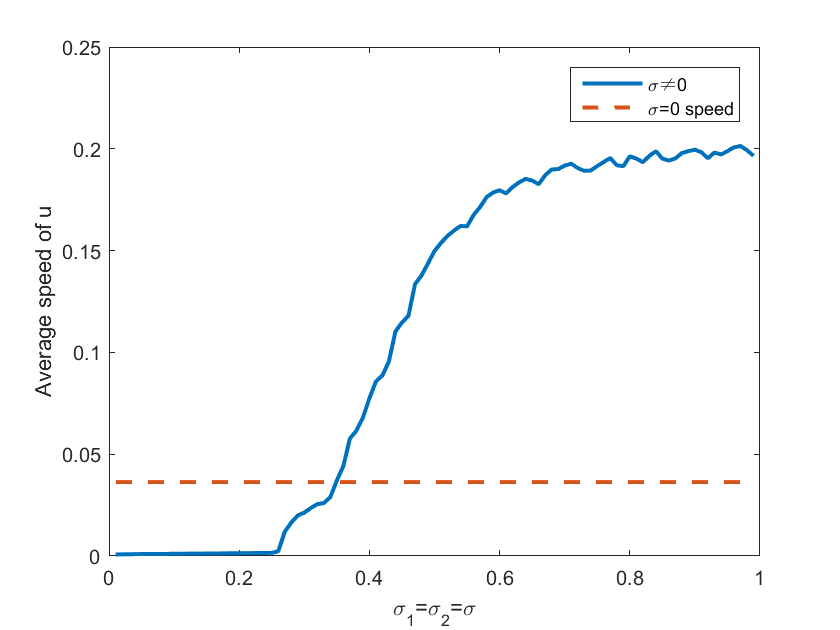}
} }
\subfigure[$c_1=0.3$]{ {
\includegraphics[width=0.35\textwidth]{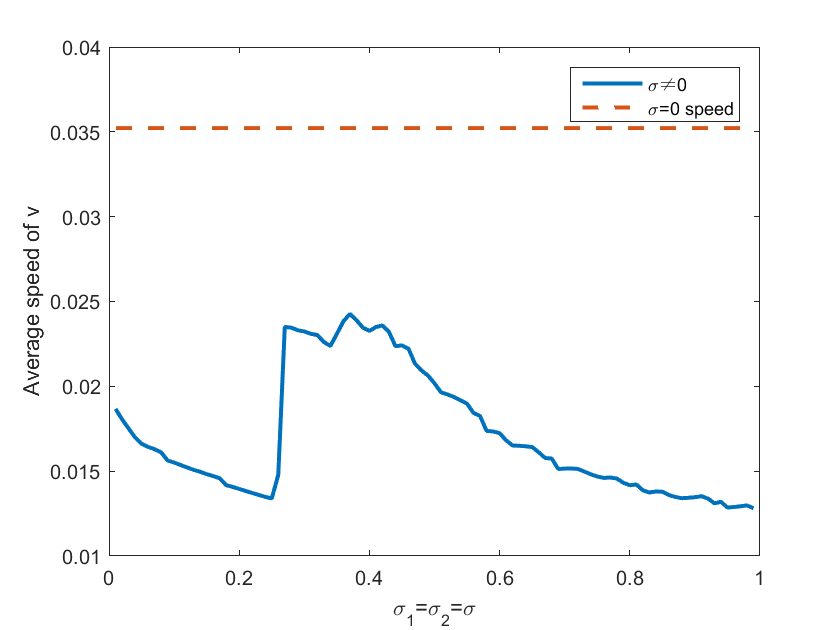}
} }
\caption{\scriptsize Average speeds of $u$ and $v$ with $\sigma_1=\sigma_2=\sigma\in[0.01,0.99]$ for impulsive effects $c_1=0.5$ in graphs $(a-b)$ and $c_1=0.3$ in graphs $(c-d)$. The spreading speeds of $u$  and $v$ in homogeneous environments are in red dashed lines. Graphs $(b)$ and $(d)$ show that the average speeds of $v$ with $\sigma\in[0.01,0.99]$ is smaller than spreading speeds in homogeneous environments, while the average speeds of $u$ with $\sigma\in[0.01,0.99]$ is bigger than spreading speeds in homogeneous environments, which is shown in graph $(a)$.}
\label{tus2}
\end{figure}

Note that the magnitudes of environmental perturbation $\sigma_1$ and $\sigma_2$ have nothing to do with threshold values $\lambda_1$ that  determine the competition outcomes. However,  they have impacts on spreading speeds of species. To investigate the effects of periodic environmental perturbations on the spreading of competitors, we will study how they affect spreading speeds of the competitors in the coexistence situations.
\begin{exm}{\small(Effects of the magnitudes $\sigma_1$ and $\sigma_2$ on spreading speeds)}
Fix  $c_2=1$ and $\tau=2$,  and  choose $\sigma_1=\sigma_2=\sigma$.

We choose  $t=40$ as a large time to illustrate  the spreading speeds of $u$ (i.e. $r'(40)$) and $v$ (i.e. $s'(40)$) with respect to $\sigma$ under different impulsive harvesting $c_1=0.5$ and $c_1=0.4$ in Fig. \ref{tus1}, where the red dashed line represents the spreading speeds of species in homogeneous environment, and the blue line stands for the spreading speeds by taking a series of magnitude $\sigma\in[0.01,0.99]$. The spreading speeds of $u$  and $v$ for $c_1=0.5$ are shown in Fig. \ref{tus1} (a-b), while corresponding speeds for $c_1=0.4$  are exhibited in Fig. \ref{tus1} (c-d). One can see from Fig. \ref{tus1} $(a)$ and $(c)$ that spreading speeds of superior $u$ affected by environmental perturbations are bigger than it in homogeneous environments, and from  Fig. \ref{tus1} $(b)$ and $(d)$ that there is not much difference in the spreading speeds of inferior $v$ with and without environmental perturbations.

We next simulate the average  speeds of $u$ (i.e. $\overline{r'(40)}$) and $v$ (i.e. $\overline{s'(40)}$) as $\sigma$ varies in Fig. \ref{tus2}, where $\overline{r'(t)}:=\frac{1}{\tau}\int_{t-\tau}^t r'(\rho)d\rho=\frac{r(t)-r(t-\tau)}{\tau}$ and $\overline{s'(t)}=\frac{s(t)-s(t-\tau)}{\tau}$ as in \cite{khan}. It is known from Fig. \ref{tus2} $(a)$ and $(b)$ that the average speed of $u$ is bigger than the speed in homogeneous environments, while average speed of species $v$ is smaller than the speed without environmental perturbations for impulsive effect $c_1=0.5$, which are similar as in \cite{khan}.   However, phenomena different from the results in \cite{khan} can occur. For example, under the case that $c_1=0.3$, the average speed of $u$ is not always bigger than the speed in homogeneous environments, see Fig. \ref{tus2} $(c)$.
\end{exm}

In Example 5.3, we can conclude from Fig. \ref{tus1} that periodic environmental perturbations speed up the spreading of the superior $u$ despite $u$ is harvested impulsively, and although the number of inferior $v$ is increased, the spreading speed of $v$ is not much affected by the environmental perturbations. Different from the results without pulses in \cite{khan}, the superior is more sensitive to environmental perturbations  due to the impulsive harvesting only on the superior. As it is shown in Fig. \ref{tus2}, the average speed of inferior $v$ is always decreased by the periodic environmental oscillations. And owing to the introduction of advection and pulses, it is different from \cite{khan} that the average speed of  superior $u$ is not always improved by periodic environmental perturbations.

\section{Discussions}
The heterogeneity of environments has received considerable attention and decreases the predictability of competition outcomes shown in several works about free boundary problems for competition models. Naturally, it is of interest to know whether impulsive intervention strategies can alter the competition outcomes, or drive the coexistence of two competitors in such heterogeneous environments.

The purpose of this paper is to understand how impulsive interventions and environmental perturbations  affect the spreading and vanishing of species and whether impulsive interventions can alter the competition outcomes.  Periodic environmental perturbations and  seasonal pulses are introduced into a free boundary problem for diffusive-advective competition model, which is more realistic for describing dynamics of two competing species and is also quite complex  mathematically.

The principal eigenvalue of impulsive periodic eigenvalue problem is first defined and the relationships respect to some parameters including diffusion, domain, pulse and advection are investigated. We show that competition outcomes for our problem can be classified into four different scenarios: co-extinction, coexistence, persistence of the species $u$ only, and persistence of the species $v$ only, see Theorems 3.1 and 3.7. Some sufficient conditions for species spreading and vanishing according to pulses are obtained, and a minimal pulse intensity for spreading of species is defined for a linear impulsive function. Meanwhile, inspired by the works in \cite{guowu2015,wu2015jde} for competition models with two different free boundaries, under the assumptions $(\mathcal{A}1)$ and $(\mathcal{A}2)$,  the minimal habitat for species spreading is introduced and the effects of expanding capability are studied as initial region belongs to different intervals. The asymptotic spreading speeds for the case that $\alpha,\beta\geq0$ are also given, which are related to advection rates and pulses.

Some interesting phenomena about the effects of pulses and environmental perturbations are last exhibited via numerical simulations. The efficacy and timing of impulsive interventions can significantly alter the competition outcomes (see Examples 5.1 and 5.2), which reveals that different impulsive strategies should be designed based on our pursued outcomes. Though  environmental perturbation is not involved in the threshold value $\lambda_1$, it plays an important role in  spreading speeds of species, which is mainly investigated in Example 5.3 and showed that the superior is more sensitive to environmental perturbations than the inferior when impulsive harvesting is carried out only on the superior.

We have overcome some difficulties caused by the combined effects of advection, pulses and heterogeneous environments. Mathematically, in order to obtain threshold values for spreading or vanishing, we introduce a  periodic eigenvalue problem with pulses and advection. Due to the dependence of  periodic eigenvalue problem on advection and pulse, we can not give the explicit expression of the principal eigenvalue. To study the long time behaviors of  our threshold values on diffusion, habitat and advection, we study a  fundamental eigenvalue problem depending only on $x$ instead by performing spatiotemporal variable separation on the periodic eigenvalue problem. Since advection rate belongs to $\mathbb{R}$, different types of suitable upper solution are constructed according to different advection rate. We have given some sufficient conditions on pulses for species spreading or vanishing. The proofs and results related to advection and pulses in this paper involves some modified  techniques, and  some new and rich perspectives.

It should be admitted that advection rates are also factors affecting the dynamics. Advection terms bring difficulties for not only analyzing the relationships between  principal eigenvalues and advection, but also on the long time behaviors of species.  We have not given the effects of advection on principal eigenvalues rigorously. And  competition outcomes have been considered only under the small advection rates (for example $|\alpha|<2\sqrt{D(\ln G'(0)/\tau+\gamma)}$) and negative advection rates (for example $\alpha\leq-2\sqrt{D(\ln G'(0)/\tau+\gamma)}$). It is known from \cite{ghjfa2015} that the dynamics of species for the large advection rate are much more complicated. We leave the large advection rate case for further study.

In this paper, we only gave some rough estimates of spreading speeds of competitors when they spread successfully under small positive advection rates. Two competitors possess two different moving regions and their boundaries may intersect. For our problem \eqref{a01} without advection and pulses, some sharper estimates of spreading of species, when they spread successfully, have been given in \cite{liuwang2019}, which  have  not be solved for our problem.
Also, noting that environmental perturbations and  seasonal pulses take the same time period in our model, we believe that for  general environmental perturbations, the time periods of environmental perturbations also have impacts on spreading speeds of species as in \cite{khan}. Meanwhile, we carry out impulsive intervention only on $u$ in our simulations. When different  impulsive interventions implemented in species $u$ and $v$, the phenomena may be different and  diverse under different parameters of $u$ and $v$.

\end{document}